\documentclass{article}
\usepackage{amssymb}
\usepackage[centertags]{amsmath}
\usepackage[dvips]{epsfig}

\def\vap{{\varphi}}
\def\vae{{\varepsilon}}

\def\dd{{\mathrm d}}

\def\Atens{\mbox{\boldmath $A$}}

\def\Atens{\mbox{\boldmath $A$}}

\def\acc{\mbox{\boldmath $a$}}

\newcommand{\de}{\partial}
\newcommand{\RR}{\mathbb{R}}

\newcommand{\ZZ}{\mathbb{Z}}

\numberwithin{equation}{section}

\begin{document}

\title{Felice Casorati's work on finite differences and its influence on Salvatore Pincherle}

\author{Riccardo Rosso \\
Dipartimento di Matematica ``F. Casorati" \\
Universit\`a degli Studi di Pavia \\
via Ferrata, 5, 27100 Pavia, Italy
}
\date{July, 9th, 2013}

\maketitle
\begin{abstract}
This paper, which is mainly based on unpublished material, focuses on the scientific influence that Felice Casorati exerted on Salvatore Pincherle. This influence
can be traced, in particular, in Casorati's work on the finite-difference calculus as conceived and published between 1879 and 1880 when Pincherle was living in
Pavia. Casorati's work has an interesting {\it back story} related to his entry to the 1880 {\it Grand Prix} of the French {\it Acad\'emie des Sciences} that
helps us in understanding Casorati's personality. Moreover, the correspondence that Casorati exchanged with other mathematicians on his work reveals that some of
the results contained in [Casorati, 1880b] had been obtained---though in a narrower context--in an early paper by Christoffel. Finally, the letters between
Casorati and Pincherle contain a short unpublished note by Pincherle on a paper by Jules Tannery [Tannery, 1875]. This note offers the first evidence of the
influence of [Casorati, 1880b] on Pincherle's work on the finite-difference calculus.
\end{abstract}

%\pacs{61.30.-v, 61.30.Pq (microconfined LC), 61.30.Jf (defects),
%64.60.Q- (nucleation), 64.60.qj 64.60.an (finite size)}

\maketitle

\section{Introduction}

Felice Casorati (1835-1890) and Salvatore Pincherle (1853-1936) are representatives of two successive generations of Italian mathematicians: Casorati, jointly
with Francesco Brio\-schi, Enrico Betti, Eugenio Beltrami, and Ulisse Dini, played an important r\^ole in the renaissance of Italian mathematics that intertwined
with the political events leading to unification of Italy. An important step in this process was  Italy's entry into the main stream of European mathematical
research as cultivated in Paris, G\"ottingen, and Berlin. In this sense, the journey through Europe that Casorati, jointly with Betti and Brioschi, performed in
1858 cannot be underestimated, since first-hand information on the main trends in mathematics were gained and then diffused in Italy [Volterra, 1902]. Casorati's
relevant efforts culminated with the publication of the volume  {\it Teorica delle funzioni di variabili com\-ples\-se} [Casorati 1868] where he expounded the
main results in complex analysis, with particular attention to Riemann's approach ({\it see} [Bottazzini, 1994], Chapts. IV and VI).

Pincherle, jointly with Luigi Bianchi, Ernesto Ces\`aro, Giuseppe Peano, Vito Vol\-ter\-ra, and Tullio Levi-Civita, belonged to the subsequent generation of
mathematicians who contributed to the strengthening of Italy's position within the international scientific community. His main contributions were in the field of
functional (linear) calculus where he worked for more than thirty years along a route distinct from the one that Vito Volterra followed in the same period.

Pincherle graduated from Pisa as a student of the {\it Scuola Normale Superiore}, and spent his academic life entirely in Bologna. His first scientific interest
was in the theory of capillarity which was the subject of his thesis ({\it Tesi di Laurea}) published in the {\it Nuovo Cimento} [Pincherle, 1874] with a short
addition [Pincherle, 1875].  In 1875, when he moved to Pavia to teach at the {\it Liceo Ginnasio} ``Ugo Foscolo", he was recommended to Felice Casorati by Ulisse
Dini:\footnote{\begin{flushleft}
Caro Casorati
\end{flushleft}

\begin{flushright}
5. 9bre 75
\end{flushright}

Occupato come sono, non posso in questo momento scriverti a lungo. Ti scrivo soltanto due parole per presentarti il D$^{\rm r}$. Salvatore Pincherle, uno dei
migliori allievi della nostra Universit\`a e della nostra Scuola, attualmente Professore al Liceo di Pavia. Troverai in Lui un giovane studioso e volenteroso di
fare e con molto ingegno. Ajutalo e consiglialo nei suoi studi.

\begin{flushleft}
Credimi
\end{flushleft}

\begin{flushright}
Tuo aff. amico,\\
U. Dini
\end{flushright}}

\medskip

\begin{flushleft}
{\it Dear Casorati}
\end{flushleft}

\begin{flushright}
{\it November 5th, 75}
\end{flushright}

{\it Busy as I am, I cannot write at length to you. Just a couple of words to introduce you D$^{\rm r}$. Salvatore Pincherle, one of the best students of both our
University and our Scuola, at present Professor at the Liceo in Pavia. You will find in him a studious young man, willing and very ingenious. Please, help him and
advise him in his studies.}

\begin{flushleft}
{\it Yours sincerely}
\end{flushleft}

\begin{flushright}
{\it U. Dini}
\end{flushright}

\medskip

This short letter opened the way to Casorati's and Pincherle's personal relationship that lasted for 15 years, up to the untimely death of Casorati in 1890. Their
correspondence amounts to some 80 letters that have been kept in  Casorati's archive.\footnote{Hereafter denoted as his {\it Nachlass}, after [Neuenschwander,
1978], it is preserved by Casorati's descendants in Pavia. It has been made available to me by Prof. Alberto Gabba, to whom I address my warmest thanks.} These
letters provide some further information about Casorati's influence on Pincherle. At the very beginning of his stay in Pavia, Pincherle had devoted his interests
to the theory of minimal surfaces that formed the subject of three notes [Pincherle, 1876a,b,c] published in 1876, but the interaction with Casorati turned his
attention to differential equations. Shortly thereafter Pincherle obtained a fellowship to spend a year in Berlin, where he met Weierstrass, whose methods were to
exert a decisive influence on his scientific activity ([Bottazzini, 1994], Chapt. VI). Direct contacts with Casorati were resumed in 1878 when Pincherle came back
to Pavia from Berlin. At the beginning of 1879 Pincherle applied for the habilitation ({\it libera docenza}) in higher calculus at the University of Pavia. A
local committee examined his application and, the judgement being positive, on May 9th the Faculty approved his habilitation, which was ratified by the Ministry
of Education at the end of 1879.\footnote{The original documents concerning Pincherle's habilitation are kept in the Historical Archive of the University of
Pavia.} Consequently, in the fall semester 1879--80 Pincherle gave a course in the higher calculus ({\it Analisi superiore}) that   Beltrami and Casorati among
others attended,\footnote{See [Bortolotti, 1937].} where he offered an exposition of Weierstrass's views on the theory of analytic functions. Pincherle's lectures
and, more importantly, his long account of them published in Battaglini's {\it Giornale di Matematiche} [Pincherle, 1880] contributed significantly to the spread
of Weierstrass's ideas in Italy\footnote{See [Bottazzini, 1994, Chapt. VI].}. In the period 1878--1880 Pincherle also attended  Casorati's lectures on higher
calculus, as he reminded in the following letter addressed to Casorati in 1888:

{\begin{flushright}
Bologna, li 25/3/88
\end{flushright}

{\it Dear Professor

\bigskip

During the researches on functional operations represented by definite integrals and on their inversion, on which I have been working for a while, it occurred to
me to study linear differential equations, whose solution is a problem that arises as a particular case in the inversion of definite integrals; on this occasion I
rejoiced at having attended the lectures you gave on this subject, sparing to your audience the reading, and the not-too-easy {\rm digestion}, of the Memoirs by
Fuchs and Thom\'e. (...)}\footnote{Chiarissimo Sig. Professore

\bigskip

Nel corso delle ricerche che ho intraprese gi\`a da parecchio tempo sulle o\-pe\-ra\-zio\-ni funzionali rappresentate da integrali definiti e sulla loro
inversione, mi \`e capitato di dover studiare le equazioni differenziali lineari, la cui risoluzione \`e un problema compreso come caso particolare
nell'inversione degli integrali definiti; ed in questa occasione ho avuto pi\`u di una volta da rallegrarmi di avere assistito alle lezioni che Ella ha fatto in
proposito, risparmiando ai Suoi uditori la lettura e la {\it digestione} non troppo agevole delle Memorie di Fuchs e Thom\'e.}

\medskip

Actually, since the mid-1870s Casorati had turned his attention to Fuchs's theory of linear differential equations in the complex plane that was consonant with
his early concern with complex analysis\footnote{For a historical account of the r\^ole of complex function theory on the research in differential equation, we
refer the reader to Chapt. 7 of [Bottazzini and Gray, 2013].}, and in the academic years 1878--79 and 1879--80 Casorati based his course in higher calculus on
this topic, too. As a comment to the opening lecture, held on November 16th 1879, Casorati wrote in his log-book:

\bigskip

{\it I declare that I will carry out the theory of linear differential equations according to Fuchs, etc. as the first topic. They [the students] will find the
original memoirs in Borchardt's Journal.}\footnote{Dichiaro che svolger\`o per primo argomento la teoria delle equazioni differenziali lineari secondo Fuchs, ecc.
Troveranno le memorie originali nel Giornale di Borchardt.} [Casorati, 1879]

\bigskip

 Casorati's influence on Pincherle had three aspects. Firstly, there was what could be called a {\it political} influence: Casorati knew personally most of the
 leading mathematicians in Europe and, together with Betti, introduced Pincherle to Kronecker and Weierstrass when he got a fellowship to spend the academic year
 1877--1878 in Berlin. Moreover, Casorati was a member of both the Committee that declared Pincherle eligible for the chair of Infinitesimal Analysis in 1880 at
 Bologna {\it ex aequo} with Cesare Arzel\`a, who was preferred to Pincherle on that occasion, and the Committee that promoted Pincherle to a full professorship
 in 1888. In addition, Casorati supported Pincherle's election as a corresponding member ({\it Socio Corrispondente}) of the Accademia dei Lincei in 1887.
 Secondly, there was an {\it editorial} influence. From 1868 Casorati was an effective member of the Royal {\it Istituto Lombardo di Lettere, Scienze ed Arti}
 based in Milan, and as such he presented some early papers by Pincherle to this local Academy for publication. More importantly, Casorati was a member of the
 editorial board of the {\it Annali di Matematica pura ed applicata} that had Francesco Brioschi as the editor in chief. In this capacity Casorati  presented
 Pincherle's first important memoir [Pincherle, 1884a] on series expansions  for publication in the {\it Annali}. Thirdly, there was a {\it scientific} influence
 that stemmed from Casorati's 1880 paper on the interpretation of finite difference [Casorati, 1880b], to which I will devote more attention in the present
 paper.

An early acknowledgement of the close connection between Casorati's and Pincherle's researches on that topic was given in the foreword of [Guldberg and
Wallenberg, 1911] where the authors stated (p. vi)

\bigskip

{\it In the meanwhile, another area of our discipline underwent a powerful development, precisely the formal side of the theory of linear difference equations,
due in particular to the analogy with algebraic equations, and to the extension of the analogies with linear differential equations. Here, after Casorati,
Pincherle with his students Bortolotti\footnote{Ettore Bortolotti, better known as a historian of mathematics, in his early days devoted several papers to the
formal theory of finite differences, ({\it see} for instance [Bortolotti, 1895]).} and Amaldi---above all---has built the foundations and, in particular, in
Chapt. 10 of his book on distributive operations\footnote{We will discuss below in detail the content of Chapt. X of the monograph [Pin\-cher\-le and Amaldi,
1901].}, he has built the frame, so to say, for the edification of such a theory.}\footnote{Inzwischen erfuhr aber ein anderes Gebiet unserer Disziplin eine
m\"achtige F\"orderung, n\"amlich die {\it formale} Seite der Theorie der linearen Differenzengleichungen, insbesondere soweit sie sich auf die Analogien mit den
algebraischen Gleichungen und auf die entsprechenden Analogien mit den linearen Differentialgleichungen erstreckt. Hier hat nach {\it Casorati} haupts\"achlich
{\it Pincherle} mit seinen Sch\"ulern {\it Bortolotti} und {\it Amaldi} der Grundstein gelegt und besonders im 10. Kapitel seiens Buches \"uber die distributiven
Operationen gewisserma{\ss}en das Ger\"ust f\"ur den Aufbau einer solcher Theorie errichtet.}

\bigskip

In his  [1880b] Casorati had proved the necessary and sufficient condition for $n$ functions of a variable that undergoes discrete increments to be linearly
independent by introducing a suitably defined determinant, now known as the Casorati determinant, or the Casoratian. However, the aspect that Casorati emphasized
there was the wide spectrum of applications of the classical finite-difference calculus by showing that, after suitable interpretations, classical results could
be translated into complex function theory: in particular, Casorati noticed that the local behaviour of the solution of a linear differential equation around a
singular point [Fuchs, 1866] could be easily obtained by adapting finite-difference calculus appropriately.

Casorati attached great importance to the results contained in his paper. Also urged by Brioschi, with great expectation he decided to apply for the 1880 Grand
Prix de Math\'ematiques announced by the French Academy of Sciences, and devoted to the theory of differential equations. However, his paper found in Halphen's
memoir [Halphen, 1884] an invincible obstacle to winning the prize, and Casorati had to content himself with a special mention shared with the young Henri
Poincar\'e. Casorati's [1880b] paper also lies at the heart of a harsh argument that Casorati had with Ludwig Stickelberger, who was working somehow on the same
topics.  The paper [Stickelberger, 1881], published a few months after [Casorati, 1880b], appeared to Casorati as an attempt to question, if not to deny, his
priority in the applications of finite differences to complex analysis. Their controversy reached its {\it climax} with the publication of a short letter in the
{\it Annali} [Casorati, 1880c] where Casorati staunchly defended his primacy in bringing new life to finite difference calculus.

In Casorati's {\it Nachlass} there is a folder where Casorati kept most of his correspondence with several mathematicians concerning both [Casorati, 1880b] and
the Grand Prix. The content of those letters provides us with a vivid account of his personality, and informs us of the unpredictable events that contributed to
hamper his work. Moreover, copies of letters that Stickelberger and Casorati sent to Brioschi after the publication of [Casorati, 1880c] also contain some
valuable information, not only because they offer an idea of Stickelberger's point of view, but also because there Stickelberger pointed out that Casorati's
criterion of linear dependence among functions of a discrete variable had already been established by Christoffel in [Christoffel, 1858].

In Italy, the main results of [Casorati, 1880b] were made known to a broader audience through a detailed account that appeared in the {\it Giornale di
Ma\-te\-ma\-ti\-che} by Paolo Cazzaniga, a relative of Casorati's, who presented them as a chapter of the symbolic calculus [Cazzaniga, 1882]. As we shall see,
Casorati had denoted the operation mapping $f(x)$ into $f(x+1)$ by $\theta$, and he had emphasized the importance of conceiving it as an operator mapping {\it
functions} into {\it functions}. This is the aspect that most influenced Pincherle's subsequent treatment of the $\theta$ operator. He first applied Casorati's
operator to linear differential equations in a short note written in 1885. Even though this note was meant to appear in the {\it Rendiconti} of the Istituto
Lombardo, it was left unpublished. There, resorting to Casorati's [1880b] approach Pincherle re-obtained a theorem proved in [Tannery, 1875] by which Tannery had
given a characterization of the integrals of a linear differential equation  with uniform coefficients in the complex plane: a theorem that is the ``pleasant
converse to Fuchs's main theorem" ([Gray, 1986], p. 90). Pincherle explained his result to Casorati who, though appreciating the simplicity of the latter's
result, suggested him to add some more applications. Actually, he asked Pincherle to correct a point in Tannery's paper but, as Pincherle could not find any
relevant mistake at the point indicated by Casorati, he did not insist in the publication, all the more so because he was diverted from this issue by the
publication of a paper [Mellin, 1886] where Mellin had bridged a differential and a difference equation by employing a transformation that was a particular case
of a more general class just obtained by Pincherle: this urged the latter to write  a note [Pincherle, 1886a] to secure his priority on this topic.

A few years later, in the academic year 1893--94 Pincherle went back to Casorati's $\theta$ operator that in his hands became the prototype of linear operators
acting upon sets of analytic functions, i.e. the functions that Pincherle considered in his theory of linear operators where the symbolic calculus rose to the
rank of a rigorous mathematical theory, and ceased to be only a useful but dubious compact notation.

The present paper is organized as follows. In Section 2 we describe Casorati scientific activity in differential equations in the period 1875--1877 when he first
interacted with Pincherle and then, in Section 3, we illustrate the content of [Casorati 1880b]. Section 4 is entirely devoted to the history of Casorati's
application to the Grand Prix des Math\'ematiques by examining the letters he exchanged with both Italian and foreign mathematicians (the original letters are
transcribed in Appendix 1, see Section \ref{sec:originals_GP}). Intertwined with this topic there is the argument with Stickelberger that will be described mainly
in its lesser known aspects through the letters that both he and Casorati addressed to Brioschi after the publication of [Casorati, 1880c]: some excerpts of these
letters have been also transcribed in Appendix 1 (see Section  \ref{sec:originals_GP}). Then, in Section 5 Pincherle will be back on stage as we  examine the
letters he exchanged with Casorati concerning his note devoted to Tannery's theorem: in particular, Pincherle's unpublished note will be translated into English,
while the original, relevant letters have been transcribed in Appendix 2 (see Section \ref{sec:originals_Tannery}). Finally, in Section 6 we will examine the
impact that [Casorati, 1880b] had on Pincherle's research in the last decade of the XIX century up to the publication of his treatise [Pincherle and Amaldi,
1901]. In a conclusive section, we will resume the outcomes of the analysis contained in this paper, and point to directions for future work.

\section{Casorati and Pincherle interactions on differential equations}\label{sec:eqdiff}

During the period that Pincherle spent in Pavia, which  spans from 1875 to 1880 with the exclusion of his Berlin year, Casorati's main scientific activity was in
the field of ordinary differential equations. In fact, Casorati had resumed scientific activity in analysis in 1874 after a long break. As he wrote in a  letter
to Gaston Darboux on August 3rd, 1879

\bigskip

{\it ...For a lot of different reasons, I seldom read [the Darboux {\it Bulletin}] in the years 1872, 73, 74. In those years, as in the previous years 70, 71, I
lived as a stranger to mathematics.}\footnote{..., per un cumulo di cause diverse, avvenne che negli anni 1872, 73, 74 non lo osservassi quasi mai. In quegli
anni, come nei precedenti 70 e 71, ho vissuto da straniero alla matematica.}

\bigskip

In the period there alluded to, Casorati was intensely involved in teaching both in Pavia and in Milan, at Brioschi's {\it Regio Istituto Tecnico Superiore} --
today the Milan Polytechnic -- but in his personal life Casorati also suffered several grievous events that impaired him  to the point of hampering his research.
He  lost his beloved father -- the renowned physician Francesco Casorati -- in 1859; then, his mother, and his life-mate. The following passage from a 1880 letter
to Joseph Bertrand (February 6th) clarifies Casorati's feelings:

\bigskip

{\it After that a cruel illness mercilessly led my father to death---he, who was the sun of my life, the man to whose goodness  I will never be able to think
without being completely affected---my heart and my health were no longer as before; that inexplicable grief and that cruel reality depressed my forces,
annihilated the poetry of life, in short, they destroyed my person both physically and morally. Moreover, I had my mother suffering from cancer for six
consecutive years, when our family was reduced to the two of us and, finally, she died; it will not be hard to believe for you that, for several years, I remained
as suspended between life and death. Finally, I improved; but my soul has remained, and will always remain, ready to be frightened, to get disheartened, to
overstate any misfortune, to suffer enormously.}\footnote{Dopoch\'e un morbo crudele trasse spietatamente alla tomba mio padre, il sole della mia vita, l'uomo
alla cui bont\`a non potr\`o mai pensare senza commuovermi tutto, l'animo mio e la mia salute non furono pi\`u quelli di prima; quell'indicibile dolore, quella
feroce realt\`a abbatterono le mie forze, annientarono la poesia della vita, demolirono insomma tutta la mia persona fisica e morale. Avuta, per giunta, mia madre
ammalata di cancro per sei anni consecutivi, mentre la famiglia riducevasi a noi due soli, ed essa infine soccombeva; non stenterete a credere come io dovessi per
lunga serie di anni rimanerne come tra questo e l'altro mondo. Finalmente migliorai; ma l'anima mia \`e rimasta e rimarr\`a sempre pronta a sgomentarsi, ad
avvilirsi, ad esagerarsi ogni contrariet\`a, a soffrire senza misura.}

\bigskip

Looking at Casorati's scientific production on differential equations, we can single out two main periods, each of them with its own {\it leitmotif}. In the first
period up to 1878, Casorati's main concern was the research on {\it singular} solutions of differential equations, a topic on which mathematicians such as Cayley
and Darboux also devoted their attention. Since 1878, Casorati was mostly influenced by Fuchs's theory of linear differential equations in the complex plane.
Since Casorati's activity in the former period had little impact on Pincherle, I will limit myself to a few remarks.

Following the definition given by  Boole in his treatise [Boole, 1859], a textbook on which Casorati often based part of his lectures on infinitesimal calculus,

\bigskip

{\it a singular solution of a differential equation is a relation between $x$ and $y$, which satisfies the differential equation ``by means of the value it gives
to the differential coefficients" $\frac{\dd y}{\dd x}$, $\frac{\dd^2 y}{\dd x^2}$, $\&$c., but is not included in the complete primitive.} ([Boole, 1859], p.
139)

\bigskip

\noindent As an example, Boole considered the equation ([Boole, 1859], p. 150)
\begin{equation}\label{Boole}
y^2-2xy\frac{\dd y}{\dd x}+(1+x^2)\left(\frac{\dd y}{\dd x}\right)^2=1
\end{equation}
whose general solution (or complete primitive) is given by
\begin{equation}\label{Boole1}
y=cx+\sqrt{1-c^2},
\end{equation}
where $c$ is an arbitrary constant. However it is also solved by
\begin{equation}\label{Boole2}
y=\sqrt{x^2+1}
\end{equation}
which is a singular solution of \eqref{Boole} since there is no value of $c$ for which \eqref{Boole1} yields \eqref{Boole2}.

Casorati focused his attention on algebraic differential equations whose coefficients are polynomials in both the dependent and the independent variables:

\medskip

A first order differential equation between two variables $u$ and $v$ is called algebraic-differential when it yields a constraint among $u$, $v$, $\dd u$ and
$\dd v$ that can be expressed by a {\it finite} number of {\it algebraic} operations.\footnote{Una equazione differenziale di primo ordine tra due variabili $u$ e
$v$ va detta {\it{algebrico-differenziale}} quando consista in un legame tra $u$, $v$, $\dd u$ and $\dd v$ esprimibile mediante un numero {\it finito} di
operazioni {\it algebriche}.} ([Casorati, 1875], [Casorati, 1952, p. 10])

\medskip

\noindent If $\vap_k(u,v)$ are $m+1$ polynomials, then it is always possible to reduce such an equation to the form
$$
\sum_{k=0}^m \left(\begin{array}{c} m \\ k\end{array}\right)\vap_{m-k}(u,v)(\dd u)^{k}(\dd v)^{m-k}=0
$$
where $\left(\begin{array}{c} m \\ k\end{array}\right)$ are binomial coefficients and $m$ determines the {\it degree} of the equation.

Dissatisfied with the procedures he found in the relevant  literature, Casorati proposed a new theory for dealing with singular solutions of algebraic
differential equations of first order and second degree [Casorati, 1876]. Specifically, the equations considered by Casorati are
\begin{equation}\label{algdiff2}
   \alpha(u,v)(\dd u)^2+2\beta(u,v)\dd u\dd v+\gamma(u,v)\dd v^2=0
\end{equation}
where the polynomials $\alpha$, $\beta$ and $\gamma$ have no common factors. Casorati made the hypothesis that \eqref{algdiff2} has an algebraic complete
primitive.

For an equation of the first order and the first degree, that is,
\begin{equation}\label{Rosanes}
    A(x,y)\dd x+B(x,y)\dd y=0
\end{equation}
 Jakob Rosanes had worked out a condition under which this is the case, in [Rosanes, 1871]. It is precisely to equations like \eqref{Rosanes} that Pincherle
 turned his attention before his leaving for Berlin. He explained his reflections on algebraic differential equations to Casorati in a meeting ({\it Conferenza})
 on February 5th 1877, a few weeks before presenting the note [Pincherle, 1877] at the Istituto Lombardo. Like Rosanes, Pincherle emphasized the geometric meaning
 of \eqref{Rosanes} that can be generated by taking a one-parameter family of algebraic curves
\begin{equation}\label{fascio}
G(x,y)+\lambda H(x,y)=0
\end{equation}
where $G$ and $H$ are irreducible polynomials of degree $n$, and $\lambda$ is a parameter.
By using the differential of this equation to get rid of $\lambda$, an equation of the class \eqref{Rosanes} is obtained. Pincherle resorted to the geometric
meaning of \eqref{Rosanes} to obtain {\it analytically} a result  geometrically established by Cremona fifteen years before:

\bigskip

{\it In a family of curves of order $m$, there are $3(m-1)^2$ double points.} ([Cremona, 1862], p.67)

\bigskip

In the second part of [Pincherle, 1877] Pincherle followed the opposite direction: under the hypothesis that \eqref{fascio} is an integral of \eqref{Rosanes}, he
aimed at finding its explicit expression in terms of $A$ and $B$. The results obtained in this part of the paper are subject to restrictions --- as Pincherle
himself stated in his 1887 curriculum vitae submitted  for the election at the {\it Lincei}. Actually, these restrictions were pointed out to Pincherle by
Casorati on September 1878, that is, almost eighteen months after the presentation of the note. These limitations concern some terms that were omitted by
Pincherle, but taken into account by Casorati in his research on this subject. All this might be surprising, as Casorati had himself read the note to the members
of the {\it Istituto Lombardo} but, presumably, on that occasion he did not take too much care over the details of Pincherle's note since he relied on the meeting
he had had with the latter a few weeks before. On the other hand in 1878, when Casorati was working on equations like \eqref{Rosanes} [Casorati, 1878b], he  read
Pincherle's paper more accurately, and pointed him out that the generality of his conclusions had to be considerably reduced\footnote{Casorati met Rosanes in
Parpan (Switzerland) on August 27, 1879, and on that occasion they discussed of research on the equation $A\dd x+B\dd y=0$.}.

\bigskip

%{\it On (August) 27th, Rosanes arrived in Parpan. I asked him if he had made further research on the algebraic integration of $A\dd x+B\dd y=0$. He answered,
%substantially, in the negative. When I told him that, apart from his and Pin\-cher\-le's paper and a few more lines on this subject, I did not know papers where
%progress on the subject had been made, he himself was unable to quote any paper. He did not know my Note on the integration by linear functions and I got the
%impression that he had not heard about Darboux's similar papers. He is in Breslau with Schr\"oter.}\footnote{Il 27 (Agosto) venne a Parpan Rosanes. Gli domandai
%se avesse fatto ulteriori ricerche sulla integ.[razione] alg.[ebrica] di $A\dd x+B\dd y=0$. Risposemi, in sostanza, di no. Ed avendogli poi detto, che, fuori del
%suo e del lavoro del Pincherle e di altre righe relative all'argomento, non conosco scritti in cui siasi fatto progredire l'argomento; anch'egli non me ne seppe
%citare. Non aveva cognizione della mia Nota sulla integ.[razione] per funzioni lineari; ne parvemi avesse sentore degli affini articoli di Darboux. Egli \`e a
%Breslavia presso Schr\"oter.}

\section{Casorati's paper on finite differences}\label{sec:Casorati80}

In the history of mathematics, Casorati is best remembered for the Casorati-Weierstrass-(Sokhotskii) theorem he published in [Casorati, 1868]. However, his name
has also been linked to a special determinant --- the Casoratian or the Casorati determinant --- that plays a  r\^ole in linear difference equations analogue to
the one played by the Wronskian in the theory of linear differential equations. This determinant was introduced in [Casorati 1880b], the paper by Casorati that
mostly influenced Pincherle. This long memoir had a history that I review in the next section,  after having examined its content.

Casorati applied the formalism of finite differences to study the behaviour of a function $y$ of a complex variable $x$ in the neighborhood of a point $x_1$. By
taking a point $x$ in a circular annulus surrounding $x_1$ within which $y$ is free of singularities, and by performing a simple closed contour starting at $x$,
Casorati set
$$
\Delta y(x):=y(x+{\rm e}^{2\pi i})-y(x)
$$
and he noted that
 $$
t(x):=\frac{\log (x-x_1)}{2\pi i}
$$
is such that $\Delta t=1$. Moreover single-valued (monodromic) functions are such that $\Delta y=0$ and so they play the same r\^ole as that played in standard
calculus of differences by periodic functions $f$ whose period is equal to the constant increment $h=1$ of the independent variable. Casorati introduced an
operator --- the $\theta$ operator --- defined by
\begin{equation}\label{theta}
\theta y(x):=\Delta y(x)+y(x)=(\Delta+1)y(x)
\end{equation}
which is clearly the counterpart of the operator mapping $f(t)$ into $f(t+1)$. In this new interpretation $\theta y$ yields the value of $y$ after that the closed
contour has been completed once. Casorati remarked:

\bigskip

{\it It is well known that, in the Calculus of differences, together with the differences $\Delta y$, $\Delta^2 y$ of a function, it is also important to consider
the values that it attains at the values $t+1$, $t+2$,.... of the variable. These values are mostly denoted as follows:
$$
y_{t+1},\qquad t_{t+2},...
$$
But there is a great advantage in denoting these values with another notation, such as the following
$$
\theta y,\qquad\theta^2 y
$$
that makes it possible to consider $\theta$ as a symbol of an operation, that is, of that operation which  yields $y_{t+1}$ when it is performed on
$y_t$.}\footnote{Si sa che nel {\it Calcolo delle differenze} importa di considerare insieme con le differenze $\Delta y$, $\Delta^2 y$ di una funzione, anche
separatamente i valori che essa prende corrispondentemente ai valori $t+1$, $t+2$,.... della variabile. Questi valori vengono per lo pi\`u significati come
segue:
$$
y_{t+1},\qquad t_{t+2},...
$$
Ma vi ha grande vantaggio a significare questi valori con altra notazione, come la seguente
$$
\theta y,\qquad\theta^2 y
$$
la quale si presta a lasciar riguardare $\theta$ come un simbolo di operazione, cio\`e di quell'operazione che eseguita su $y_t$ d\`a per risultato $y_{t+1}$.}
([Casorati, 1880b], p.13; [Casorati 1951], p. 320)

\bigskip

Such an {\it operator} interpretation is in the spirit of Pincherle's subsequent views as  clearly stated, for instance, in his major memoir on Laplace transform
[Pincherle 1887] where he considered the function
$$
f(x):= \int{\rm e}^{xy}\vap(y)\dd y,
$$
the integration being performed along a line in the complex plane. Pincherle observed that $f(x)$ is an analytic function if $\vap(y)$ is analytic as well, and
hence

\medskip

\noindent{\it in this paper we will consider this transformation under a new point of view: that is, as a functional operation (...) defined by some
characteristic properties.}\footnote{In questo lavoro si vuole considerare questa trasformazione sotto un punto di vista nuovo: si vuole cio\`e riguardarla come
una operazione funzionale (...) definita da alcune sue propriet\`a caratteristiche.} ([Pincherle 1887, p. 125])

\medskip

\noindent The most remarkable property of the operator $\theta$ is that it is {\it eminently distributive} ({\it eminentemente distributivo}), ([Casorati 1880b],
p. 13; [Casorati 1951], p. 321) since it is not only linear but it also obeys
\begin{equation}\label{distributive1}
\theta(u(x)v(x))=\theta(u(x))\theta(v(x))\qquad\theta\left(\frac{u(x)}{v(x)}\right)=\frac{\theta u(x)}{\theta v(x)}
\end{equation}
and
\begin{equation}\label{distributive2}
\theta F(u(x), v(x), w(x),...)=F(\theta u(x), \theta v(x), \theta w(x),....),
\end{equation}
for any operator $F$ that yields a unique result when applied on functions $u(x)$, $v(x)$, $w(x)$,... . The main result contained in [Casorati 1880b] is the
following theorem, proved in Chapt. 2:

{\bf Theorem} {\it Among the functions $y_1(x)$, $y_2(x)$,...$y_n(x)$ of a complex variable $x$ there exists or not a linear, homogeneous relation with
coefficients that are single-valued (monotropi) in an annulus centered at $x=x_1$ according to whether the determinant}
 $$
 \Theta=\left|\begin{array}{cccc}
 y_1 & y_2 & \dots & y_n \\
 \theta y_1 & \theta y_2 & \dots & \theta y_n \\
 \dots & \dots & \dots & \dots \\
 \theta^{n-1} y_1 & \theta^{n-1} y_2 & \dots &\theta^{n-1} y_n
 \end{array}\right|
 $$
{\it vanishes or not.} ([Casorati 1880b], p. 19; [Casorati 1951, p. 329])

The direct theorem is easy to prove since, assuming that
$$
\sum_{i=1}^n \vap_i y_i=0
$$
and observing that $\theta\vap_i=0$ as these functions are single-valued, then $\theta(\vap_i y_i)=\theta\vap_i\theta y_i=0$. Thus, also
$$
\sum_{i=1}^n \vap_i\theta^k y_i=0\qquad \forall k=0,...,n-1
$$
and so the determinant $\Theta$ vanishes. Casorati proved the converse by induction on the number of functions $y_i$.
When $n=2$, from
$$
\det\left( \begin{array}{cc} y_1 & y_2 \\ \theta y_1 & \theta y_2 \end{array}\right)=0
$$
it follows that
$$
\frac{\theta y_1}{\theta y_2}=\frac{y_1}{y_2}
$$
and since $\theta\left(\frac{y_1}{y_2}\right)=\frac{\theta y_1}{\theta y_2}$, he concluded that the ratio $\frac{y_1}{y_2}$ goes back to its original value, when
its argument has performed a complete tour. Hence it is a single-valued function, say $-\frac{\vap_2}{\vap_1}$, and the relation
$$
\vap_1y_1+\vap_2 y_2=0
$$
holds. Induction from $n-1$ to $n$ is then performed by resorting to a transformation of determinants given in [Hermite, 1849]. We note that, as a result of the
new interpretation of finite difference formalism, the linear, homogeneous relation between the functions $y_i(x)$ has {\it variable} coefficients, contrary to
its counterpart based upon the Wronskian.

As a first application in Chapt. 3 Casorati employed his formalism to study the structure of the roots of the algebraic equation
$$
z^n+\psi_1 z^{n-1}+\dots +\psi_{n-1}z+\psi_n=0,
$$
where $\psi_i$ are single-valued functions of a complex variable $x$. In doing this, Casorati simplified considerably the analysis contained in Puiseux's long
memoir [Puiseux, 1850]. In Chapt. 4 Casorati applied the theorem on the $\Theta$ determinant --- that is the Casoratian --- to the linear differential equation
with single-valued coefficients $p_i(x)$ in the complex plane
\begin{equation}\label{fuchs}
D^my+p_1 D^{m-1}y+\dots +p_{m-1}Dy +p_m y=0,
\end{equation}
where $D:=\frac{\dd}{\dd x}$. Given a point $x_1$, Casorati associated to the equation \eqref{fuchs} a homogeneous difference equation that is linear and of the
$m$-th order as well but with {\it constant} coefficients, whose general solution he had recalled in Chapt. 6. Actually, Casorati observed that, when $y$ denotes
the general integral of \eqref{fuchs}, then \eqref{fuchs} itself can be viewed as a linear, homogeneous relation with single-valued coefficients $p_i$ among the
functions $y, Dy,....,D^m y$, and so the corresponding $\Theta$ determinant vanishes, that is
$$
H=\left|\begin{array}{cccc}
 y & Dy & \dots & D^my \\
 \theta y & \theta Dy & \dots & \theta D^my \\
 \dots & \dots & \dots & \dots \\
 \theta^{m} y & \theta^{m} Dy & \dots &\theta^{m} D^m y
 \end{array}\right|=0.
 $$
However, since $\theta$ commutes with $D$,
$$
H=\left|\begin{array}{cccc}
 y & Dy & \dots & D^my \\
 \theta y & D\theta y & \dots & D^m\theta y \\
 \dots & \dots & \dots & \dots \\
 \theta^{m} y & D\theta^{m} y & \dots &D^m\theta^{m} y
 \end{array}\right|=0
 $$
follows, that can be viewed as the Wronskian determinant of the functions $y, \theta y,...., \theta^m y$. Hence, it follows that these functions obey a linear
difference equation with constant coefficients
\begin{equation}\label{difference}
A_0\theta^m y+A_1\theta^{m-1}y+\dots+A_{m-1}\theta y+A_my=0
\end{equation}
that can be  used in turn to gain information on the behaviour of the function $y$ around the point $x=x_1$. When $x=x_1$ is a singular point of \eqref{fuchs},
Casorati applied the correspondence between \eqref{fuchs} and \eqref{difference} to revisit \S 3 of [Fuchs 1866], and to obtain quickly the local form of the
canonical fundamental system of solutions of \eqref{fuchs} with a cursory mention to the fact that also the structure of the integrals forming Hamburger's
subgroups [Hamburger, 1873] can be similarly obtained, an issue that will be the core of the subsequent paper [Casorati 1881a].

The operator $\theta$ admits other interpretations. For instance, by setting
$$
\theta f(x)=f(x+\omega)
$$
it is possible to view periodic functions with period $\omega$ as the fixed points of $\theta$. By introducing another operator $\theta'$ such that
$$
\theta' f(x)=f(x+\omega')
$$
it is possible to study systems of difference equations
\begin{equation}\label{Picard}
\left\{
\begin{array}{l}
A_0\theta^m y(x)+A_1\theta^{m-1}y(x)+\dots+A_{m-1}\theta y(x)+A_my(x)=0 \\
B_0\theta'^{n} y(x)+B_1\theta'^{n-1}y(x)+\dots+B_{n-1}\theta' y(x)+B_ny(x)=0
\end{array}
\right.
\end{equation}
where the $mn$ coefficients $A_i$, $B_j$ are constants.  As a special case equations \eqref{Picard} embody those arising in a problem just tackled by Picard in
[Picard, 1879] where the latter generalized the concept, due to Hermite, of doubly periodic functions of the second kind with periods $2K$ and $2\iota K'$, by
looking for uniform functions $f(x)$ such that
$$
f(x+4K)=Af(x)+Bf(x+2K)\quad\mbox{and}\quad f(x+4\iota K')=A'f(x)+B'f(x+2\iota K').
$$
The same path followed to associate \eqref{difference} to \eqref{fuchs} also works for a linear differential equation with periodic or doubly periodic
coefficients, and Casorati devoted Chapt. 7 of [Casorati, 1880b]  to investigate the structure of the general solution of these equations that will been studied
intensively also by Gaston Floquet [Floquet, 1883, 1884].

As a final interpretation of the $\theta$-operator, in Chapt. 8 of his paper Casorati considered functions in $p$ independent variables that are multiply
periodic, so that
$$
f(x_1,x_2,...,x_p)=f(x_1+\omega_1,x_2+\omega_2,...,x_p+\omega_p)
$$
holds for a set of periods $\omega_1,...,\omega_p$ and now $\theta$ is interpreted as the operator
$$
\theta:\quad f\mapsto \theta f(x_1,x_2,...,x_p)=f(x_1+\omega_1,x_2+\omega_2,...,x_p+\omega_p).
$$

\section{Chronicle of a missed prize}\label{sec:Grand Prix}

From the introduction of [Casorati 1880b] we realize the importance that Casorati attached to his paper:

\bigskip

{\it Researches that employ complex variability rely---at least in part---in an essential way upon  relations that stem from considering the way in which
functions of a complex variable behave when the variable turns around a particular value. This remark prompted me to investigate and determine once and for all,
in advance and independently of any particular study or aim, the properties and the general formulae that it is important,  and possible, to obtain to the benefit
of further, particular researches; in such a way to form a theory in itself from which all those who begin investigations on functions with complex variables from
now on can avail themselves as a common instrument.}\footnote{Le ricerche che si fanno colla variabilit\`a complessa si basano, almeno in parte, essenzialmente
sulle relazioni che scaturiscono dal considerare i modi di comportarsi delle funzioni di una variabile complessa al girare di questa intorno a suoi valori
particolari. Questa osservazione mi fece nascere il pensiero di investigare e stabilire, una volta per sempre in an\-ti\-ci\-pa\-zio\-ne ed indipendentemente da
ogni studio o scopo particolare, le propriet\`a e le formole generali che da siffatte relazioni fondamentali importa ed \`e possibile di ricavare a beneficio
delle singole ricerche ulteriori; cos\`{\i} da costituire una teorica a s\`e, della quale, come di strumento comune, possano valersi tutti coloro che d'ora
innanzi vorranno intraprendere stud\^{\i} sulle funzioni di variabili complesse.} ([Casorati, 1880b], p.10.$\equiv$ [Casorati, 1951], p. 317)

\bigskip

At the beginning of November 1879 Casorati, influenced by Brioschi, decided to apply for the {\it Grand Prix des Sciences Math\'ematiques} that the French Academy
of Sciences proposed for the year 1880, whose theme perfectly matched the topic of his researches:

\bigskip

{\it Improve in some important point the theory of linear differential equations in only one variable.}\footnote{\it Perfectionner en quelque point important la
th\'eorie des \'equations diff\'e\-ren\-tiel\-les lin\'eaires \`a une seule variable ind\'ependante.}

\bigskip

Casorati kept all the correspondence concerning [Casorati 1880b] and the prize, and in this section we outline the content of part of these letters whose original
texts are displayed in Appendix 1 (see Section \ref{sec:originals_GP} below)\footnote{We did not transcribe Fuchs's letter to Casorati, dated March 30th, 1881 on
the Stickelberger controversy, as it has been already published in [Neuenschwander, 1978].}. We should note that at the beginning of the correspondence Casorati
had not published yet either [Casorati 1880b] or anything else on the interpretation of finite differences. We also recall that the letters sent by Casorati are
known through his minutes, whereas those addressed to him are the originals. They cover more than one  year: the first letter is dated December 27th, 1879 and the
last letter was written on April 17th, 1881. Letters 1--5 contain the correspondence between Casorati and Joseph Bertrand, the {\it S\'ecretaire perpetuel} of the
French Academy. Casorati asked for general information on the competition and, curiously enough, at the beginning of the correspondence he pretended that the
applicant was a fellow-countryman of his, probably to secure anonymity: it was only after the long-waited reply from Bertrand that Casorati disclosed the identity
of the applicant. From these letters we learn that the topic for the prize was chosen by Bertrand who wished to reward Edmond Nicolas Laguerre for the results on
linear differential equations that had impressed him so much. As an aside, we note that, apparently, Laguerre did not apply for the prize. Casorati was not
willing to wait the end of the competition to publish his results but, from Bertrand replies, we also realize that published papers could be eligible for the
prize, provided that no member of the judging Committee raised formal objections.

A second group of letters (Letters 6-9) are concerned with the turbulent editorial process related to the publication of [Casorati, 1880b], whose content had been
presented to the {\it Istituto Lombardo} in February 1880. In fact, a strike of the typographers was unduly delaying the publication of his notes: besides
[Casorati 1880b], Casorati illustrated applications both to Fuchs's theory and to Floquet's theory of differential equations with periodic coefficients in
[Casorati, 1880a]. He then wrote to his friend Luigi Cremona, asking him to present a revised version of the manuscript to the {\it Annali} at the Accademia dei
Lincei in Rome. Casorati not only called attention to applications to differential equations, but also to algebraic equations that might seem more interesting to
Cremona who, in turn, presented the manuscript on March 7th\footnote{{\it Memorie della Reale Accademia dei Lincei}, {\bf 5}, S.III, 195-208, (1880)}. The long
strike in Milan went on, and the deadline of June 1st for applications to the {\it Grand Prix} was approaching: Casorati withdrew the first version of [Casorati
1880b] from the {\it Annali}, and worked on the manuscript that was finally sent to the Acad\'emie on May 22nd. As to the epigraph, he chose a passage from Blaise
Pascal's {\it Pens\'ees}:

\medskip

\noindent{\it Nous sommes si malheureux, que nous ne pouvons prendre plaisir \`a une chose qu'\`a condition de nous f\^{a}cher si elle r\'eussit mal.} ([Pascal,
1835], Chapt. IV, Art. V. p.153)

\medskip

This epigraph was inserted in the {\it billet cachet\'e} together with an explanatory note where Casorati summarized the history of the manuscript that could not
be sent in a printed form. We remark that the judging committee for the prize was formed by Bertrand, together with Ossian Bonnet, Charles Hermite, Victor
Puiseux, and Jean-Claude Bouquet, as communicated in the  {\it Comptes Rendus} of April 12th, 1880.

\medskip

Casorati applied to another competition, the second edition of the {\it Premio Bressa} proposed by the Royal Academy of Sciences in Turin (Letter 10). With
respect to the {\it Grand Prix}, the Bressa prize had a broader spectrum, as it aimed at rewarding the best work written by an Italian scholar in the period
1877--1880 on a subject ranging from pure and applied mathematics to pathology and geology. To have a chance, it was then important that the mathematicians in the
Turin Academy could appreciate the importance of Casorati's work. For this reason, since the principal analyst in Turin at that time was Angelo Genocchi, Casorati
wrote him a letter (June 21st, 1880, not published here) asking for an appointment, since he was going to visit Turin with Beltrami. In this occasion he would
have explained the importance he ascribed to his findings. The meeting with Genocchi took place at the end of June and, according to Casorati's short account of
it to Brioschi, Genocchi was favorably impressed. Thus Casorati asked Brioschi and Betti, both members of the Turin Academy, to support his paper: in spite of
these efforts, however, the 2nd Bressa prize was awarded to the natural scientist Luigi Maria d'Albertis.

At the beginning of August, with his manuscript finally published in the {\it Annali}, Casorati sent offprints of it to many mathematicians, both in Italy--for
instance to Francesco Siacci---and abroad, including Fuchs, Klein, Lindemann, Mittag-Leffler, Dillner, Kronecker, Picard, and Weierstrass. In all his accompanying
letters, he stressed the importance of the theorem on the $\Theta$ determinant. As an example, letters 11 and 12 of Section \ref{sec:originals_GP} are those that
Casorati sent to, and received by, \'Emile Picard.

So far, we followed only the correspondence concerning the {\it Grand Prix}, but at this point the events concerning the prize intertwined with a polemics
involving Casorati and Ludwig Stickelberger. By the end of 1880 Casorati received a memoir --- an {\it Akademische Antrittsschrift} [Stickelberger, 1881] --- that
accompanied Stickelberger's entry as a professor at the University of Freiburg. There the author  employed finite differences to study the behaviour of the
solutions of a linear differential equation in the neighborhood of its singular points, thus obtaining results in the same vein as those contained in Chapt. 4 of
[Casorati 1880b]. Casorati was really hurt by a sentence placed in the Introduction to the {\it Antrittschrift}, where Stickelberger ([Stickelberger, 1881], p. 4)
had referred to Casorati's procedure as {\it inadvisable} ({\it nicht rathsam}). Moreover, since Stickelberger referred to papers by Riemann [Riemann, 1857] and
Hamburger [Hamburger, 1873] where the use of the logarithm of the independent variable had been made, Casorati felt that Stickelberger was attempting to diminish
the importance of his own findings. Hence, at the end of 1880 he sent a letter to Brioschi to be published in the {\it Annali} [Casorati 1880c]. This letter is
firm in making it clear that Casorati's memoir had a broader aim than Stickelberger's, and that the applications to differential equations were just an example to
show, in a significant topic, the versatility of his ideas.

Admittedly, Casorati reply was excessive and, in my opinion, the reasons for his attitude could have been the following: {\it a}) he was concerned with the
priority of the theorem, now known as the Casorati-Weierstrass theorem, which  he had proved it in \S 88 of [Casorati 1868], and had then been published by
Weierstrass in 1876 [Weierstrass, 1876] without any reference to Casorati's work [Neuenschwander, 1978] (it should be noted, however, that there is some evidence
that Weierstrass had known the theorem since 1863 ({\it see} [Neuenschwander, 1978]; [Bottazzini and Gray, 2013], p. 436)); {\it b}) the great expectations he had
about the impact of his work that were now brought into question by Stickelberger's remark; {\it c}) the familiar concerns ({\it see} Letter 6 of Appendix 1 (see
Section \ref{sec:originals_GP} below) he had to face when working on the final version of [Casorati, 1880b] and that interfered strongly with his scientific
activity. The short letter [Casorati, 1880c] was probably also written  in a hurry, and Casorati did not properly separate relevant scientific comments from
personal feelings. A few months later, on February 6th 1881, Casorati wrote a long letter to Ossian Bonnet (Letter 13) where he separated personal animosity from
scientific judgements. There Casorati was appreciative of some aspects of Stickelberger's paper so that his judgement of [Stickelberger, 1881] contained in the
central part of the letter to Bonnet resembles a standard referee's report. First, he praises Stickelberger's use of the results obtained by Cauchy in [Cauchy,
1827], a paper devoted to symbolic calculus. There, using the notation
$$
\Delta y=\Delta_xy:=f(x+\Delta x)-f(x)
$$
for the finite difference of a function $y=f(x)$, and $Df$ for $\frac{\dd f}{\dd x}$, Cauchy had also proved ([Cauchy, 1827], p.228-230) that, if $F(x)$ is a
polynomial with {\it constant} coefficients having only {\it simple} roots $r_1, r_2, ..., r_n$, then the general solution of the differential [difference]
equation
\begin{equation}\label{cauchy}
F(D)y(x)=f(x)\qquad [F(\Delta)y(x)=f(x)],
\end{equation}
could be rewritten as
\begin{equation}\label{cauchyrepr}
y(x)= \sum_{k=1}^n \frac{1}{F'(r_k)}\frac{f(x)}{D-r_k}
\end{equation}
where $F'(r_k)$ is the shorthand for $\left.\frac{\dd F}{\dd r}\right|_{r_{k}}$ and $\frac{f(x)}{D-r_k}$ denotes symbolically the solution of the {\it first
order} equation $(D-r_k)y(x)=f(x)$, for which the representation
$$
y(x) ={\rm e}^{r_k x}\int {\rm e}^{-r_k t}f(t)\dd t
$$
exists. On the contrary, if $r_1=r_2=...=r_m=\varrho$, the first $m$ terms in the expression \eqref{cauchyrepr} had to be replaced by the {\it residue} of the
function
$$
\frac{1}{(r-z)F(z)}
$$
at $z=\varrho$, so that the contribution of these terms to the solution of \eqref{cauchy} was
$$
\sum_{k=1}^m R_{m-k}\frac{f(x)}{(D-\varrho)^k}
$$
where
$$
R_{m-k}=\frac{1}{(m-k)!}\frac{\de^{m-k}}{\de\vae^{m-k}}\left(\frac{\vae^m}{f(\varrho+\vae)}\right)_{\vae=0}
$$
and, again, $\frac{f(x)}{(D-\varrho)^k}$ represented symbolically the solution of $(D-\varrho)^ky(x)=f(x)$, given by
$$
{\rm e}^{r_k x}\int\dots\int {\rm e}^{-r_k t} f(t)\dd t^k.
$$
Actually Cauchy's version of the theorem was more general, but this application is enough for our purposes, as it is this classification of the solution according
to the presence or not of coincident roots in the characteristic equation $F(x)=0$ that bridges Cauchy's and Fuchs's theories. However, in all the equations
considered by Cauchy the coefficients of $F(x)$ are constant, and so the worlds of differential and difference equations remained distinct, though parallel, in
his work. As noted before, it is the change of variable $u:=\frac{\log x}{2\pi \iota}$ that makes it possible to export methods of finite difference calculus to
linear differential equations with variable coefficients.

The second remark made by Casorati on [Stickelberger, 1881] concerned the resort to Weierstrass's theory of elementary divisors ({\it Elementartheiler}) to
classify solutions of \eqref{fuchs}. Casorati did not think that Stickelberger's approach is optimal, and this view was the starting point of the notes he alluded
to in this letter to Bonnet, and that also appeared in the {\it Comptes Rendus} [Casorati, 1881a].

Casorati's disappointment with Stickelberger's behaviour reappears strongly in a letter, dated February 22nd, 1881, he intended to write to Frobenius, with whom
Stickelberger had scientific connections, but that was never sent. This sudden change of mind might have been caused by the the final decision of the Committee
charged to award the Grand Prix, communicated on February 14th: Georges Halphen was the recipient of the prize, and two special mentions---{\it mentions tr\`es
honorables}---were granted to two other memoirs: to Poincar\'e's essay on Fuchsian functions that remained unpublished until 1923 [Poincar\'e, 1923] ({\it see} Chapt. 3 of [Gray, 2013], for a detailed account of Poincar\'e achievements)--- and to Casorati's who, unlike Poincar\'e, decided to remain
anonymous ({\it see} vol. {\bf 92} of the {\it Comptes Rendus}, pp. 551--554 for the reports on the three best memoirs among the six manuscripts sent to the
Academy). The outcome of the competition did not satisfy Casorati who feared that his prestige as a professor could suffer a serious blow (Letter 15), despite of
the appreciation that his work received from Bertrand (Letter 14) and Darboux. In particular, the attitude of Casorati, who was prone to be bewildered in front of
misfortune, appears in his letter to Bertrand. We note that, on writing the report on behalf of the Committee, Hermite meant to point out a preference for
Poincar\'e's memoir over Casorati's since the order of the memoires by Casorati (N. 3) and Poincar\'e (N. 5) was shifted from the natural one: Casorati asked
Bertrand (on February 26th) to correct it {\it unless} it meant a difference in quality between the two memoirs.

From the historical point of view, Letters 16 and 17 are the more significant. Stickelberger replied to [Casorati, 1880c] by writing to Brioschi twice to have his
rebuttal published in the {\it Annali}: Letter 16 is the second one written by Stickelberger. Since these letters provoked another bitter reply from Casorati,
Brioschi finally decided not to proceed with publishing Stickelberger's defence. Stickelberger stated he received the paper [Casorati, 1880b] from Ferdinand von
Lindemann in August 1880, when he had already obtained his main results but, more importantly, he quoted a 1858 paper by Christoffel [Christoffel, 1858] where
Casorati's condition of linear independence---with {\it constant} coefficients---among functions depending on a variable undergoing discrete increments had
already been obtained. Casorati admitted to Brioschi that he was unaware of Christoffel's result, and promised to quote it as soon as possible --- a promise he
did not fulfill as he ceased  to work on this topic.

Christoffel's paper --- which also had a continuation first published in his Collected papers [Christoffel, 1910] --- is divided in two parts: for our purposes,
the discussion in \S\S 1-2 of [Christoffel, 1858] is the most important. There he considered $n+1$ functions $\{f(m), f_1(m),...f_n(m)\}$ with $f_i:\ZZ\mapsto
\RR$ and, supposing that they are linear dependent for the values of $m$ ranging from $m_0$ to $m_0+n+p$, with $p\ge 0$, he distinguished three cases according to
whether, for these values of $m$,

{\it a)} only a relation $Af(m)+A_1f_1(m)+\cdots+A_nf_n(m)=0$ holds;

{\it b)} at least two relations $Af(m)+A_1f_1(m)+\cdots+A_nf_n(m)=0$, $Bf(m)+B_1f_1(m)+\cdots+B_nf_n(m)=0$ hold for the {\it same} values of $m$, and so the {\it
rank} of the matrix
$$
\Atens:=\left(\begin{array}{cccc}
f(m) &f_1(m) & \cdots & f_n(m) \\
f(m+1) &f_1(m+1) & \cdots & f_n(m+1) \\
\cdots & \cdots & \cdots & \cdots \\
f(m+n) &f_1(m+n) & \cdots & f_n(m+n)
\end{array}
\right)
$$
is strictly lower than in case $a)$;

{\it c)} it is possible to single out $q+r$ intervals, each consisting of at least $n$ elements, such that in $q$ of them, case $(a)$ holds, whereas in remaining
$r$ intervals case $(b)$ holds.

\noindent[Casorati, 1880b] is consistent with case $(a)$ of Christoffel's treatment. Here Christoffel wrote down the linear system
$$
\left\{
\begin{array}{l}
Af(m)+A_1f_1(m)+\cdots+A_nf_n(m)=0 \\
Af(m+1)+A_1f_1(m+1)+\cdots+A_nf_n(m+1)=0 \\
......................................\\
Af(m+n)+A_1f_1(m+n)+\cdots+A_nf_n(m+n)=0:
\end{array}
\right.
$$
denoting the determinant of the matrix $\Atens$ by $\Delta(m)$, and the minor obtained by removing the $\mu$-th column and the last line of $\Atens$ by
$\Delta_\mu(m)$, he obtained
$$
A_\mu=\omega\Delta_\mu(m)
$$
where the factor $\omega$ is undetermined, while $\Delta(m)=0$ guarantees existence of the solution in the unknowns $A_\mu$. To prove the converse statement
Christoffel ([Christoffel, 1858], \S 2) supposed that $\Delta(m)$ vanishes at $m=m_0, m_0+1,\cdots, m_0+p$, but that it does not vanish at $m_0-1$ and at
$m_0+p+1$: then, the functions $\{f, f_1,...,f_n\}$ are linearly dependent at $m=m_0, m_0+1,\cdots,m_0+p+n$. Christoffel set
$$
\vap(m):=Af(m)+A_1f_1(m)+\cdots+A_nf_n(m)
$$
and selects the $n+1$ coefficients $A_i$ such that the $n$ equations
$$
\vap(m')=\vap(m'+1)=\cdots=\vap(m'+n-1)=0
$$
hold, for a yet unspecified value $m'$ of $m$. By adding to these equations either
$$
Af(m'-1)+A_1f_1(m'-1)+\cdots+A_nf_n(m'-1)=\vap(m'-1)
$$
or
$$
Af(m'+n)+A_1f_1(m'+n)+\cdots+A_nf_n(m'+n)=\vap(m'+n),
$$
then
\begin{equation}\label{Christoffel}
A_\mu\Delta(m'-1)=(-1)^n\vap(m'-1)\Delta_\mu(m'),\quad\mbox{or}\quad A_\mu\Delta(m')=\vap(m'+n)\Delta_\mu(m')
\end{equation}
follows, respectively. Hence, when $m'$ is set equal to one of the values $m_0$, $m_0+1$,...,$m_0+p$, the left-hand side of \eqref{Christoffel}$_2$ vanishes by
hypothesis and so
$$
\Delta_\mu(m')\vap(m'+n)=0:
$$
if at least one of the minors $\Delta_\mu(m')$ does not vanish, then $\vap(m'+n)=0$, and linear dependence of $\{f, f_1,...,f_n\}$ is guaranteed. As a consequence
of his analysis, Christoffel obtained that

\bigskip

{\it The $n+1$ assigned functions $\{f(m), f_1(m),...,f_n(m)\}$ are linearly independent for any value of $m$ if and only if their determinant $\Delta(m)$ is
different from zero for all values of $m$.}\footnote{Damit gegebene $n+1$ Functionen $\{f(m), f_1(m),...,f_n(m)\}$ f\"ur alle Werthe von $m$ linearunabh\"angig
sind, ist erforderlich und hinreichend, dass ihre Determinante $\Delta(m)$ f\"ur alle Werthe von $m$ von Null verschieden ist.} ([Christoffel, 1858], p. 288)

\bigskip

 Christoffel's paper is not mentioned even in treatises on finite differences that have rather extensive bibliographies like [N\"orlund, 1924], [Guldberg and
 Wallenberg, 1911], and [Pascal, 1897a]. It was quoted in passing by Boole ([Boole, 1880] p. 232) as a {\it remarkable paper} though for a different reason, and
 it did not escape the encyclopedic history of determinants by Muir who actually reviewed both Christoffel's [Muir, 1911, pp. 225--228] and Casorati's papers
 [Muir 1923, pp. 242--243], though he did not feel the need to note any interrelation between them, apart from inserting both of them in chapters on Wronskians.
 We note in passing that Ernesto Pascal made a short reference to Christoffel's paper in another book on determinants, for a theorem on Wronskians [Pascal,
 1897b]. Looking at the outline [Butzer, 1981] of the main achievements by Christoffel in mathematics, [Christoffel, 1858] does not occupy the first rank, and the
 episode discussed in this Section certifies once more a certain difficulty in appreciating the content of Christoffel's paper by the scientific community.

We can conclude that the novelty of Casorati's approach is the {\it interpretation} he gave to the formalism of finite differences. Curiously enough, his theorem
on the $\Theta$ determinant has been later applied as a tool in the {\it standard} theory of finite-difference calculus where, however, Casorati had been
partially anticipated by Christoffel.

\section{An unpublished note by Pincherle}\label{sec:Tannery}

In 1875 Jules Tannery wrote a long memoir [Tannery, 1875] in which he aimed at clearly explaining the main principles upon which the theory of linear differential
equations rests. In doing this, he had as models the relevant papers by Thom\'e and, more importantly, by Fuchs. The memoir was actually Tannery's thesis, which
was discussed at the Faculty of Sciences in Paris on November 28th 1874. It is divided into five parts: the first and the second part illustrate general basic
properties of the solutions of linear differential equations while the third part, on which we will focus attention, is concerned with the behaviour of a
fundamental set of solution in the neighbourhood of a singular point; finally, the fourth part is concerned with a special class of differential equations, namely
those that only admit integrals with poles of finite order at the singular points of the equations and at infinity, and the closing fifth section is devoted to
applications. Casorati, who was given a copy of the thesis by Tannery himself, expounded parts of [Tannery 1875] in his 1878--79 course devoted to differential
equations [Casorati, 1878a] and, as we shall see later on, he made some critical remarks on [Tannery 1875] that played a r\^ole in the history of Pincherle's
unpublished note. Tannery's paper was noted also by Goursat in [Goursat, 1883] where he generalized the representation of the hypergeometric series as definite
integrals depending on parameters. Actually, both Goursat and Pincherle were interested in the same theorem ([Tannery, 1875], pp. 130--132):

\medskip

{\it Let $y_1,y_2, \dots, y_m$ be $m$ functions of $x$ that are continuous everywhere, except at some singular points isolated from each other, and that are
uniform within the portions of the plane (or of the sphere) with a simple boundary that do not contain singular points: if, whenever the variable makes a complete
tour around a singular point, the new values
$$
[y_1]', [y_2]',\dots, [y_n]'
$$
of these functions are related to the previous ones by linear equations with constant coefficients, (.....), these functions are the integrals of a linear
differential equation with uniform coefficients.}\footnote{Soient $y_1,y_2, \dots, y_m$ $m$ fonctions de $x$ continues, sauf pour des points singuliers isol\'es
les uns des autres, uniformes dans les portions de plan (ou de la sph\`ere) \`a contour simple qui ne contiennent pas de points singuliers: si, lorsque la
variable fait le tour d'un point singulier, les nouvelles valeurs
$$
[y_1]', [y_2]',\dots, [y_n]'
$$
de ces fonctions sont li\'ees aux pr\`emieres par des \'equations lin\'eaires \`a coefficients constants (...), ces fonctions sont les int\`egrales d'une
\'equation diff\'erentielle lin\'eaire \`a coefficients uniformes.}

\medskip

\noindent Since the converse of this theorem is true ([Tannery, 1875], pp. 129--130) the property analyzed by Tannery {\it characterizes the integrals of a
linear, homogeneous differential equation, with uniform coefficients.} ([Goursat, 1883], p. 21). With this remarks we can now follow the letters between Casorati
and Pincherle on this topic in the period 1885-1886: the original texts are available in Appendix 2 (see Section \ref{sec:originals_Tannery} below).

On July 21st, 1885 Pincherle communicated to Casorati (see Letter 1) that he could prove Tannery's theorem quickly by resorting to the formalism of the $\theta$
operator proposed in [Casorati, 1880b].  We note that Casorati had lectured on some results then published in [Casorati, 1880b] during the course of Higher
Calculus in the academic year 1879-80, when Pincherle was in Pavia. Pincherle had only an indirect knowledge of [Tannery, 1875] via Goursat's memoir [Goursat,
1883] that he was studying during his researches {\it on integrals}.\footnote{Pincherle started reading Goursat's memoir on June 18th, 1885 as noted in Vol. V of
his {\it Ricerche e Saggi}, an impressive collection of notebooks on which Pincherle annotated his reflections on the scientific literature, and wrote down drafts
of his own papers and letters. The notebooks are kept in the library of the Mathematical Department of the University of Bologna.} Actually, Pincherle's work on
the $\theta$ operator is an {\it impromptu} since at that time he was more concerned with the representation of analytic functions by definite integrals that form
the content of [Pincherle, 1885b] and [Pincherle, 1886b] (Letters 3 and 5). In any case, Casorati (Letter 2) suggested that Pincherle read [Tannery, 1875]
because he might find new applications therein for the $\theta$ operator. While in Letter 2 Casorati seemed to consider the presence of these applications as not
so stringent a requirement if he were to present Pincherle's note for publication, he then changed his mind (Letter 4) and his doubts cooled Pincherle who,
however, was still convinced that Casorati's method provided the right way to obtain an easier proof of the theorem (Letter 5). At this point Casorati felt that
Pincherle had taken his doubts too seriously, and encouraged him in writing  the note (Letter 6). As a result, Pincherle asked Casorati for the volume of the {\it
Annales} containing [Tannery, 1875] that was unavailable in Bologna, and that might give him new ideas for a possible expansion of his note (Letter 7).

Letter 8, dated January, 9th 1886, is the most important letter of this correspondence as it contains the note that, as we shall see, will remain unpublished.
Pincherle proved Tannery's theorem under slightly weaker hypotheses, since he assumed that the coefficients could also be multi-valued functions, though of a
specific type: uniform functions of an analytic point.\footnote{These functions are uniform on their Riemann surface which, in turn, consists of a finite number of sheets. Here Pincherle follows Paul Appell who introduced this concept in his [1882a, b]. Giulio Vivanti, a student of Picherle's, devoted [Vivanti, 1887] to study their properties. In this paper, Vivanti also gave this informal definition (p. 54): {\it Uniform functions of an analytic point are uniform variables in two variables that are related by an algebraic equation.}}.  Here is the translation of Pincherle's note attached to that letter.

\medskip

Analysis: On a theorem of Mr. Tannery

\bigskip
\bigskip

This short note aims to show, by a new application, the advantages in the study of analytic functions, and of their behaviour in the neighborhood of given points,
of the method envisaged by prof. Casorati and exposed in his Memoir: ``Il calcolo delle differenze finite interpretato\footnote{{\it see a detailed account of
this Memoir in ``Bulletin de Darboux, S. II, T. VI, 1882''}} \&tc, Annali di Matematica, S.II, T. X.'' Here I apply this method to prove a theorem stated by Mr.
Tannery in his Memoir ``Propri\'et\'es des int\'egrales des \'equations diff\'erentielles lin\'eaires, Annales de l'Ec. Normale, S. II, T. IV, p. 130;'' this
theorem is remarkable both in its own, and for an important application given by Mr. Goursat.\footnote{in the Memoir ``Sur une classe des fonctions represent\'es
per des int\'egrales d\'efinies. Acta Mathematica, T.II''.}

The theorem is as follows, though in a modified form: ``Let $E_1$, $E_2$,....$E_n$ be $n$ analytic functions that are regular in the neighborhood of any point of
a simply connected set $T$, apart from a finite number of points $a_1$, $a_2$,...$a_p$, and among which there is no linear relation with constant coefficients;
if, when the variable moves around any of the points  $a_1$, $a_2$,..$a_p$ without leaving the set $T$, the new values of the functions are related to the former
ones by linear relations with constant coefficients, then these functions are the solutions of a linear differential equation with coefficients that are
single-valued in the set $T$.''

\bigskip

\noindent In fact, let us consider the function $E$
\begin{equation}\tag{1}
E=c_1E_1+c_2E_2+\dots +c_nE_n,
\end{equation}
where the coefficients $c$ are arbitrary constants; if the variable varies following a contour enclosing $a_1$, and if we denote by $\theta$ the Casorati
operator, it follows that
$$
\theta E=c_1\theta E_1+c_2\theta E_2+\dots+c_n\theta E_n
$$
and, by hypothesis,
\begin{equation}\tag{2}
\theta E=k_{1,1}E_1+k_{1,2}E_2+\dots +k_{1,n}E_n;
\end{equation}
similarly
\begin{equation}\tag{3}
\left\{
\begin{array}{l}
\theta^2 E=k_{2,1} E_1+k_{2,2} E_2+\dots+k_{2,n} E_n,\\
\dots\dots\dots \\
\theta^n E=k_{n,1}E_1+k_{n,2} E_2+\dots+k_{n,n} E_n;
\end{array}
\right.
\end{equation}
if we eliminate $1, E_1$, $E_2$,,...$E_n$ in equations (1), (2), and (3), an equation
$$
K_n\theta^nE+K_{n-1}\theta^{n-1}E+\dots K_0E=0
$$
is obtained; this equation, however (Casorati, loc. cit. \S 10), shows that $E$ obeys a linear differential equation
$$
\varphi_0 E^{(n)}+\varphi_1 E^{(n-1)}+\dots+\varphi_n E=0,
$$
 with coefficients that are single-valued in a neighborhood of $a_1$, and having $E$ as its general integral. From the same proof it also follows that
 $\varphi_0$, $\varphi_1$,...,$\varphi_n$ are single-valued in a neighborhood of  $a_2$, $a_3$,... $a_p$, whence the theorem follows.

If the set $T$ covers all the sphere once, the functions $\varphi_0$, $\varphi_1$,...$\varphi_n$ are uniform. The case in which the set $T$ is a  Riemann surface,
reduced to be simply connected by means of appropriate cuts, is not ruled out; in this case, the functions $\varphi_0$, $\varphi_1$,...$\varphi_n$ would be
uniform functions of an analytic point.

\smallskip

\noindent S. Pincherle

\bigskip

The important application by Goursat alluded to here is Theorem I of [Goursat, 1883]. There Goursat considered functions $f(x,u)$, where both $x$ and $u$ are
complex variables such that, for any fixed value of $x$, $f(x,u)$ admits only a finite number $m=n+p$ of singular values $v_1,...,v_m$, some of which
($a_1,...a_n$) are independent of $x$, whereas the remaining ones ($u_1,...,u_p$) do depend on $x$. For instance, the function
\begin{equation}\label{Goursat}
f(x,u):= (u-a_1)^{b_1-1}(u-a_2)^{b_2-1}\cdots (u-a_n)^{b_n-1}(u-x)^{\lambda-1}
\end{equation}
has $v_1=a_1$,..., $v_n=a_n$, $v_{n+1}=u_1=x$ as its critical values. Goursat considered the functions
$$
(v_i v_h)(x):=\int_{v_{i}}^{v_{h}}f(x,u)\dd u
$$
where the integration is performed along the straight line joining $v_i$ to $v_h$. To define these functions one has first to make sure that none of the singular
points $v_i$ fall on the path of integration : if, in \eqref{Goursat}, $b_i$ and $b_h$ have a positive real part, in defining
$(a_ia_h)=\int_{a_{i}}^{a_{h}}f(x,u)\dd u$ we need to suppose that none of the remaining points $a_j$ lie on the segment joining $a_i$ with $a_h$, and the
$x$-plane must have cuts ({\it coupures}, after [Hermite, 1881]) along the segment $(a_i,a_k)$, to ensure that the singularity $u_1$ does not enter the domain of
integration . However, this is not enough to guarantee that $(a_i,a_h)$ is a uniform (single-valued) function of $x$ since, when $x$ makes a complete turn along a
close curve surrounding the segment $(a_i,a_h)$, then the integral $(a_ia_h)$ turns into $(a_i,a_h)e^{\pm 2\pi\iota\lambda}$ that, in general, does not coincide
with $(a_ia_h)$: new {\it coupures} must be introduced to guarantee that $(a_ia_h)$ is uniform. Goursat applied Tannery's theorem to prove that {\it all} the
integrals $(v_iv_h)(x)$ that can be formed from a function $f(x,u)$ obey one and the same linear differential equation with uniform coefficients of order $m-1$.

Pincherle's proof simplified Tannery's proof, and Casorati himself had considered the possibility of applying the results of [Casorati, 1880b] to simplify
Tannery's analysis, though at another point. In fact, in a footnote on p. 336 of [Casorati, 1951], the editors of Casorati's collected papers --- either Luigi
Berzolari or Silvio Cinquini --- report on an unpublished manuscript where Casorati commented upon the advantages of \eqref{difference} to reobtain both the main
results of [Fuchs, 1866] as well as Tannery's results on the local behaviour of the general solution of \eqref{fuchs} as exposed on pp. 139-140 of [Tannery,
1875]. However, as we mentioned in Section \ref{sec:Grand Prix}, after the outcome of the Grand Prix, and the stress induced by Stickleberger's controversy,
Casorati avoided any publication somehow related to [Casorati, 1880b].

With Pincherle's note at his disposal, Casorati changed his mind again. On the envelope containing Letter 8, he wrote:

\medskip

{\it On going to Modena, I have to take Tannery [memoir] with me, to tell Pincherle that it would be a good thing that the note to be published should contain the
correction of pp. 132--133 of Tannery's [memoir]. Without this or without something more, I would not publish the Note.}\footnote{Andando a Modena, portar meco il
Tannery, per dire al Pincherle che starebbe bene far s{\`\i} che la Nota a pubblicarsi contenesse la correzione delle pag. 132-133 del Tannery. Senza di ci\`o o
di qualche altra cosa, io non pubblicherei la Nota.}

\medskip

\noindent Casorati and Pincherle met shortly in Modena, not far from Bologna where Pincherle taught, at the beginning of 1886, and it is clear from Letter 9 that
Pincherle followed Casorati's advice, though he was unable to find errors on the pages indicated by Casorati, apart from a trivial misprint. Casorati remarked
(Letter 10) that he had in mind something more serious than a typo. To understand his doubts, we have to summarise the content of [Tannery, 1875, pp. 132-133],
where Tannery aimed at obtaining the linear differential equation satisfied by the roots of an algebraic equation $f(x,y)=0$, $f$ being a polynomial of the $m$-th
degree in $y$, with the coefficient of $y^m$ a constant, for simplicity. The $m$ roots of this equation are functions of $x$ that, when $x$  moves around a
singular point, are simply shuffled so that by Tannery's theorem they have to satisfy a linear differential equation of order $m$ at most, the order being
strictly less than $m$ when  the roots of $f(x,y)=0$ satisfy a linear relation with constant coefficients. To find this differential equation Tannery eliminated
$y$ between the equations
$$
f(x,y)=0\qquad\mbox{and} \qquad\frac{\de f}{\de y}=0
$$
thus obtaining
$$
\vap(x)=f(x,y)A(x,y)+\frac{\de f}{\de y}B(x,y)=0,
$$
for suitable polynomials $A$ and $B$ of degree $m-1$ and $m-2$ in $y$, respectively. By taking $y$ as a function of $x$, the implicit function theorem yields
$$
\frac{\dd y}{\dd x}=-\frac{\frac{\de f}{\de x}}{\frac{\de f}{\de y}}=-\frac{B\frac{\de f}{\de x}}{B\frac{\de f}{\de y}}=-\frac{B\frac{\de f}{\de x}}{Af+
B\frac{\de f}{\de y}}=-\frac{B\frac{\de f}{\de x}}{\vap(x)}=0
$$
where $f(x,y)=0$ has been used to write $B\frac{\de f}{\de y}=Af+ B\frac{\de f}{\de y}$ and to arrange $B\frac{\de f}{\de x}$ so that it contains $y$ up to the
power $m-1$. Iteration of the procedure makes it possible to conclude that
\begin{equation}\label{reduction}
\frac{\dd^k y}{\dd x^k}=\frac{P_k(x,y)}{[\vap(x)]^k}\qquad\forall k=1,\dots, m
\end{equation}
where all the polynomials $P_k(x,y)$ contain $y$ up to the power $m-1$. From equation \eqref{reduction} Tannery concluded that {\it clearly} ({\it visiblement})
the following equation is satisfied by $y(x)$:
\begin{equation}\label{Tannery}
\frac{\dd^m y}{\dd x^m}+\sum_{k=1}^m\frac{Q_k(x)}{[\vap(x)]^k}\frac{\dd^{m-k} y}{\dd x^{m-k}}=0,
\end{equation}
where $Q_k$ are polynomials in $x$ only, since all the powers of $y$ from $y^2$ to $y^{m-1}$ can be eliminated. Casorati criticised this argument by producing a
counterexample, in the case $m=2$, by setting
$$
f(x,y)=ay^2+by+c
$$
where $a$, $b$, and $c$ are polynomials in $x$. Hence by standard manipulations he wrote
$$
\varphi(x)=ac-b^2=(ay^2+2by+c)\times a +(2ay+2b)\left(-\frac a2y-\frac b2\right)
$$
whence, following Tannery's recipe he obtained
$$
\frac{\dd y}{\dd x}=\frac{\alpha y+\beta}{a\vap(x)}\,,\qquad\frac{\dd^2y}{\dd x^2}=\frac{\gamma y+\delta}{a^2\vap^2}
$$
where $\alpha$, $\beta$, $\gamma$, $\delta$ are polynomial functions of $a$, $b$, $c$, and of their first and second derivatives with respect to $x$. If the terms
not containing $y$ are eliminated, the equation
$$
\beta\frac{\dd^2y}{\dd x^2}-\frac{\delta}{a\vap}\frac{\dd y}{\dd x}+\frac{\alpha\delta-\beta\gamma}{a^2\vap^2}y=0
$$
arises where, contrary to Tannery's claim, the coefficient $\beta$ fails to be a constant.

Actually, Casorati had had the same doubt a few years before in 1879, when he studied Tannery's memoir (p. 1286 of his unpublished mathematical diary). We note
that  Casorati's example concerns the case $m=2$ for which there are no powers of $y$ with exponent larger than 1 in \eqref{reduction}. Essentially, to keep the
order of the equation equal to 2 he eliminated from the expressions of $\frac{\dd^2y}{\dd x^2}$ and $\frac{\dd y}{\dd x}$ the terms independent of $y$. In
Pincherle's {\it Memorie e Saggi} for the year 1891 Pincherle went back to Tannery's method to obtain a differential equation satisfied by an algebraic function.
He repeated all the passages up to \eqref{reduction} and then he remarked that:

{\it It is possible to eliminate  $y^2$, $y^3$,...,$y^{m-1}$ among the $m$ equations
$$
\frac{\dd y}{\dd x}=\frac{P_1}{\vap(x)},\qquad \frac{\dd^2 y}{\dd x^2}=\frac{P_2(x,y)}{[\vap(x)]^2},\quad\dots\qquad \frac{\dd^{m-1} y}{\dd
x^{m-1}}=\frac{P_{m-1}(x,y)}{[\vap(x)]^{m-1}},
$$
whence an equation like
$$
\vap(x)^{m-1}y^{(m-1)}+a_1\vap(x)^{m-2}y^{(m-2)}+....+a_my=R.
$$
is obtained. To eliminate $R$ one has to further differentiate this equation.}\footnote{\it Fra le $m$ equazioni
$$
\frac{\dd y}{\dd x}=\frac{P_1}{\vap(x)},\qquad \frac{\dd^2 y}{\dd x^2}=\frac{P_2(x,y)}{[\vap(x)]^2},\quad\dots\qquad \frac{\dd^{m-1} y}{\dd
x^{m-1}}=\frac{P_{m-1}(x,y)}{[\vap(x)]^{m-1}}
$$
si pu\`o eliminare $y^2$, $y^3$,...,$y^{m-1}$ e viene un'equaz.[ione] della forma
$$
\vap(x)^{m-1}y^{(m-1)}+a_1\vap(x)^{m-2}y^{(m-2)}+....+a_my=R.
$$
Per eliminare $R$ bisogna derivare ancora.
}
The outcome, however, is similar to Casorati's remark because, after differentiations to get rid of the polynomial $R(x)$ one is led to an equation with a
constant coefficient in front of the the highest derivative, as in [Tannery, 1875], but now the structure of the denominators in the remaining terms differs from
that suggested by Tannery. More importantly, the order of the differential equation is larger than $m$.

In Letter 9 Pincherle attempted to prove a generalization of Tannery's result: any algebraic function satisfies a linear differential equation but, again, his
results eluded the question asked by Casorati that concerned the {\it structure} of the coefficients.

In my opinion, the example proposed by Casorati played an important r\^ole in Pincherle's decision not to publish his note: in fact, he could not find the error
in Tannery's paper, and looking for other interesting applications of Casorati's method would have probably taken him too much time. Moreover, Pincherle was soon
urged to write another short note [Pincherle 1886a], and he again resorted to Casorati (Letter 11) to have it published, if possible. In this note, Pincherle
introduced the function
\begin{equation}\label{mellin}
f(x):=\int_{(\varrho)} \vap(y)y^{x-1}\dd y
\end{equation}
where $\vap(y)$ is an analytic function without singularities along the curve $\varrho$ of the complex plane. By performing an integration by parts, and selecting
$\varrho$ so that the boundary terms disappear, he obtained
$$
-(x-1)f(x-1)=\int_{(\varrho)} \vap'(y)y^{x-1}\dd y
$$
and, by iterating the procedure, he arrived at
$$
\int_{(\varrho)} \vap^{(r)}(y)y^\lambda y^{x-1}\dd y=(-1)^r(x+\lambda-1)(x+\lambda-2)\cdots(x+\lambda-r)f(x+\lambda-r).
$$
When \eqref{mellin} is applied to the linear differential equation
$$
\sum_{r=0}^n\sum_{\lambda=0}^m a_{\lambda,r}y^\lambda \vap^{(r)}(y)=0,
$$
then it is transformed into the linear difference equation
$$
\sum_{r=0}^n\sum_{\lambda=0}^m (-1)^ra_{\lambda,r}(x+\lambda-1)(x+\lambda-2)\cdots(x+\lambda-r)f(x+\lambda-r)\theta^{\lambda-r}f(x)=0
$$
that Pincherle wrote down by resorting to Casorati's operator $\theta$.  Pincherle noted that, if the boundary terms do not drop out in \eqref{mellin}, then the
latter transformation embodies that adopted by Mellin in [Mellin, 1886 pp. 79-80], in which he obtained a correspondence between fundamental sets of solutions of
a linear homogeneous differential equation and of a functional equation. Actually, the publication of [Mellin, 1886] prompted Pincherle to have [Pincherle, 1886a]
published in order to secure his priority (Letter 12).

For a while, the interest in [Casorati, 1880b] remained quiescent, but in the next section we will see how Pincherle expounded and then extended Casorati's
approach to finite differences in his papers written after Casorati's death.

\section{The influence of [Casorati, 1880b] on Pincherle's work}\label{sec:findiff}

Pincherle's interest in finite difference equations and, more generally, in functional equations, remained alive during all his long career. In  a conference held
in 1926, and published in [Pincherle, 1926], considering the lack of attention to the calculus of finite differences in textbooks of mathematical analysis, he
wrote:

\bigskip

{\it The development of differential and integral calculus, with its numberless applications; the supremacy of the concept of limit and the criticism that, raised
by the need to make more precise this concept, led to an unavoidable revision of the foundations of calculus; the interest originated by the seminal results of
the theory of functions of one real variable on one hand, and of analytic functions on the other hand, and the appealing questions posed by the theory of
aggregates [set theory]: all this left little room to a branch of science in which the concept of limit is almost foreign at a first sight, and also those
branches of analysis in which that concept plays a minor r\^ole, I mean higher algebra, Galois's theory, the study of finite groups, did not seem to offer any
relation with the calculus of finite differences.

However, a more careful examination not only suggests that this branch of science is still worthy of attention, but [it suggests] also the expectation that,
thanks to the new paths recently opened, it can tend to a more brilliant future.}\footnote{Lo sviluppo del calcolo differenziale ed integrale, colle infinite sue
applicazioni; il si\-gno\-reg\-gia\-re del concetto di limite e le critiche, che, nate dalla necessit\`a di precisare questo concetto, conducevano alla necessaria
revisione dei fondamenti del Calcolo; l'interesse destato dai fecondi risultati della teoria delle funzioni di variabile reale d'una parte, delle analitiche
dall'altra, e dalle attraenti questioni della teoria degli aggregati: tutto ci\`o lasciava poco spazio ad un capitolo della scienza in cui il concetto di limite
appariva a prima vista pressoch\`e estraneo, mentre anche quei rami dell'analisi in cui quel concetto ha minore parte, intendo dire l'algebra superiore, la teoria
di Galois, gli studi sui gruppi finiti, non parevano offrire, col calcolo delle differenze finite, alcun addentellato. Per\`o, un esame pi\`u attento fa nascere
non solo il pensiero che questo ramo della scienza sia ancora degno di interesse, ma la fiducia che, in grazia delle nuove vie che di recente gli vennero aperte,
possa giustamente aspirare ad un pi\`u brillante avvenire.} ([Pincherle, 1926], pp. 233-234)

\bigskip

In my opinion, Casorati's views on the calculus of finite differences belong to these {\it new paths} that were given to this discipline. In this section I will
consider only those of Pincherle's results that can be related to [Casorati, 1880b]: essentially, they were published in the last decade of the 19th century and
they form the backbone of Chapt. X of the monograph on functional calculus written by Pincherle and Ugo Amaldi in 1901. It should be clear that Pincherle also
contributed to other areas of finite differences. Notably, the close relation between difference and differential equations that had been revealed new aspects
after the publication of [Poincar\'e, 1885]  soon captured Pincherle's attention, leading him to the concept of a {\it distinguished integral} of a difference
equation [Pincherle, 1892]: for the relations between Poincar\'e and Pincherle I refer the reader to Chapt. VI of [Bottazzini, 1994] and to [Dugac, 1989], where a
couple of letters between the two mathematicians have been published. Moreover, it is from finite difference equations that Pincherle was led to generalize the
algorithm of continued fractions ({\it see}, for instance, [Pincherle, 1890]), a topic on which he had a correspondence with Charles Hermite.

Turning back to the $\theta$ operator, we start from the Academic year 1893--94, when Pincherle gave a course on hypergeometric functions, a long account of which
was then printed in [Pincherle 1894]. There, for pedagogical reasons, Pincherle inserted several sections devoted to difference and differential equations where
he gave a detailed exposition of [Casorati, 1880b], reproducing part of the results and qualifying Casorati's method as {\it extremely ingenious} (genialissimo)
because

\bigskip

{\it [it] offers an unrivalled simplicity and an interest that is not only scientific but also pedagogical.}\footnote{offre una semplicit\`a senza pari ed
interesse scientifico non meno che didattico.} ([Pincherle, 1894], p. 210-211.$\equiv$ [Pincherle, 1954], p. 275)

\bigskip

In Chapt. 2 of [Pincherle, 1894], devoted to linear difference equations, Pincherle emphasized the importance of the theorem on the $\Theta$ determinant, of which
he gave a proof that coincides, apart from a slight change of notation, with the original one proposed by Casorati. Then in Chapt. IV (\S 27) Pincherle reviewed
the content of Chapter 2 of [Casorati 1880b] where applications to Fuchs's theory of differential equations had been illustrated. It is interesting to remark that
whenever Pincherle considered the distinctions of solutions in subgroups, he never mentioned [Casorati, 1881a], neither in [Pincherle, 1894] nor in [Pincherle and
Amaldi, 1901] whose Chapt. XII is devoted to this subject.

If in [Pincherle, 1894] there are no original results on the $\theta$ operator, in 1895 Pincherle published a series of papers [Pincherle, 1895a,b,c] where the
{\it operatorial} nature of $\theta$ plays a major r\^ole. In these papers Pincherle laid the basis of his functional calculus, based upon the properties of
linear operators ({\it operazioni distributive}) acting on a set of {\it analytic} functions, or to a properly delimited subset.

In [Pincherle, 1895a] Pincherle, soon after the definition of $\theta$, quoted Casorati to emphasize the large spectrum of distributivity endowed by this
operator, recalling equations \eqref{distributive1} and \eqref{distributive2}. Then, in \S 2 he started to study the properties of the set of operators
\begin{equation}\label{lineari}
    F:= a_0(x)+a_1(x)\theta+ a_2(x)\theta^2+\dots+a_m(x)\theta^m,
\end{equation}
where $\{a_{k}(x)\}$ are given functions. Pincherle called operators like \eqref{lineari} linear difference-forms ({\it forme lineari alle differenze}) of the
$m$-th order. They map a function $f(x)$ into
$$
f\mapsto F(f):= a_0(x)f(x)+a_1(x)f(x+1)+ a_2(x)f(x+2)+\dots+a_m(x)f(x+m).
$$
Up to \S 20, [Pincherle, 1895a] is devoted to forms of finite order and there Pincherle did not insist on the nature of the functions $a_i(x)$: $a_0(x)$ is only
supposed not to vanish at the points $x+k$, so that after a rescaling, it can be set equal to unity. On the contrary, in the second part of [Pincherle, 1895a],
devoted to series in $\theta$, he sketched the elements of the systematic presentation of the theory contained in [Pincherle and Amaldi, 1901] where much
attention is devoted to characterize a suitable general functional space ${\cal D}$ where $\theta$ is defined. Pincherle considered a region $\acc_0$ of the
complex plane lying in the strip $\Re(x)\in [0,1]$ and then he defined the region
$$
\acc:=\bigcup_{k=-\infty}^\infty\acc_k
$$
where $\acc_k$ is obtained by translating $\acc_0$ of $k$ units along the real axis. In any $\acc_k$ Pincherle considered a meromorphic function $\alpha_k(x)$ and
defineed another meromorphic function $\alpha(x)$ in $\acc$ such that
$$
\alpha(x)= \alpha_k(x)\quad \mbox{if $x\in\acc_k$}:
$$
the set ${\cal D}$ of all meromorphic functions defined in this way is, by construction, an invariant set for the operator $\theta$ and Pincherle considered forms
\eqref{lineari} with coefficients belonging to this set. (In particular, if $\alpha_k(x)$ is independent of $k$, we have linear difference-forms with constant
(periodic) coefficients.)

Pincherle was careful about the algebraic properties of the sets of operators he considered, not only  those in the class \eqref{lineari} but, in general,
throughout his functional calculus. As an early example of this attitude, we can quote the papers [Pincherle, 1884b, 1885a] that were markedly influenced by the
works of Kronecker [Kronecker, 1882] and Dedekind and Weber [Dedekind and Weber, 1882]. To give a flavour of his attitude in the present context, we now look at:
$a)$ the extension of the {\it Ruffini rule} to form the quotient and the remainder obtained when a form $F$ is divided by the first-order form
$E_a:=\theta-a(x)$; $b)$ a criterion to decide when a form $B$ divides another form $A$. We  follow the presentation as given in [Pincherle and Amaldi, 1901],
which is more complete than the original treatment in [Pincherle, 1895a].

After defining the non commutative product of a form $A$ of order $m$ time a form $B$ of order $n$ as the form $AB$ mapping $f(x) \mapsto A[B(f(x))]$,
Pincherle aimed at extending the classical Euclidean algorithm to linear difference-forms. By assuming $n\le m$ he looked for a third form $\Gamma$, of order
$m-n$, such that
$$
A-\Gamma B
$$
is of order $n-1$, at most. By setting
$$
A:=\sum_{h=0}^m \alpha_{h}(x)\theta^h\,,\qquad B:=\sum_{k=0}^n \beta_{k}(x)\theta^k\,,\qquad\mbox{and}\qquad \Gamma:=\sum_{i=0}^{m-n} \gamma_{i}(x)\theta^i\,,
$$
the requirement on $\Gamma$ yields the following linear system for the $m-n+1$ coefficients $\gamma_i(x)$:
\begin{equation}\label{quotient}
\left\{
\begin{array}{l}
\alpha_m(x)-\gamma_{m-n}(x)\beta_n(x+m-n)=0 \\
\alpha_{m-1}(x)-\gamma_{m-n}(x)\beta_{n-1}(x+m-n)-\gamma_{m-n-1}(x)\beta_n(x+m-n-1)=0 \\
............................................\\
\alpha_{m-r}(x)-\sum_{j=0}^r\gamma_{m-n-j}(x)\beta_{n-r+j}(x+m-n-j)=0 \\
............................................
\end{array}
\right.
\end{equation}
that can be solved provided that
\begin{equation}\label{solvability}
\beta_n(x)\beta_n(x+1)\cdots\beta_n(x+m-n-1)\beta_n(x+m-n)\neq 0,
\end{equation}
a condition automatically satisfied in [Pincherle, 1895a] where $\beta_n(x)\equiv 1$. When \eqref{solvability} is satisfied, the form $\Gamma$ is uniquely
determined, and it can be called the {\it quotient} of the forms $A$ and $B$. When the remainder vanishes one has $A=\Gamma B$ so that $B$ {\it divides} $A$. When
$B=\theta-\gamma$ the conditions \eqref{quotient} become
\begin{equation}\label{Ruffini}
\left\{
\begin{array}{l}
\beta_{m-1}(x)=\alpha_m(x) \\
\beta_{m-2}(x)=\alpha_{m-1}(x)+\alpha_m(x)\gamma(x+m-1) \\
\beta_{m-3}(x)=\alpha_{m-2}(x)+\alpha_{m-1}(x)\gamma(x+m-2)+\alpha_{m}(x)\gamma(x+m-1)\gamma(x+m-2) \\
............................................ \\
\beta_0(x)=\alpha_1(x)+\alpha_1(x)\gamma(x+1)+\dots+\alpha_m(x)\gamma(x+1)\cdots\gamma(x+m-1),
\end{array}
\right.
\end{equation}
that is the analogue of the classical Ruffini's rule to obtain the coefficients of the quotient between a polynomial and binomial $x-a$.

In terms of higher algebra, Pincherle defined irreducibility of a linear difference-form $A$ with coefficients in a set $\Omega$ that  is a field ${\cal D}$ say,
({\it campo di razionalit\`a}, i.e the Italian term for {\it Rationalit\"atsbereich}): $A$ is irreducible in $\Omega$ if there is no linear difference-form with
coefficients in $\Omega$ that divides $A$ ([Pincherle, 1895a, \S6], [Pincherle and Amaldi, 1901, \S 295]. In \S 296 of [Pincherle and Amaldi, 1901], the analogue
of a lemma by Galois is stated:

\bigskip

{\it Let $\Phi$ be a form of the order $n$, irreducible in a certain field, and let $\Phi_1$ another form of the order $m\ge n$ that shares a root $\vap$ with
$\Phi$, and whose coefficients belong to the same field. Then, $\Phi_1$ is necessarily divided by $\Phi$.}\footnote{Sia una forma $\Phi$, d'ordine $n$,
irriducibile in un determinato campo di razionalit\`a, e sia $\Phi_1$ una forma di ordine $m\ge n$, i cui coefficienti appartengono a quel campo, e che abbia una
radice $\vap$ comune con $\Phi$. La $\Phi_1$ \`e allora necessariamente divisibile per $\Phi$.}

\bigskip

We recall that in Pincherle's language a root of a linear-difference form \eqref{lineari} is an element $\omega(x)$ of ${\cal D}$ such that $F(\omega(x))=0$.
Finally, in \S 297 of [Pincherle and Amaldi, 1901] the dependence of irreducibility of a form on its domain is mentioned, together with the fact that irreducible
forms can be made reducible by adjoining suitable functions to the field where the coefficients lie.

We now turn to the second part of [Pincherle, 1895a] where Pincherle introduced {\it series} in $\theta$. There he clearly stated a difference between the
traditional way of conceiving symbolic calculus and his own. As for series like
\begin{equation}\label{serietheta}
\sum_{n=0}^\infty \alpha_n(x)\theta^n
\end{equation}
he remarked ([Pincherle, 1895a] pp. 107-108) that

\medskip

{\it These series (...) have been conceived more as a concise style of symbolic representation rather than a means apt to yield the element of a genuine theory:
here, we aim at showing how, by placing suitable constraints on the functions that are acted upon by these series, they can be taken as the basis of a truly
rigorous calculus.}\footnote{Tali serie (...) sono state considerate pi\`u come un modo conciso di rappresentazione simbolica che come atte  fornire l'elemento di
una vera teoria: noi qui ci proponiamo di mostrare come, assoggettando a convenienti limitazioni le funzioni su cui si opera, queste serie siano suscettibili di
essere prese a fondamento di un calcolo perfettamente rigoroso.}

\medskip

As a consequence of this programme, Pincherle explored sufficient conditions for the expressions \eqref{serietheta} to make sense. He considered functions
$\alpha_n(x)$ that belong to ${\cal D}$, and defined the {\it functional field of convergence} ({\it campo funzionale di convergenza}) of \eqref{serietheta}
referred to the point $x_0\in\acc$ as the subset of ${\cal D}$ such that \eqref{serietheta} converges absolutely in a neighbourhood of $x_0$ proving that any
series \eqref{serietheta} admits a non-empty field of convergence. As an example of an elementary condition guaranteeing the convergence of \eqref{serietheta}
Pincherle  [1895a, p.109] ({\it see} also \S 311 of [Pincherle and Amaldi, 1901], p. 250) stated that:

\bigskip

{\it Let $\{\lambda_n(x)\}$ be a sequence of functions that are finite in the set $\acc_0$. If the series
$$
\sum_{n=0}^\infty \alpha_n(x)\lambda_n(x)
$$
is absolutely convergent, and if the element $\vap(x)$ is such that, $\forall x\in\acc_0$,
$$
\vap(x+n)\le \gamma(x)|\lambda_n(x)|
$$
for $n\ge n_0$, where $\gamma(x)$ is both positive and finite in a neighbourhood of a point $x=x_0\in\acc_0$, then $\vap(x)$ belongs to the functional field of
convergence of \eqref{serietheta}, relative to the point $x_0$.}

\bigskip

Let us now turn attention to the concept of {\it functional derivative} $A'$ of a linear operator $A$ defined by Pincherle in [Pincherle, 1895b] as
\begin{equation}\label{funcder}
A'(\vap):=A(x\vap)-xA(\vap),
\end{equation}
provided that $\vap(x)$ and $x\vap(x)$ belong to the domain of $A$. The Casorati operator $\theta$ coincides with its derivative or, equivalently, it satisfies
$$
\theta'=\theta
$$
and so it plays a distinguished r\^ole, analogue to the exponential function in classical calculus. Actually, in [Pincherle, 1895b], Pincherle showed that an
operator $A$ satisfying $A'=A$ has the form $A=M\theta$, where $M$ is the multiplication by a suitable function. In fact, by setting $\vae(x)=A(1)$ and by using
\eqref{funcder}, the requirement $A'=A$ implies that
$$
A(x)=A(x\cdot 1)=(x+1)\vae(x),...,A(x^n)=(x+1)^n\vae(x)
$$
so that, {\it if} $\vap(x)$ belongs to the set ${\cal S}^1$ of functions whose series expansions centered at the origin have radius of convergence larger than 1,
then
$$
A(\vap(x))=\vae(x)\vap(x+1)=\vae(x)\theta.
$$
As we have already seen, Casorati emphasized that the operator $\theta$ is distributive also with respect to other operations besides sum. Pincherle deepened this
property in [Pincherle, 1895b,c] where he introduced the operation of {\it substitution} (sostituzione), denoted with $S$ or $S_a$, that amounts at replacing
$\vap(x)$ with $\vap(a(x))$, where $a(x)$ is an assigned function.  By definition, the operation $S_a$ satisfies
$$
S(\vap\psi)=S(\vap)S(\psi)\qquad\mbox{and}\qquad S\left(\frac\vap\psi\right)=\frac{S(\vap)}{S(\psi)}
$$
like the operator $\theta$ that clearly corresponds to the particular choice $a(x)=x+1$ ({\it see} also \S 122 of [Pincherle and Amaldi, 1901]). In fact, it
satisfies the first order {\it symbolic} differential equation
\begin{equation}\label{firstsymb}
S_a'=(a(x)-x)S_a
\end{equation}
which, as we have just seen for the operator $\theta$, is a characteristic property of substitutions since any solution of \eqref{firstsymb} is of the form $MS$,
where $S$ is a substitution, and $M$ is the multiplication by a certain function. Proceeding in  parallel with Fuchs's theory of linear differential equations,
Pincherle considered $n$-th order symbolic differential equations
\begin{equation}\label{nsymb}
\lambda_0(x)A^{(n)}+\lambda_1(x)A^{(n-1)}+\dots +\lambda_n(x)A=0
\end{equation}
and he proved that $n$ operators $A_1$, ..., $A_n$ solving \eqref{nsymb} are linearly dependent if and only if the determinant
\begin{equation}\label{Grevy}
\Delta:=\left|
\begin{array}{llll}
A_1 & A_2 & ... & A_n \\
A_1' & A_2'&...& A_n' \\
... & ... & ...& ... \\
A_1^{(n-1)} & A_2^{(n-1)} & ... & A_n^{(n-1)}
\end{array}
\right|
\end{equation}
is equal to zero ({\it see} also Chapt. VI of [Gr\'evy, 1894] for an analogue statement).  Then he proved that it is possible to build a fundamental system of
solutions of \eqref{nsymb} starting from a particular set of substitution operators, defined as $S_{a_{n}}$, where the $n$ functions $a_1(x),...,a_n(x)$  are the
roots of the characteristic equation
$$
f(z):=\lambda_0(x)(z-x)^n+\lambda_1(x)(z-x)^{n-1}+\dots+\lambda_{n-1}(x)(z-x)+\lambda_n(x)=0.
$$
Chapter XIV of [Pincherle and Amaldi, 1901] is devoted to a detailed study of the expressions
$$
\sum_{k=0}^m \alpha_k(x)S_\mu^k
$$
called linear substitution forms ({\it forme lineari alle sostituzioni}). There,  a parallel of the theory to the one already illustrated for linear-difference
forms is drawn by referring to results by Gr\'evy [Gr\'evy, 1894].

 Given that, besides the Wronskian and the Casoratian determinants, there is also \eqref{Grevy}, the question arises {\it of looking for larger classes of
 operations for which a theorem, analogous to that of the Wronskian, holds}\footnote{ricercare se esistano classi pi\`u ampie di operazioni, per le quali sussista
 un teorema analogo a quello del Wronskiano.} ([Pincherle and Amaldi, 1901], p. 430)

Chapter XV of [Pincherle and Amaldi, 1901] is devoted to answering this question by looking at the common structure of the ``Wronskian-like" theorems: to find out
a necessary and sufficient condition for $n$ regular enough functions $\vap_1,...,\vap_n$ to satisfy a relation
$$
c_1\vap_1+c_2\vap_2 +\cdots+c_n\vap_n=0
$$
where not all the coefficients $c_i$ vanish. The coefficients $c_i$ are {\it constant} with respect to a suitable {\it linear} operator $A$: if
$A=D=\frac{\dd}{\dd x}$
then $c_i$ are really constants; if $A=\theta$ then $c_i$ are periodic functions with period equal to 1; if $A=S_\mu$, then $c_i(x)=c_i(\mu(x))$. Pincherle and
Amaldi note that the crucial steps to obtain a Wronskian theorem are: {\it a}) to show that
$$
A(c_1)=0
$$
and {\it b}) to show that a {\it multiplication theorem} holds, that is, a particular relation among $ A(\vap\psi)$, $A(\vap)$ and $A(\psi)$ whose meaning,
however, slightly varies on passing from  $D$ to $\theta$ and to $S_\mu$:
$$
D(\vap\psi)=D(\vap)\psi +\vap D(\psi)\qquad \theta(\vap\psi)=\theta(\vap)\theta(\psi)\qquad S_\mu(\vap\psi)=S_\mu(\vap)S_\mu(\psi).
$$
As a consequence,  Pincherle and Amaldi restricted their attention to single-valued linear operators such that
\begin{equation}\label{multiplication}
A(\vap\psi)=\alpha_{00}\vap\psi+\alpha_{01}\vap A(\psi)+\alpha_{10}\psi A(\vap)+\alpha_{11}A(\vap)A(\psi)+\beta_1\vap A(\vap)+\beta_2\psi A(\psi),
\end{equation}
where $\alpha_{ij}$ and $\beta_i$ are functions that are not completely independent. In fact, by imposing that $A(\vap\psi)=A(\psi\vap)$ and then testing
\eqref{multiplication} on $\psi=1$ they got the simpler expression
\begin{equation}\label{multiplication2}
A(\vap\psi)=\xi(\alpha\xi-1)\vap\psi+(1-\alpha\xi)[\vap A(\psi)+\psi A(\vap)]+\alpha A(\vap)A(\psi)
\end{equation}
where $\alpha:=\alpha_{11}$ and $\xi:=A(1)$ are independent from one another. If $\alpha=0$, the most general linear operator satisfying \eqref{multiplication2}
is $$
A= (\xi_1-\xi x)D+M_\xi,
$$
where $\xi_1=A(x)$ and $M_\xi$ represents pointwise multiplication by $\xi$. On the other hand, if $\alpha\neq 0$ such a r\^ole is played by
$$
A=\frac{1}{\alpha}S_\mu+M_{\xi-\frac{1}{\alpha}}:
$$
hence, linear differential forms and substitution forms exhaust the class of linear operators satisfying \eqref{multiplication2}. For this class of operations, an
analogue of the $\Theta$-determinant holds ({\it see} [Pincherle and Amaldi, 1901], \S 455). With Pincherle's work the results of [Casorati, 1880b]  found their
natural place within functional calculus.

\section{Conclusions}

The main obituaries of Pincherle, namely [Amaldi, 1937] and [Tonelli, 1937] refer to the influence exerted by Casorati on some aspects of Pincherle's research.
Resorting to unpublished documents kept in the Casorati {\it Nachlass}, we have examined some aspects of their relationship. In particular, we have found an
unpublished note by Pincherle centered on a paper by Tannery [Tannery, 1875] where Casorati's $\theta$ operator is used in the new interpretation envisaged by
Casorati himself to export finite-difference methods in the realm of complex analysis. Casorati's paper [Casorati, 1880b] had the most influence on Pincherle, as
we showed in Section \ref{sec:findiff} through an overview of Pincherle's work on this topic, culminating in the publication of his monograph [Pincherle and
Amaldi, 1901]. Moreover, other documents of the {\it Nachlass} made it possible to investigate the complex history behind Casorati's participation in the 1880
{\it Grand Prix de Sciences Math\'ematiques}, and to add some detail on the controversy he had with Stickelberger. In particular, we have seen that Stickelberger
claimed that the theorem on the now-called Casorati's determinant had already been proved by Christoffel more than 20 years before the publication of [Casorati,
1880b], though in the context of classical finite-difference calculus. We plan to deepen other aspects of Pincherle's early career in future publications.

\section{Appendix 1: Original letters concerning the 1880 Grand Prix des Math\'ematiques}\label{sec:originals_GP}

{\it Letter 1}: {\underline{Casorati to Bertrand}}

\begin{flushright}
Pavia, 27 dicembre 1879\\
\end{flushright}

\medskip

Caro ed illustre Amico

\medskip

 Un mio compaesano, discontinuo ma passionato cultore dei nostri stud\^{\i}, \`e riuscito a scoperte analitiche di grande momento, segnatamente efficaci nelle
 moderne ricerche basate sulla variabilit\`a complessa. (...) Egli ha gi\`a potuto impossessarsi dei lavori [che sono] stati fatti in Germania, e vi ha subito
 introdotto belle ed importanti semplificazioni, cos\`{\i} da avere la certezza di figurare con onore tra concorrenti anche fortunatissimi. Ma poich\'e i suoi
 concepimenti hanno un'importanza generale, e ripeto grande, e non subordinata al tema su citato, cos\`{\i} egli, giustamente nell'interesse dei suoi stud\^{\i} e
 nel suo, non vorrebbe tenerli segreti sino alla fine del 1880. Per\`o d'altra parte, egli vorrebbe comportarsi in modo da conseguire la grande soddisfazione di
 un premio dell'Accademia di Parigi, la quale vedr\`a con piacere uscire in luce le nuove idee per mezzo delle sue pubblicazioni\footnote{Casorati added the
 following comment, that was not included in the letter: {\it Qualunque sia per essere l'esito del concorso, deploro di aver scritto questi pensieri. Essi erano
 naturali quando un fatale inganno mi fece credere commutativi due simboli d'operazione e mi condusse quindi a conseguenze d'incredibile importanza. Ma non dovevo
 esternarli subito, come feci, ad uno straniero; dovevo lasciar passare almeno qualche giorno, secondo le buone regole, che non si violano impunemente. Che cosa
 potrebbe aver detto Bertrand, se avesse letta bene la mia lettera? E perch\'e, dopo riconosciuto il mio traviamento, non gli scrissi per avvertirvelo? Ci\`o che
 rimaneva di accertato nelle mie ricerche era pur sempre meritevole, a mio giudizio, di molta considerazione.}}. Interessatissimo pel mio compaesano, pensai di
 ricorrere in segreta confidenza a voi per consiglio; fidente nell'interessamento che voi avete sempre altamente mostrato per i nostri stud\^{\i}, e nella
 benevolenza che avete sempre accordato a chi vi scrive. Io dunque chiedo non al Segretario Perpetuo, ma all'amico benevolo, come dovrebbe il mio compaesano
 comportarsi? Se stampasse anonime una parte delle cose sue, potrebbero essere ancora valutate in un manoscritto presentato al Concorso? [...]

\begin{flushright}
Il vostro devotissimo\\
F. Casorati\\
della Universit\`a di Pavia
\end{flushright}

\medskip

{\it Letter 2}: {\underline{Casorati to Bertrand}}

\begin{flushright}
Pavia, 28 genn. 80
\end{flushright}

\medskip

Ch. mo Prof. (Bertrand)

(...) Questo silenzio, ch'io credo non meritato, ha gi\`a naturalmente prodotto il suo triste effetto, di scoraggiare chi forse meritava incoraggiamento, di
dissuaderlo dal continuare in quelle applicazioni che aveva cominciato felicemente, con quell'ardore che guida sempre a qualche risultato non volgare.

Ho sempre creduto che, come la grande maggioranza delle persone colte guarda a Parigi con ammirazione ed affetto, come alla gran madre d'ogni progresso
dell'umanit\`a, cos\`{\i} voi Parigini, e voi particolarmente membri dell'Istituto, dobbiate guardare con benevolo affetto tutti gli operai del progresso, pi\`u o
meno distinti che sieno, giovani o no, purch\'e onesti, purch\'e sinceramente intenti al meglio dell'umanit\`a.

Ma non aggiunger\`o altro, perch\'e, rattristato come sono, vi riuscirei troppo molesto.

\begin{flushright}
Obbl.$^{\rm{mo}}$ vostro \\
F. Casorati \\
dell'Universit\`a di Pavia
\end{flushright}

\medskip

{\it Letter 3}: {\underline{Bertrand to Casorati}}

\begin{flushright}
Paris, 31 janvier 1880
\end{flushright}

\medskip

\begin{flushright}
Tr\`es cher monsieur
\end{flushright}

(...) Je ne puis malhereusement engag\'e \`a l'avance ni l'acad\'emie, ni la Commission, qui n'est pas m\^{e}me nomm\'ee et qui sera souveraine. Elle fera tout,
j'en ai la convinction, pour courronner une m\'emoire importante et les questions du forme, celles m\^{e}me du publications ant\'erieure, n'ont pas l'habitude du
pr\'evaloir sur le m\'erite du fond. Je peut vous citer l'example du Kummer qui a obtenu le prix proponi pour le th\'eor\`eme de Fermat, {\it sans avoir
concouru}, son m\'emoir publi\'e en langue allemande depuis plusieurs ann\'es, et traduit en Francais depuis plus d'un an, a \'et\'e jug\'e digne du prix et l'a
obtenu au grand \'etonnement du l'auteur qui n'y songeait pas. Je conseillerai donc \`a votre ami d'imprimer en langue italienne les principaux r\'esultats et
d'envoyer un m\'emoire d\'etaill\'e qui, j'en ai la convinction, sera lu et jug\'e avec le m\^{e}me empressement que si tout etait in\'edit. Le secret qui doit
\^{e}tre gard\'e sur le nom de l'auteur, sera, il est vrai, viol\'e: un juge formaliste pourrait y voir une difficult\'e et je n'ai pas le droit {\it d'affirmer}
que cela n'arrivera pas, mais je ne pense pas, d'apr\`es mon habitudes et les ide\'es que je connais \`a mon Confr\`eres, que cela soit \`a craindre.

Je ne veux pas vous cacher que le sujet a \'et\'e propos\'e pour moi et que j'avais pens\'e en le r\'edigeant \`a la possibilit\'e de r\'ecompenser les travaux de
Mr. Laguerre sur les \'equations differenti\'elles lineaires et particuli\`erement sur les invariants. Mr. Laguerre lui m\^{e}me ignore d'ailleurs cette arri\`ere
pens\'ee dont je n'ai pas fait part \`a mes Confr\`eres, mais si lui vient l'id\'ee de concourir, il sera, vous le voyer, dans le m\^{e}me cas que votre ami.
(...)

\begin{flushright}
J. Bertrand
\end{flushright}

\medskip

{\it Letter 4}: {\underline{Casorati to Bertrand}}

\begin{flushright}
Pavia, 6 febbraio 1880
\end{flushright}

\begin{flushleft}
Carissimo Signore (Bertrand)
\end{flushleft}

(...) Ora vi faccio una confessione. Il mio compaesano \`e quegli stesso che vi scrive. Ma venendo ormai all'argomento, vi dir\`o, che, ritornato con lena agli
stud\^{\i} dopo un lungo soggiorno sulle Alpi, fui condotto ad interpretare il calcolo delle differenze finite, diretto ed inverso, in modo da farlo diventare un
potente ausiliario delle moderne ricerche basate sulla variabilit\`a complessa. Naturalmente, a questa interpretazione aggiunsi alcuni teoremi affatto nuovi.
Allorch\'e poi l'utilit\`a di queste mie cose mi parve largamente dimostrata dall'applicazione allo studio, ora tanto in voga, delle equazioni differenziali
lineari, pensai che avrei potuto concorrere al premio proposto dalla vostra Accademia su quest'ar\-go\-men\-to, e mi risolvetti a scrivervi la mia prima lettera.

Per\`o voi avete interpretato questa lettera in un senso che io era lontanissimo dal darle. Dite "Je ne puis malhereusement engag\'e \`a l'avance ni l'Acad\'emie,
ni la Commission, qui n'est pas m\^{e}me nomm\'ee et qui sera souveraine." Io lo credo bene! E non avrei mai imaginato [sic!] di pretendere cose ingiuste. Ma
tutto ci\`o sia come non detto. L'essenziale ora si \`e che la vostra lettera mi ha rialzato e soddisfatto pienamente.

Accettando il vostro consiglio, far\`o stampare subito negli {\it Annali di Ma\-te\-ma\-ti\-ca} di Milano una parte della interpretazione suddetta con
l'applicazione per ora alle propriet\`a fondamentali degli integrali delle equazioni differenziali lineari a coeff.[icienti] monodromi. E non mander\`o alcun
manoscritto; parendomi non necessario giacch\'e mi basta di sapere dalla vostra lettera, che la Comm.[issione] potrebbe anche prendere in considerazione lavori
pubblicati per la stampa. Io non avrei nessun gusto di fare concorrenza al sig. Laguerre. Sono abbastanza pago di poter credere che ove il sig. Laguerre non si
risolvesse a concorrere, potrebbe essere presa in considerazione anche la Memoria che io adesso pubblicher\`o, e quelle altre sulle equazioni lineari ancora e su
altre ricerche con variabilit\`a complessa che avessi agio di pure redigere in netto e stampare.

Perdonatemi tutto questo disturbo, pur troppo grande per voi certamente accerchiato da mille cose, ed aggradendo l'attestazione della profonda mia stima e del mio
affetto, vogliate conservarmi la vostra preziosa benevolenza.

\begin{flushright}
Dev.$^{\rm{mo}}$ vostro \\
F. C.
\end{flushright}

\medskip

{\it Letter 5}: {\underline{Bertrand to Casorati}}

\begin{flushright}
Paris, 15 f\'evrier 1880 \\

\medskip

Cher monsieur,
\end{flushright}

J'ai recu votre lettre avec  grand plaisir tr\`es heureux de voir qu'il ne reste rien du malentendue dont j'ai \'et\'e, bien sans int\'entions, la cause.

Je ne puis, comme je vous j'ai dit, prendre aucun engagement au nom d'une Commission qui n'existe pas encore et vous l'avez parfaitement compris mais je crois de
vous avoir laiss\'e prendre pour une habitude de l'acad\'emie ce qui a \'et\'e fait une fois or deux seulement dans des circumstances exceptionelles. Les
commissions {\it peuvent} d\'ecerner le prix \`a des m\'emoirs imprim\'e et non envoy\'e au concours mais elles ne l'ont fait que tr\`es rarement, en l'absence de
plis r\'egulierement envoy\'es au concours et digne de disputer le prix et lorsqu'aucune membre de la Commission n'envoquait la lettre des programmes pour s'y
opposer.

Si donc vous avez obtenu ....\footnote{It was hard to decode Bertrand's manuscript here. It seems to me: {\it comme j'en suis content}} des r\'esultats importants
et si votre d\'esir bien natur\'el est de les faire courronner par l'acad\'emie des Sciences de Paris, il est prudent de les lui envoyer et s'ils ont \'et\'e
imprim\'es par entier, vous pouver citer le pubblications, mais en mettant sous les yeux des commissaires un\footnote{Same remark as before. It looks like: {\it
essai}.}... complet de votre travail. En procedent autrement, vous ne rendrez pas le succ\`es impossible parce qu'il y a des examples qui sont tout semblables,
mais vous en diminuerez singulierement les chances. [....]

\begin{flushright}
J. Bertrand
\end{flushright}

\medskip

{\it Letter 6}: {\underline{Casorati to Cremona}}

\begin{flushright}
Pavia, 25 febb. 1880
\end{flushright}

Caro Cremona,

\medskip

Sono scampato anch'io da una grave disgrazia. La mia Eugenia veniva assalita da pleurite con minaccia di tifo, cos\`{\i} da metterne per qualche giorno in
pericolo la vita. Fortunatamente il tifo non si svilupp\`o e la pleurite va attenuandosi. Bench\'e sventata, questa minaccia ha tanto prostrato le mie forze che
non sono peranco capace di applicarmi per un'ora allo studio. L'essersi ammalata in Milano, e il dovervisi naturalmente fermare ancora molto tempo fa s\`{\i} che
non posso pi\`u vedermi nella solitaria casa di Pavia.

Al principio di questo mese avevo fatto comunicar una mia lettera all'Ist.[ituto] Lombardo alla quale tenevo assai. Ma nella seduta ultima mi mancava ogni lena
per spiegare come avrei desiderato ai colleghi mat.[ematici] le cose mie. Ora poi lo sciopero degli operai tipografi impedisce la stampa dei Rendiconti. Perci\`o
avrei piacere di fare ai Lincei la Comunicazione suddetta, e per mezzo tuo, che sai sostenere cos\`{\i} bene i tuoi clienti.

La comunicazione ha per titolo: ``Il calcolo delle differenze finite interpretato ed accresciuto di nuovi teoremi a sussidio principalmente delle odierne ricerche
basate sulla variabilit\`a complessa".

In virt\`u della medesima moltissime ricerche che costarono gran fatica a matematici valenti diventano ovvie traduzioni di quanto gi\`a da tempo \`e stato fatto
nel calcolo delle differenze finite. In questa prima Comunicazione dimostro la importanza della mia interpretazione facendone applicazione a due argomenti: allo
studio delle equazioni algebriche a coefficienti monodromi; ed a quello ora di moda delle eq.[uazioni] diff.[erenziali] lineari.

Potrai vedere nella 1$^{\rm a}$ applicazione come riducasi a poche linee una met\`a delle celebri e lunghe ``Recherches sur les fonctions alg.[\'ebriques] del
sig. Puiseux" e nella $2^{\rm a}$ come scendano immediatamente dalle note formole d'integrazione delle equazioni alle differenze finite lineari e coi
coeff.[icienti] costanti le propriet\`a e le espressioni degli integrali delle equaz.[ioni] differenziali lineari con coefficienti monodromi che costarono tanta
fatica al sig. Fuchs (1$^{\rm a}$ sua Memoria) e diedero da fare ai signori Thom\'e Frobenius Hamburger, Jurgens.

Quel teorema sul determ.[inante] $\Theta$ parvemi veramente bello e di grande uso. Utile pure il determ.[inante] $H$.

Ero molto soddisfatto di questa interpretazione che sembrami di poter qua\-li\-fi\-ca\-re come una scoperta analitica di grande fecondit\`a, e m'avviavo a farne
molte altre applicazioni, quando la malattia dell'Eugenia venne a spaventarmi ed istupidirmi.

Se non ti dispiace di fare per me questa comunicazione nella prossima seduta Lincea ti mander\`o il manoscritto, che in parte feci gi\`a leggere a Beltrami.

Non t'avrei dato questo disturbo se il ministero t'avesse trascinato gi\`a sino d'ora, come si cominciava a dire, nel suo vortice.

Mille affettuosi e rispettosi saluti alla sig. Elisa ed alle tue care figlie, ed a te una cordiale stretta di mano.

\begin{flushright}
Tuo F. C.
\end{flushright}

{\it Letter 7}: {\underline{Casorati to Bertrand}}

\begin{flushright}
Pavia, 3 aprile 80
\end{flushright}

Car.$^{\rm mo}$ ed. illustr. Signore

\medskip

(......) Mi sarei deciso di ritirare il manoscritto destinato agli Annali, ampliarlo e riscriverlo meglio che potr\`o, e spedirlo al Segretario della vostra
Accad.[emia] entro il maggio, conformemente alle riflessioni della vostra ultima lettera. Ma devo pregarvi di una risposta alla seguente domanda: Posso mandare un
manoscritto in lingua italiana, o devo tradurlo in francese? Se aveste la bont\`a di darmi questa risposta, ve ne sar\`o sempre pi\`u obbligato, e spero che
sar\`a l'ultimo disturbo che vi avr\`o arrecato.

Per\`o non posso chiudere la lettera senza pregarvi vivamente, per quando riceverete il mio articolo stampato dai Lincei, di volergli dare un'occhiata. Voi
troverete tutto facilissimo, e mi \`e caro sperare che il teorema pel determin.[ante] $\Theta$, e le applicazioni che ho indicate alla fine, e la evidente
molteplicit\`a delle altre applicazioni possibili, non abbiano a parervi indegne della vostra attenzione.

\begin{flushright}
Vostro devotissimo \\
F. C. \\
Univ.[ersit\`a] di Pavia
\end{flushright}

{\it Letter 8}: {\underline{Bertrand to Casorati}}

\medskip

\begin{flushright}
Paris, 7 avril 1880 \\
Cher monsieur,
\end{flushright}

Il n'y a pas de temps perdu, aucun m\'emoire n'a \'et\'e envoy\'e jousqu'ici et la Commission qui doit les juger n'est pas encore nomm\'ee.

Un manuscrit en langue italienne sara certainement re\c{c}u et la lecture, j'en suis certains, sera facile \`a tous les commissaires. Je pense cependant que si
cela ne vous est trop p\`enible, il vaut mieux \'ecrire en Francais; vous serait jug\'e plus facilement et sans porter atteinte ou secret qui, \'exig\'e pour les
noms des concurrents, ne l'est pas pour leur nationalit\'e, mais reste cependent pr\'eferible. [...]

\begin{flushright}
Votre amis tr\`es d\'evou\'e, \\
J. Bertrand
\end{flushright}

{\it Letter 9}: {\underline{Casorati in the {\it billet cachet\'e}}}

\medskip

\begin{flushright}
ce 22 mai 1880
\end{flushright}
F\'elix Casorati\\
Professeur \`a l'Universit\'e de Pavie-Italie

\bigskip

Le manuscrit concernant les interpr\'etations du {\it Calcul aux diff\'erences} dont j'ai fait mention dans la pr\'eface a \'et\'e comuniqu\'e \`a l'{\it Institut
Lombard} {\it de sciences et lettres} dans sa s\'eance du 19 f\'evrier, et devait \^{e}tre publi\'e tout-de-suite dans les {\it Annali di Matematica (Milano)}.
Mais une longue tr\`eve des ouvri\`ers interrompit tout travail dans l'imprimerie des Annali. Alors j'ai comuniqu\'e une copie de ce manuscrit \`a l'{\it
Acad\'emie R. des Lincei} dans sa s\'eance du 7 mars. Mais l'exc\`es des travaux acad\'emiques \`a imprimer et une grave malheur domestique, qui d\'etournait mon
esprit des \'etudes, retard\`erent aussi l'impression de cette comunication. Les deux manuscrits, maintenant quelque peu amplifi\'es, vont enfin para\^{\i}tre
publi\'es dans peu de jours. Mais leur lecture serait tout-\`a-fait inutile pour qui conna\^{\i}t le manuscrit pr\'esent\'e avec ce billet au concours.

\medskip

{\it Letter 10}: {\underline{Casorati to Brioschi}}

\begin{flushright}
Madesimo presso Pianazzo \\
sullo Spluga, 31 luglio 80
\end{flushright}

Carissimo sig. Direttore (Brioschi)

\medskip

Gi\`a le dissi in Milano di aver fatto qualche pensiero sul premio Bressa da conferirsi stavolta all'Italiano che nel quadriennio 1877-80 a giudizio
dell'Accad.\-[emia] di Torino avr\`a fatto il miglior lavoro ``sulle scienze fisiche e sperimentali, storia naturale, matematiche pure ed appl.[icate], chimica,
fisiologia e patologia, non escluse la geologia, la storia, la geografia e la statistica".

In  verit\`a, per la moltitudine delle materie ammesse a concorso, \`e assai poco probabile che il premio venga dato ad un matematico. Per\`o, dappoich\'e ci\`o
\`e possibile, e dappoich\'e a me sembra che la Memoria, che ora le dirigo (Il calcolo delle diff.[erenze] ecc.), sia la pi\`u importante tra le apparse nel
quadriennio, in primo luogo per la estesa applicabilit\`a delle idee che vi sono esposte, in secondo luogo per la importanza dei teoremi (specialmente quello del
$\Theta$ o $\Delta$) e delle applicazioni che gi\`a in questo primo saggio potei abbastanza minutamente indicare, penso che sarebbe quasi colpevole negligenza in
me che ho famiglia, il nulla fare per mettermi in vista all'Accademia. Sentito il parere di Beltrami, ne parlai con Genocchi, il quale aggrad\`{\i} moltissimo
l'esposizione che gli feci di una parte della Memoria, allora non per anco uscita dalla Tipografia.  Ma a portare vie pi\`u l'attenzione dell'Accademia sulla
medesima, servirebbe egregiamente la proposta di Lei, quale socio dell'Accad.[emia], come pure quella di Betti. E pertanto, lo scopo della presente lettera \`e
pregarla a voler dare un'occhiata alla Memoria, ed a voler poi, nel caso che ne ricevesse impressione abbastanza favorevole, scrivere due semplici righe
all'Accademia, come gi\`a altri fecero pel 1$^{\rm o}$ premio (...)

Ed ora devo pregarla di scusarmi, se, per questa eccezionale circostanza dovetti lodare un mio lavoro. L'assicuro che questa necessit\`a mi spiace assai, e che
non penso di prendere da qui le mosse per seguire le consuetudini dei francesi. (...)

\begin{flushright}
il suo aff.$^{\rm mo}$ F. Casorati
\end{flushright}

{\it Letter 11}: {\underline{Casorati to Picard}}

\begin{flushright}
Pavie, ce 11 d\'ecembre 1880 \\
\end{flushright}

Monsieur (E. Picard)

Vos travaux, et tout ce que m'a dit des vous M. Mittag-Leffler, en me visitant il y a quelque mois, et la nouvelle de v\^{o}tre prochain marriage avec une fille
du grand G\'{e}om\`{e}tre\footnote{Picard was going to marry Charles Hermite's daughter.}, qui est une des plus pures et plus hautes gloires de la France, tout
cela m'ispire une vive estime et sympathie et, par cons\'equent, le d\'esir d'\'etablir avec vous quelques rapports personnels. C'est pour cela que je me suis
d\'ecid\'e \`a vous envoyer tous les exemplaires qu'il m'est possible, des mes travaux de math\'em.[atique] pure, (deux paquets) avec cette lettre, en vous priant
d'avoir la bont\'e de les agr\'eer et de vous souvenir possiblement de moi, lorsque vous dispenserez les exemplaires de vos propres travaux.

Je vous enverrai dans une autre occasion quelques autre M\'emoires, dont une partie est d\'ej\`a imprim\'ee, contenants bon nombre d'applications du {\it Calcul
des diff\'erences finies}, lesquelles vous interesseront, j'\`espere, beaucoup. Dans pluiseurs d'entre elles joue un r\^{o}le tr\`es important et
tr\`es-\'el\'egant le det\'erminant
$$
\left|
\begin{array}{cccc}
f_1 & f_2 & ... & f_n \\
\Delta f_1 & \Delta f_2 & ... & \Delta f_n \\
 ...& ...& ...&...\\
 \Delta^{n-1} f_1 &  \Delta^{n-1} f_2 & ... & \Delta^{n-1} f_n
\end{array}
\right|
$$
form\'e avec les diff\'erences successives de pluiseurs fonctions $f_1,...,f_n$, d'une m\^{e}me variable $t$, le quel j'\'ecris aussi comme il suit
$$
\left|
\begin{array}{cccc}
f_1 & f_2 & ... & f_n \\
\theta f_1 & \theta f_2 & ... & \theta f_n \\
 ...& ...& ...&...\\
\theta^{n-1} f_1 &  \theta^{n-1} f_2 & ... & \theta^{n-1} f_n
\end{array}
\right|
$$
en \'ecrivant $\theta f$, $\theta^2 f$,... au lieu de $f(t+1)$, $f(t+2)$,.... De telle mani\`ere, on peut regarder $\theta$ comme un symbole d'op\'eration, ce qui
est, comme vous verrez, tr\`es utile.

En envisageant la variabilit\'e de $t$ de plusieurs point des vue diff\'erents et en relations avec d'autres variables qui jouent elles-m\^{e}mes le r\^{o}le de
variables ind\'{e}pendantes, on obtient plusieurs interpr\'etations diff\'erentes du symbole $\theta$, et autant de mani\`eres corr\'espondantes de traduire
toutes les formules ou propositions, qui ont \'et\'e acquis\'ez jusqu'ici dans l'analyse directe et inverse des diff\'erences finies, dans des r\'esultats
concernants d'autres branches de l'analyse math\'ematique. Le Calcul des diff\'erences finies si neglig\'e dans ces derni\`ers temps acquiert par l\`a une
importance tr\`es-grande et inattendue m\^{e}me pour l'analyse infinit\'{e}simale.

Mais je ne veux pas \^{e}tre trop indiscret, en vous entretenant plus longtemps sur mes travaux. Veuillez donc agr\'ee, je vous en prie de nouveau, les
s\'entiments d'\'estime et de sympathie qui m'ont fait \'ecrire cette lettre, et soyez si bon de me pas m'oublier tout-\`a-fait.

\begin{flushright}
Votre serviteur \\

F. Casorati \\

Prof. \`a l'Universit\'e de Pavie \\

Italie
\end{flushright}

Je vous serai tr\`es-oblig\'e si vous voudrez m'assurer par quelques mots d'avoir re\c{c}u les paquets.

\medskip

{\it Letter 12}: {\underline{Picard to Casorati}}

\begin{flushright}
Toulouse, ce 17 d\'ecembre 1880.
\end{flushright}

\begin{center}
Monsieur
\end{center}
J'ai \'et\'e tr\`es flatt\`e de votre aimable lettre et je serai tr\`es heureux d'\'etablir avec vous de relations amicales. Je vous remercie vivement des
exemplaires de vos travaux que vous m'avez envoy\'es; la plupart d'entre eux m'\'etaient d\'eja connus, et quant \`a votre trait\'e sur les fonctions d'une
variable complexe c'est un livre excellent et que j'ai pu appr\'ecier d'autant mieux que les sujets qui y sont trait\'es me sont extr\^{e}mements familiers. Tout
ce que vous me dites des diff\'erences finies et notamment du d\'eterminant des $\Delta$, m'int\'eresse extr\`emement, et j'attends avec impatience les
r\'esultats complets de vos \'etudes sur ces importantes questions.

Je vous enverrai tr\`es prochainement les m\'emoires ou notes qu'ai publi\'es et dont j'ai encore entre le mains de tirages et vous pouvez \^{e}tre assur\'e que
dans l'avenir je me ferai un devoir et un plaisir de vous r\'eserver toujours un exemplaire de ce que je pourrai faire. Ce sont surtout les questions de la
th\'eorie g\'en\'erale des fonctions de variables complexes qui m'ont beaucoup occup\'e depuis quelque temps. Les beaux r\'esultats de M. Weierstrass ont \'et\'e
le point de d\'epart des mes propres r\'echerches et vous trouverez dans mon m\'emoire sur les fonctions enti\`eres un th\'eor\`eme qui, si je ne me trompe, a
peut \^{e}tre quelque importance. Je viens de corriger en ce moment les \'epreuves d'un m\'emoire qui va prochainement para\^{\i}tre dans le Journal de
Crelle\footnote{It is the paper [Picard, 1881].}, et qui est relatif aux \'equations lin\'eaires \`a coefficients doublement p\'eriodiques; j'esp\`ere qu'il vous
int\'eressera. Ce sont les recherches de M. Hermite sur l'\'equation de Lam\'e qui m'ont jet\'e sur la voie des r\'esultats auxquels je suis parvenu.

En ce moment, monsieur, mon mariage est la grand action qui me pr\'eocupe; j'ai du revenir \`a Toulouse pour faire quelques cours \`a la Facult\'e et j'attends
avec impatience mon d\'epart pour Paris, qui aura lieu dans trois jours.

Soyez persuad\'e que je recevrai toujours vos lettres avec joie et je serai heureux de toutes les communications math\'ematiques que vous voudrez bien me faire.
Permettez de vous assurer de sentiments d'estime et de sympathie que je conserve pour vous.

\begin{flushright}
Emil Picard \\
professeur \`a la Facult\'e des Sciences \\
de Toulouse.
\end{flushright}

{\it Letter 13}: {\underline{Casorati to Bonnet}}

\begin{flushright}
Pavie, ce 6 f\'evrier 1881
\end{flushright}

Cher Monsieur et ami (Ossian Bonnet)

\medskip

Avec ces lignes je vous envoye la traduction fran\c{c}aise d'une lettre que j'ai fait ins\'erer dans les Annali de M. Brioschi, \`a cause d'une Note publi\'ee
r\'ecemment par M. Stickelberger\footnote{This French version of [Casorati, 1880c] is [Casorati, 1881b].}; et je profite de l'amiti\'e dont vous m'avez toujours
honor\'e pour vous prier de vouloir l'observer, et de compatir mon irritation, en \'egard \`a la peine que j'ai d\^{u} essayer en voyant malveillance et,
peut-\^{e}tre, mauvoise foi de la part d'un coll\`egue.

Mon m\^{a}itre bien-aim\'e M. Brioschi me conseillez d'adresser \`a M. Hermite une Note, pour \^{e}tre ins\'er\'ee aux Comptes Rendus, dans le but de faire
toujours mieux constater ma priorit\'e dans l'application du Calcul des diff\'erences \`a la th\'eorie des \'equations diff\'erentielles lin\'eaires. Cependant
dans cette Note je n'ai rien dit qui p\^{u}t avoir l'aspect de personnalit\'e, en me bornant \`a rappeller mon Memoire, et du rest completant la d\'etermination
des sous-groupes d'int\'egrales faite par M. Hamburger; ce qui peut int\'eresser ind\'ependamment de toute question personnelle.

Mais dans une lettre confidentielle \`a un ami, quelque manifestation de m\'e\-con\-ten\-te\-ment personnel pourra \^{e}tre pardonne\'e.

L'ann\'ee derni\`ere je m'\'etais occup\'e avec ardeur de la theorie des \'equations diff\'erentielles lin\'eaires, mais une maladie tr\`es-grave de ma fille
ain\'ee, \'etant survenue [\`a] troubler ma t\^{e}te et deranger mes plans, je n'ai pu ne m'a permis de rediger mes recherches qu'en partie d'une mani\`ere
convenable pour \^{e}tre imprim\'ees. Comme si cela n'etait pas assez pour me f\^{a}cher, M. Stickelberger vient maintenent essayer d'approprier \`a soi m\^{e}me
et \`a Riemann et a M. Hamburger une partie des id\'ees qu'il m'a \'et\'e cependant possible de publier dans cette ann\'ee malheureuse. Mais sa tentative ne peut
qu'\'echouer. Je n'ai pas la pretension que personne ne s'occupe de l'application du Calcul des diff\'erences \`a la dite th\'eorie ou \`a d'autres branches de
l'analyse de la variabilit\'e continue. Au contraire je dois voir de bon gr\'e para\^{\i}tre des travaux fond\'e sur les id\'ees que j'ai con\c{c}ues et chercher
\`a expliquer. Mais je ne puis pas \^{e}tre content que l'on cherche \`a m'\^{o}ter ce qui m'appartient.

En redigeant sa Note avec calme, sans l'id\'ee injuste de la faire passer pour contemporaine de mon M\'emoire, il aurait pu la livrer au public sous une forme
beaucoup meilleure et il n'aurait pas certainement suscit\'e de r\'ecriminations de ma parte. La m\'ethode de Cauchy \`a la r\'esolution du syst\`eme d'equations
aux diff\'erences, qui est l'object de la premi\`ere partie du \S 1 est bien faite; et la combination des formules de r\'esolution ainsi obtenues avec la formule
de Weierstrass, effectu\'e dans la premi\`ere partie du \S 3, vient \`a-propos; quoique ce tour ne soit necessaire, comm'il semble croire, pour arriver \`a la
forme la plus simple de la distinction des integrales en sous-groupes. Mais, dans l'empressement de faire para\^{\i}tre sa Note, il ne s'est pas aper\c{c}u que le
seconde partie de ces deux paragraphes et les paragraphes restants (hormis le dernier qui se rapport \`a un autre point) sont trop remplis de r\'ep\'etitions et
des calculs superflus. Le \S 2 particulierment, tout plein de calculs assez ennuyeuse, pouvait se remplacer par quelques observations br\`eves, claires et sans
calculs.

La personne \`a laquelle je fais allusion, au commencement de ma lettre imprim\'ee, est M. Frobenius, qui a travaill\'e et publi\'e plusieurs articles en
communion avec M. Stickelberger. C'est pour cela que j'ai termin\'e cette lettre en rappellant mon chap\^{\i}tre sur le derteminant
$$
\Theta=\left|\begin{array}{ccc}
y_1 & ... & y_n \\
... & ... & ... \\
\theta^{n-1}y_1 & ... & \theta^{n-1}y_n
\end{array}\right| \qquad\mbox{ou son \'egal}\qquad \Delta=\left|\begin{array}{ccc}
y_1 & ... & y_n \\
... & ... & ... \\
\Delta^{n-1}y_1 & ... & \Delta^{n-1}y_n
\end{array}\right|.
$$
M. Frobenius tient fort beaucoup \`a ce qu'on le regarde comme l'auteur qui a le plus soign\'e le determinant
$$
D=\left|\begin{array}{ccc}
y_1 & ... & y_n \\
... & ... & ... \\
D^{n-1}y_1 & ... & D^{n-1}y_n
\end{array}\right|
$$
et M. Stickelberger d\'eclare pompeusement dans son dernier \S\ qu'il d\'esignera, d'apres M. Frobenius, par $D(y_1, y_2,...,y_n)$. N'est ce donc pas \'etonnant
qu'il puisse oublier mon chap\^{\i}tre susdit, qui comprend, comme cas particulier, tout ce que l'on peut dire sur le determinant $D$? Mais le temps fera justice,
car les applications nouvelles, qu'il m'arrive de faire chaque fois que je reviens \`a ce determinant $\Theta$, ne me laissent plus douter qu'il soit pour prender
une place distingu\'ee parmi les instruments les plus utiles et les plus \'el\'egants de l'analyse. Permettez-moi d'ajouter, qui je partage la r\'epugnance de la
plus part des g\'eom\`etre contre l'utilisation de notations ou d\'enominations nouvelles. Mais l'introduction d'un symbole (tel que $\theta$), pour d\'esigner
par $\theta f$ ce que une fonction $f$ devient par effet de la variation d'une nature prefix\'ee, de certaines \'el\'ements dont elle d\'epend, m'a paru trop
utile pour y renoncer. Parmi les advantages ce n'est pas petit celui-ci, de pouvoir comprendre sous un m\^{e}me enonc\'e des propositions qui se rapportent \`a
des sph\`eres des recherches tr\`es diff\'erentes. Dans les unes, la variation, suppos\'ee par $\theta$ pourra consister dans le tour qu'une variable ou plusieurs
variables imaginaires, dont nos fonctions d\'ependront, doivent faire autour d'un ou plusieurs syst\`emes des points singuliers; dans d'autres, la variation
pourra consister dans un syst\`eme d'accroissement que les \'el\'ements variables doivent prendre dans d'autres, enfin, la variation pourra consister dans tout
autre chose. Et il va sans dire qu'il est tr\`es utile, de pouvoir aussi employer plusieurs symboles $\theta$, $\theta'$,.. simultan\'ement.

Je voudrais bien vous entretenir un peu plus sur cet argument, mais je crois d'avoir absorb\'e d\'eja trop de temps \`a vos occupations. Et la note de M.
Stickelberger aurez vous eu le temps et l'envie de la lire? Dans le cas n\'egatif je serais encore plus coupable d'indiscretion! Mais votre bont\'e a \'et\'e
toujours tr\`es grande; j'\`espere donc, m\^{e}me dans ce cas, votre pardon. Et plus ancore, je compte aussi sur votre bienveillance pour l'avenir, de sorte que
pourrai toujours me sous\'egner

\begin{flushright}
votre ami d\'evou\'e \\
F. Casorati
\end{flushright}

\begin{flushleft}
Soyez aussi indulgent pour les \\
fautes que j'aurais commis\'es \\
en ecrivant en francais
\end{flushleft}

\medskip

{\it Letter 14}: {\underline{Bertrand to Casorati}}

\begin{flushright}
Paris 15 f\'evrier 1881
\end{flushright}

Je vous remercie bien cordialement, tr\`es cher Monsieur, de la nouvelle, pour moi si flattante, que vous voulez bien m'annoncer. Je ne suis pas digne de succeder
\`a mon illustre ami M$^{\rm r.}$ Chalses et n'en dois que plus de recoinnaissance \`a ceux qui ont bien voulu songer \`a moi.\footnote{Bertrand refers to his
election to the {\it Istituto Lombardo} on Febraury 10th, 1881.} Vous voudrez bien leur transmettre l'expression de mes sentiments les plus sinc\`erement et
enti\`erement devou\'es.

L'Acad\'emie a d\'ecern\'e hier le Prix de Math\'ematiques pour la question des \'equations lin\'eaires. Un tr\`es beau et ingenieux m\'emoire dans laquel
l'\'etude des point critiques est \'eclairci et simplifi\'ee par l'emploi des \'equations aux diff\'erences, avait attir\'e l'attention de la Commission qui n'a
cependant lui accorder qu'une mention tr\`es honorable. Que le savant auteur veuille bien attendre pour condamner ses juges, la lecture du m\'emoire tr\`es
consid\'erable qui a \'et\'e couronn\'e et dont l'auteur est M$^{\rm r}$ Halphen; notre compatriote s'y montre beaucoup plus grand g\'eom\`etre encore que les
excellents debuts ne l'avaient promis. Il \'etudie les changements d'une \'equation lin\'eaire quand on remplace la variable independante par une autre ou quand
on multiplie l'inconnue par une fonction choisie de mani\`ere \`a simplifier l'\'equation qui ne cesse pas d'\^{e}tre lineaire. Il me semble qu'il \'equisse le
sujet en tirant un tr\`es grand partie des invariants introduits par M. Laguerre qui je crois n'a pas concours.

Veuillez recevoir, cher Monsieur Casorati, l'assurance de mes sentiments les plus affectueuse.

\begin{flushright}
J. Bertrand
\end{flushright}

{\it Letter 15}: {\underline{Casorati to Bertrand}}

\begin{flushright}
Pavie, ce 21 f\'evrier 1881
\end{flushright}

Cher Monsieur

je vous remercie de votre comunication, adress\'ee \`a Milan m'est parvenue en ritard, du 15 f\'evrier, et particulierement des expressions affectueuses dont vous
confortez mon esprit dans une circonstance o\`u il en a vraiment tout le besoin. Car, je l'avait dit avec Pascal dans ma devise, nous ne pouvons prendre plaisir
\`a une chose qu'\`a condition de nous f\^{a}cher si elle ne r\'eussit pas. Et je ne puis \`a pr\'esent m'emp\^{e}cher d'\^{e}tre f\^{a}ch\'e douloureusement
d'autant plus que j'ignore les terms pr\'ecis du rapport. Vos Comptes Rendus ne sont pas seulement la pubblication la plus observ\'ee par les professeurs de mon
Universit\'e, mais ils sont lus aussi par les \'etudiants qui sont inscrits dans l'\'ecole normale annex\'ee maintenent \`a cette Facult\'e des Sciences
math.[\'ematiques] etc.

Un professeur de cette m\^{e}me Universit\'e, feu M. Codazzi, qui avait concouru au prix relatif \`a la th\'eorie des surfaces applicables l'une sur l'autre, a pu
rester satisfait de la simple mention qu'on lui avait accord\'ee en Mars 1861. Mais alors le rapporteur, que vous connaissez tr\`es-bien, avait \'ecrit: "Les
trois autres M\'emoires, inscrits sous les n$^{\rm .os}$ 1,2, et 5, {\it remplissent compl\`etement le programme trac\'e par l'Acad\'emie. Si l'un quelconque des
trois avait \'et\'e pr\'e\-sen\-t\'e seul \`a notre examen, nous lui aurions sans h\'esiter accord\'e le prix}."\footnote{In fact, Bertand had written this
report.}

Maintenent, qui sera le rapporteur?

Et conna\^{\i}tra-t-il assez une position scolastique pour prendre bon soin de la m\'enager par quelques expression de la valeur de celles que je vient de citer?
Ce serait pour lui d'autant plus facile que mon travail et celui de M. Halphen n'admettent pas, il me semble, une comparaison rigoreuse; leur points de d\'epart
et leur directions \'etant differents; et mon travail contienant aussi incontestablement des perfectionnements importants. Je ne doute pas de la haute valeur du
M\'emoire du jeune g\'eom\`etre que vous avez couronn\'e; mais son m\'erite n'exclut pas que mes applications du Calcul aux differences soient tout-\`a-fait
nouvelles, et que leur port\'ee ne soit n\'ee \`a l'\'etude des int\'egrales autour des points critiques, faite par M. Fuchs, comme on peut per\c{c}evoir par mon
manuscrit et comme on verra toujours mieux dans la suite.

A la v\'erit\'e l'expression {\it tr\`es honorable} contenue dans votre lettre devrait me rassurer; mais je vous prie de compatir ces craintes, dont je ne sais
pas me delivrer sur le moment. Du reste, je vous prie aussi de croire que je ne voudrai pas mal \`a mes juges, ni \`a mon ma\^{\i}tre bien aim\'e, M. Brioschi,
qui m'a inspir\'e l'id\'ee du concours. Je ne ferai que r\'ep\'eter ce que j'ai pens\'e dans d'autres moments de ma vie: d'\^{e}tre n\'e sous un \'etoile assez
malheureuse.

Pardonn\'ez moi, cher Monsieur Bertrand, cette lettre trop longue et trop triste et veuillez me conserver votre bienveuillance toujours pr\'ecieuse.

\begin{flushright}
F. Casorati
\end{flushright}

{\it Letter 16}: {\underline{Stickelberger to Brioschi.}}\footnote{In the {\it Nachlass} there is a copy written by Casorati's daughter, Eugenia.}

\begin{flushright}
Schaffhause [sic!], 9 avril 1881
\end{flushright}

Monsieur le directeur

Puisque vous avez accord\'e quelque pages de vos annales \`a la lettre de M. Casorati, ecrite \`a l'occasion de mon m\'emoire ``Zur Theorie der linearen
Differentialgleichungen", j'\'esp\`ere que vous voudrez bien ins\'er\'e dans le m\^{e}me recueil une response de ma part, qui, du reste, ne sera pas longue.

M. Casorati part de la supposition que j'aie commenc\'e mon travail ayant conaissance de son m\'emoire; il me reproche d'avoir employ\'e sa m\'ethode sans en
indiquer l'origine. Voici ce que je pourrais lui r\'epondre. C'est seulement apr\`es avoir trouv\'e les principaux r\'esultats des paragraphes 1 \`a 4 de mon
\'ecrit, que j'en parlai avec M. Lindemann en faisant mention sp\'eciale du calcul aux differences finies, et c'etait pr\'ecisement cette mention qui l'engagea
\`a me faire parvenir le m\'emoire de M. Casorati que celui-ci venait de lui envoyer. Mais je crois ne pas avoir besoin de repondre au [sic!] long \`a une
imputation qui ne trouvera aucun credit aupr\`es  de tout ceux qui me connaissent.

Au reste il est evident que M. Casorati n'a pas pris la peine de lire mon m\'emoire en entier avant d'\'ecrire sa lettre. S'il avait fait, il aurait remarqu\'e
que ce qui, pour lui, est une id\'ee cardinale, n'est pour moi qu'un moyen pour arriver \`a un but determin\'e; il ne m'aurait pas demand\'e, pourquoi j'ai
pass\'e sous silence les recherches contenues dans son chap\^{\i}tre deuxi\`eme et, en particulier, le th\'eor\`eme relatif \`a la condition pour l'existence
d'une \'equation lin\'eaire \`a coefficients p\'eriodiques entre plusieurs fonctions d'une seule variable ind\'ependante. D'ailleurs, s'il avait en lieu
d'employer ce th\'eor\`eme, je ne l'aurais pas emprunt\'e \`a M. Casorati, mais \`a un M\'emoire de M. Christoffel ``\"Uber die lineare Abh\"angigkeit von
Functionen einer einzigen Ver\"andlichen" publi\'e en 1858 dans le Tome 55$^{me}$ du journal de Borchardt, o\`u ce suject est trait\'e d'une mani\`ere aussi
elegante que compl\`ete. J'ai coutume de citer, dans toutes les questions importants, les auteurs qui, \`a ma connaissance, ont propos\'e les pr\`emiers les
th\'eor\`emes dont je fais usage.

M. Casorati m'accuse d'avoir manqu\'e de politesse vis-\`a-vis avec lui. Je laisse aux lecteurs des Annali le soin de juger, s'il est plus courtois d'employer,
dans une critique purement objective, les mots ``nicht rathsam", ou bien de se servir, dans une r\'esponse toute personnelle, d'\'expressions telles que
``disgusto, malvolere, cattivi pensieri".

Agr\'ez, Monsieur, l'expression de mon plus profond respect

\begin{flushright}
Votre d\'evu\'e \\
L. Stickelberger
\end{flushright}

{\it Letter 17}: {\underline{Casorati to Brioschi}}

\begin{flushright}
Pavia, 17 aprile 1881
\end{flushright}

Carissimo sig. Direttore (Brioschi)

Ho letto la Memoria di Christoffel del 1858 e vi trovai infatti il teorema che la eguaglianza
$$
\sum\pm f(m)f_1(m+1)\cdots f_n(m+n)=0
$$
\`e condizione necessaria e sufficiente affinch\'e le funzioni
$$
f(m), \quad f_1(m),.....,f_n(m)
$$
della variabile discreta $m$ abbiano tra loro una relazione lineare a coefficienti costanti. Se avessi conosciuto questa Memoria, non avrei mancato di citare il
Christoffel nell'occasione del mio determinante $\Theta$ come non mancher\`o di ricordarlo in una ventura occasione. Ma il Christoffel non allude mai ai {\it
coefficienti periodici} e molto meno a coefficienti che sieno {\it funzioni di una variabile complessa, monodrome intorno ad un punto}. Non pu\`o essere che il
{\it malvolere} dello Stickelberger capace di insinuare che queste idee fossero implicite in quella Memoria. (....)

\section{Appendix 2: Original letters concerning Pincherle's unpublished note on Tannery's theorem}\label{sec:originals_Tannery}

{\it Letter 1}: Pincherle to Casorati.

\begin{flushright}
Bologna, li 21/7/85.
\end{flushright}

Chiarissimo Sig. Professore

\medskip

Nella Memoria del Sig. Goursat (Sur une classe de fonctions represente\'es par des int\'egrales d\'efinies---Acta Math., 1883) ho trovato l'enunciato di una
proposizione data dal Sig. Tannery, che \`e la seguente:

{\it Se $n$ funzioni sono tali che facendo qualunque giro colla variabile intorno a $p$ punti $a_1$, $a_2$,...$a_p$, si ritrovano sempre combinazioni lineari a
coefficienti costanti delle $n$ funzioni, esse funzioni soddisfano ad un'equazione differenziale lineare omogenea a coefficienti monotropi rispetto ai punti
$a_1$, $a_2$,...$a_p$.}

Non ho potuto vedere quale sia la dimostrazione che il Sig. Tannery d\`a di questo teorema, essendo esso pubblicato nelle [sic!] {\it Annales de l'Ecole [Sic!]
Normale,} che non si trovano nella Biblioteca di Bologna.

Per\`o questa dimostrazione risulta come conseguenza immediata dal di Lei metodo per lo studio dei valori di una funzione nell'intorno di un punto singolare,
contenuto nella Sua memoria {\it Il calcolo delle differenze finite ecc.} (Annali di Matematica, t. X); e spero che non le sar\`a discaro di leggere questa
dimostrazione. Eccola in due parole:

Formiamo colle $n$ funzioni date, e con $n$ costanti arbitrarie, la funzione
$$
E= c_1E_1+c_2E_2+\cdots+c_nE_n;
$$
eseguendo un giro intorno ad $a_1$, ed essendo $\theta$ l'operazione da Lei introdotta nella citata Memoria, sar\`a
$$
\theta E =\sum_{i=1}^n c_i\theta E_i
$$
e per l'ipotesi fatta:
$$
\theta E =\sum_{i=1}^n k_{1i}E_i
$$
e cos{\`\i}
$$
\theta^r E=\sum_{i=1}^n k_{ri}E_i\qquad, (r=2,3,\dots,n):
$$
da queste eliminando le $E_1$, $E_2$,...$E_n$, viene:
$$
K_n\theta^n E+K_{n-1}\theta^{n-1}E+\cdots+K_1\theta E+K_0E=0,
$$
dove le $K_i$ sono costanti. Ma questa equazione (Casorati, mem.[oria] citata, \S 10) indica che la $E$ soddisfa ad un'equaz.[ione] differenziale lineare a
coefficienti monotropi in $a_1$, sia
$$
\vap_0E^{(n)}+\vap_1E^{(n-1)}+\cdots+\vap_{n-1}E'+\vap_nE=0;
$$
e di questa sar\`a $E$ l'integrale generale. Ma lo stesso ragionamento dimostrerebbe che la $E$ \`e l'integrale generale di un'equaz.[ione] a coefficienti
monotropi in $a_2$, di un'equaz.[ione] a coeff.[icienti] monotropi in $a_3$, ecc., il che non pu\`o essere che se
$$
\vap_0,\ \vap_1,\ \dots \vap_n
$$
sono monotropi in $a_1$, $a_2$,...$a_p$; c.d.d.

Io credo che questa dimostrazione cos{\`\i} semplice di una proposizione cos{\`\i} importante non sar\`a senza interesse. In ogni modo La lascio giudice di
decidere se sia il caso di farne una comunicazione a qualche Accademia; a me basta che risulti una volta di pi\`u l'importanza di quella Sua bella Memoria.

Mi conservi la Sua preziosa benevolenza, e mi creda, ch$^{\underline{\rm mo}}$ Sig. Professore

\begin{flushright}
di Lei dev$^{\underline{\rm mo}}$ \\
S. Pincherle
\end{flushright}

\medskip

{\it Letter 2}: Casorati to Pincherle.

\medskip

\begin{flushright}
H\^otel Piora sopra Airolo (Svizzera) \\
30 luglio 1885
\end{flushright}

Caro Pincherle

\bigskip

Ricevo adesso da Pavia la sua lettera del 21 e rispondo come posso.

Sono assai distratto dagli stud{\^\i} e lontano dai libri, e non ebbi occasione di vedere la Mem.[oria] di Tannery contenente la proposizione enunciata nella
Mem.[oria] del Goursat. Nondimeno convengo con Lei che tale proposizione non si possa dimostrare in maniera pi\`u semplice dell'indicata nella sua lettera.
Naturalmente, a me fa piacere di veder presente in Lei la mia Mem.[oria] ``Il calc.[olo] delle differ.[enze] ecc." e farebbe pur piacere di vederla per di Lei
mezzo ricordata agli studiosi che credo potrebbero farne facilmente molte applicazioni. Ma sem\-bre\-reb\-be\-mi conveniente ch'Ella desse un'occhiata a quella
Mem.[oria] Se fossi a Pavia, gliela manderei; ma cos{\`\i} bisognerebbe ch'Ella dimandasse il vol.[ume] degli Annali de l'Ecole [sic!] normale a qualche collega
di Milano, che potrebbe ritirarlo da quel collegio degli Ingegneri. Tale occhiata potrebbe suggerirle anche altre applicazioni della stessa mia Mem.[oria] da
presentarsi insieme con la gi\`a scrittami ai Lincei o ad altra Accad.[emia]

Quando avr\`o il piacere di vederLa, Le dir\`o dei non pochi progetti che feci anch'io di correggere, semplificare, ricordare cose nostre [ai] matem.[atici]
fran\-ce\-si in questi ultimi anni.

Ma Ella, ch'\`e giovane, fa bene a non contentarsi del semplice {\it fare progetti}.

Se non potr\`a avere il Tannery da Milano, me lo scriva, che glielo mander\`o poi io da Pavia.

Continui ad amar
\begin{center}
l'aff.[ezionatissimo] suo F.C.
\end{center}

{\it Letter 3}: Pincherle to Casorati.

\medskip

\begin{flushright}
Bologna, li 3 Ottobre 1885
\end{flushright}

\begin{center}
Ch$^{\rm \underline{mo}}$ Sig. Professore
\end{center}

Secondo il suo Consiglio, ho aspettato a redigere la Nota sulle equazioni differenziali, per cercare prima la Memoria del Tannery: ma non mi \`e stato possibile
avere quel volume degli {\it Annales de l'Ec.}, e non ho voluto disturbarla per questo finch\'e Ella si trovava presumibilmente in campagna.

Ora ella sar\`a forse tornato a Pavia: in questo caso Le sarei gratissimo se mi potesse far avere quel volume (Serie 2, Tomo 4) degli Annali; e dopo letta la
memoria del Tannery, rediger\`o (se sar\`a ancora del caso) e Le spedir\`o la noticina, rimettendo a Lei di decidere se convenga o no presentarla a qualche
Accademia.

In queste vacanze ho lavorato ad una Memoria sugl'integrali definiti atti a rappresentare funzioni analitiche, nella quale spero di dare alcuni risultati nuovi:
ma ci vorr\`a forse del tempo perch\'e sia in ordine per la stampa, e l'argomento mi si fa sempre pi\`u vasto. Mi consiglia Ella a sospendere per qualche tempo
ogni pubblicazione, lavorando intanto per il Premio dei Lincei dell'89, o \`e troppa presunzione la mia?\footnote{Pincherle shared this prize with Luigi
Bianchi.}

Sarei ben lieto che si presentasse quell'occasione, che Ella accenna nella sua lettera, di poterla vedere qui; \`e un desiderio che nutro da molti anni, perch\'e
da molti anni non ho potuto discorrere dei miei stud{\^\i} con persona competente e perch\'e voglio ringraziarla a voce della benevolenza di cui Ella mi \`e stato
prodigo in questo periodo di tempo.

E con questa speranza, La prego di credere all'affetto ed alla devozione del suo

\begin{flushright}
Salvatore Pincherle
\end{flushright}

\medskip

{\it Letter 4}: Casorati to Pincherle.

\begin{flushright}
Pavia, 17 ott. 1885
\end{flushright}

\begin{center}
Caro Pincherle
\end{center}

\noindent
Circa il premio linceo io non trovo niente affatto presuntuoso il progetto suo di concorrervi. Di pi\`u non vorrei dire; se non altro, perch\'e impedito dalla
lontananza di seco conversare distesamente.

\noindent Rileggendo la Sua lettera del 21 luglio, vedo che la proposizione ivi con\-si\-de\-ra\-ta appartiene all'antica (1874) [sic!] Memoria del Tannery sulle
equaz.[ioni] diff.[erenziali] lineari; Memoria che a suo tempo lessi con attenzione. E per\`o, pur seguitando a credere che la di Lei dimostrazione sia della
massima semplicit\`a non oserei consigliarle il sacrificio del tempo occorrente alla redazione di una Nota per la stampa, se nella Nota non dovesse entrare che la
pura dimostrazione suddetta. Faccia dunque Ella ci\`o che le pare conveniente; \`e certo ch'io ne rester\`o in ogni modo soddisfatto. Quanto al volume contenente
la Mem.[oria] del Tannery, glielo spedir\`o appena sar\`a arrivata una delle tre persone che tengono le chiavi della libreria della Scuola Normale di questa
Universit\`a.

Le vacanze hanno ristorato assai le mie forze fisiche; ma sono tuttora sba\-lor\-di\-to dall'inaspettata perdita di mio fratello\footnote{Luigi Casorati was a
jurist. He was born in Pavia on January 26th, 1834 and died in Rome on August 4th, 1885.}. Era il solo superstite della mia famiglia paterna, aveva per me
indicibile affetto ed io per lui tenerezza e ve\-ne\-ra\-zio\-ne ad un tempo.

Mi ricordi ai colleghi e mi voglia sempre bene.

\begin{flushright}
Suo F. Casorati
\end{flushright}

\medskip

{\it Letter 5}: Pincherle to Casorati.

\begin{flushright}
Bologna, li 19/10/85
\end{flushright}

Chiarissimo Sig. Professore

\medskip

Con vivo rincrescimento rilevo dalla Sua lettera la notizia della perdita del di Lei fratello: io non ho avuto l'onore di conoscerlo, ma ne ho spesso sentito
parlare con deferenza e stima. Le auguro che il tempo e l'affetto dei suoi possano lenire il Suo dolore, cui mi associo.

In questa circostanza non voglio pi\`u tediarLa con cose mie, che non possono avere per Lei che un meschino interesse, e ringraziandoLa dell'offerta di spedirmi
il volume degli {\it Annales}, credo che Ella possa risparmiarsi questo di\-stur\-bo se non vede l'opportunit\`a di pubblicare quella dimostrazione. Certamente
quella Memoria, essendo del 1874 (cosa che non sapevo) non pu\`o contenere risultati che non siano noti dopo i pi\`u recenti lavori sulle equazioni
differenziali.

Come gi\`a Le dissi, ho trovato quel Teorema citato dal Goursat, in una Memoria che studiavo per i miei lavori sugl'integrali; e mi venne subito alla mente l'idea
che la dimostrazione di quella proposizione si dovesse ottenere nel modo pi\`u semplice applicando il Suo metodo (l'operazione che Ella chiama $\theta$.)
D'altronde questa osservazione non ha attinenza diretta coll'argomento che ora sto studiando, e di cui Ella vedr\`a fra poco un primo saggio in una ``Note sur une
int\'egrale d\'efinie" che uscir\`a in uno dei prossimi fascicoli degli {\it Acta Ma\-the\-ma\-ti\-ca.}

Mi auguro che questo lavoro sia per incontrare la Sua approvazione, alla quale tengo sopra tutto; e pregandola di riverire per me i Sig. Professori Beltrami,
Bertini e Maggi\footnote{Eugenio Bertini was a distinguished algebraic geometer who taught in Pavia in the period 1880-1892. Bertini theorems are still of use in
current research in algebraic geometry. Gian-Antonio Maggi was Casorati's son in law and taught Rational mechanics and Mathematical physics in Modena, Messina,
and Milan. A set of equations to study non-holonomic systems is named after him.}, La prego, ch$^{\underline{\rm mo}}$ signor Professore, di tenermi sempre per

\begin{flushright}
Suo dev$^{\underline{\rm mo}}$ e aff$^{\underline{\rm mo}}$\\
S. Pincherle
\end{flushright}

\medskip

{\it Letter 6}: Casorati to Pincherle.

\begin{flushright}
24 dic. 85
\end{flushright}

\begin{center}
Mio caro Pincherle
\end{center}

Voglia aggrad.[ire] coi miei aug.[uri] l'artic.[olo] che le invio.

Avrei voluto scriverle assai prima per dirle che ha dato repentinamente troppo peso alla mia opin.[ione] circa quella tale dimostr.[azione] Del resto, io non
potevo vederne che con piacere la pubbl.[icazione], che Ella ad ogni modo avrebbe potuto fare dicendo che proponevasi ricordare agli studiosi una Memoria forse
meritevole di qualche attenzione (e non di essere negletta, come si fa dai Francesi, non ostante la lunga recensione fattane da non so chi nel Bulletin di
Darboux\footnote{The review appeared in 1882 [Darboux, Hou\"el, Tannery, 1882].}).

\begin{flushright}
Ami sempre il suo aff. F.C.
\end{flushright}

\medskip

{\it Letter 7}: Pincherle to Casorati.

\begin{flushright}
Bologna, li 29 X$^{\rm bre}$ 85
\end{flushright}

Chiarissimo Sig. Professore.

\medskip

Prima di tutto la Prego di accogliere i miei pi\`u sinceri auguri di felicit\`a per il prossimo anno: ed i miei ringraziamenti per la Memoria che Ella mi ha
favorito\footnote{Pincherle presumably refers here to the preliminary version of [Casorati, 1886a, b] that was printed in Milan in 1885.}.

Questa Memoria \`e destinata certamente ad avere una grande importanza, e tocca uno degli argomenti pu\`i vitali della Analisi. Anche a me---per quanto poco possa
valere la mia opinione,---era sempre sembrato naturale che si dovesse poter trovare propriet\`a della $\int_0^z\frac{\dd z}{\sqrt{R(z)}}$ o della sua inversa
anche se $R(z)$ \`e di grado superiore al $4^{\underline{\rm to}}$; e credo che se ci\`o non \`e stato fatto si deve attribuire in parte alla grande fecondit\`a
delle idee di Abel e Jacobi, di considerare le inverse come funzioni di pi\`u variabili, che non ha fatto cercare una estensione della Teoria delle f.[unzioni]
Ellittiche nel campo di una variabile sola. Sono persuaso che in questo ordine d'idee Ella giunger\`a  a dei risultati di somma importanza. E cos{\`\i}, la Teoria
delle f.[unzioni] Ellittiche verr\`a anch'essa ad acquistare maggiore importanza, trovandosi per cos{\`\i} dire all'{\it intersezione} delle trascendenti Abeliane
da una parte, delle Sue nuove trascendenti dall'altra.

Io avevo rinunziato all'idea di pubblicare quella Nota sulla Mem.[oria] di Tannery, pi\`u di tutto per non darLe il disturbo di spedirmi il volume degli Annales
de l'Ec.[ole] Normale. Se Ella per\`o ha la compiacenza di farmelo avere, vedr\`o di stendere quella Nota, cercando anche d'applicare il Suo metodo ad altri
Teoremi se \`e possibile, e quando la Nota sia redatta, gliela spedir\`o ed Ella giudicher\`a se sia o no da pubblicare.

Desidererei pure dalla Sua Cortesia un'altro [sic!] favore, che spero Ella non mi vorr\`a negare; e sarebbe (con tutto suo comodo) l'indicazione del tomo del
American Journal nel quale trovansi le memorie d'un certo Daniell sulla Teoria delle f.[unzioni] Ellittiche secondo il Weierstrass.\footnote{The author's name is
misprinted: the notes on Weierstrass theory of elliptic functions were written by A.L. Daniels [Daniels, 1883a,b, 1884].} Se potessi avere quella indicazione,
farei venire quel volume per la via di questa Biblioteca, dalla Biblioteca Universitaria di Pavia, dove si trova a quanto credo, l'American Journal.

Ho avuto il piacere di vedere poco fa il prof. Maggi, che mi ha dato le di Lei notizie e mi ha fatto sperare una Sua gita a Modena e a Bologna. Sarei lietissimo
di questa circostanza, che mi auguro gi\`a da 5 anni.

Pregandola d'accogliere di nuovo i miei auguri ed i miei ringraziamenti, e di riverire per me il Prof. Beltrami colla Sua Signora, ed i di Lei colleghi di
cost{\`\i}

\begin{flushright}
mi creda, Ch$^{\underline{\rm mo}}$ Signor \\
Professore, di Lei dev$^{\underline{\rm mo}}$ \\

S. Pincherle
\end{flushright}

\medskip

{\it Letter 8}: Pincherle to Casorati.

\begin{flushright}
Bologna, li 9 Gennajo 1886
\end{flushright}

\begin{center}
Ch$^{\rm\underline{mo}}$ Signor Professore
\end{center}

Le accludo la nota relativa al teorema di Tannery, lasciandola giudice dell'op\-por\-tunit\`a di presentarla all'Istituto Lombardo. L'enunciato che ho dato \`e un
po' pi\`u generale di quello di Tannery, ed i coefficienti dell'equaz.[ione] differenziale possono anche essere f.[unzioni] polidrome di specie determinata
(f.[unzione] uniforme d'un punto analitico) riguardando il campo $T$ come una superficie di Riemann ridotta, con tagli; semplicemente connessa.

Le spedisco pure, raccomandato, il volume degli {\it Annali}: la Memoria di Tannery, che ho letto in questa occasione, non mi \`e sembrata nulla di straordinario,
una compilazione dei lavori di Fuchs fatta discretamente e nulla pi\`u. Non vi \`e di notevole che il Teorema che mi ha ispirato questa Nota.

La ringrazio nuovamente per il disturbo che Ella s'\`e preso di spedirmi il libro, e La prego di credermi sempre

\begin{flushright}
il Suo dev.$^{\rm\underline{mo}}$ ed obblig$^{\rm\underline{mo}}$ \\
S. Pincherle
\end{flushright}

\bigskip
\bigskip

\noindent Analisi: Sopra un teorema del Sig. Tannery

\bigskip
\bigskip

Scopo di questa breve Nota \`e di mostrare con una nuova applicazione, i vantaggi che pu\`o arrecare nello studio delle funzioni analitiche e del loro modo di
comportarsi nell'intorno di punti determinati, il metodo ideato dal prof. Casorati ed esposto nella Sua Memoria: ``Il calcolo delle differenze finite
interpretato\footnote{V. un estratto particolareggiato di questa Memoria nel ``Bulletin de Darboux, S. II, T. VI, 1882''} ecc., Annali di Matematica, S.II, T.
X.'' Ne faccio qui l'applicazione alla dimostrazione di un teorema enunciato dal Sig. Tannery nella Memoria ``Propri\'et\'es des int\'egrales des \'equations
diff\'erentielles lin\'eaires, Annales ce l'Ec. Normale, S. II, T. IV, p. 130;'' teorema notevole per s\'e, e per una applicazione importante che ne ha fatto il
Sig. Goursat.\footnote{Nella Memoria ``Sur une classe des fonctions represent\'es per des int\'egrales d\'efinies. Acta Mathematica, T.II''.}

Il teorema \`e il seguente, sotto un enunciato un po' modificato: ``Siano $n$ funzioni analitiche $E_1$, $E_2$,....$E_n$ a carattere regolare nell'intorno di ogni
punto di un campo $T$ semplicemente connesso, eccettuati i punti $a_1$, $a_2$,...$a_p$ in numero finito e fra le quali non passa alcuna relazione lineare  a
coefficienti costanti; se quando la variabile gira senza uscire dal campo $T$ intorno ad uno qualunque dei punti $a_1$, $a_2$,..$a_p$, i nuovi valori delle
funzioni sono legati ai primitivi da relazioni lineari a coefficienti costanti, queste funzioni sono gl'integrali di un'equazione differenziale lineare a
coefficienti monodromi in $T$.''

\bigskip

\noindent Infatti, formiamo la funzione $E$ i cui coefficienti $c$ siano costanti arbitrarie
\begin{equation}\tag{1}
E=c_1E_1+c_2E_2+\dots +c_nE_n ;
\end{equation}
eseguendo colla variabile un giro intorno ad $a_1$, ed essendo $\theta$ l'operazione di Casorati, viene
$$
\theta E=c_1\theta E_1+c_2\theta E_2+\dots+c_n\theta E_n
$$
e per l'ipotesi fatta
\begin{equation}\tag{2}
\theta E=k_{1,1}E_1+k_{1,2}E_2+\dots +k_{1,n}E_n;
\end{equation}
analogamente
\begin{equation}\tag{3}
\left\{
\begin{array}{l}
\theta^2 E=k_{2,1} E_1+k_{2,2} E_2+\dots+k_{2,n} E_n,\\
\dots\dots\dots \\
\theta^n E=k_{n,1}E_1+k_{n,2} E_2+\dots+k_{n,n} E_n;
\end{array}
\right.
\end{equation}
eliminando $1, E_1$, $E_2$,,...$E_n$ fra le (1), (2) e (3), viene una equazione
$$
K_n\theta^nE+K_{n-1}\theta^{n-1}E+\dots K_0E=0;
$$
ma questa (Casorati, loc. cit. \S 10) indica appunto che la $E$ soddisfa ad un'equazione dif\-fe\-ren\-zia\-le lineare a coefficienti monodromi nell'intorno di
$a_1$
$$
\varphi_0 E^{(n)}+\varphi_1 E^{(n-1)}+\dots+\varphi_n E=0,
$$
di cui $E$ \`e l'integrale generale. La stessa dimostrazione far\`a conoscere che $\varphi_0$, $\varphi_1$,...$\varphi_n$ sono monodromi nell'intorno di $a_2$,
$a_3$,... $a_p$; onde il teorema \`e dimostrato.

Se il campo $T$ ricopre tutta la sfera una sol [sic!] volta, le funzioni $\varphi_0$, $\varphi_1$,...$\varphi_n$ sono uniformi. Non \`e escluso che il campo $T$
possa essere una superficie di Riemann ridotta semplicemente connessa coi tagli opportuni; nel qual caso le $\varphi_0$, $\varphi_1$,...$\varphi_n$ sarebbero
funzioni uniforme [sic!] d'un punto analitico.

\medskip

\noindent S. Pincherle

\bigskip

{\it Letter 9}: Pincherle to Casorati.

\begin{flushright}
Bologna, li 27/2/86
\end{flushright}

Chiarissimo Signor Professore

\medskip

Mi permetto di manifestarle di nuovo il mio compiacimento per aver potuto passare qualche ora nella Sua compagnia: e sono ben grato al prof. Maggi di avermi
procurato col suo gentile invito l'occasione di quella visita che ho forse prolungata oltre i limiti della discretezza, ma che a me \`e sembrata troppo breve.
Spero che nell'entrante Marzo Ella far\`a, come ci ha fatto sperare, la sua gita a Roma passando per Bologna e mi procurer\`a cos{\`\i} il piacere di rivederLa.

Ho letto attentamente il paragrafo della Memoria di Tannery che Ella mi ha indicato, ma non ho potuto trovarvi una inesattezza. La dimostrazione \`e per\`o
redatta poco bene, e vi \`e un errore (di stampa probabilmente) alla linea 10 (p. 133) dove va letto $m$ invece di $m-1$.

Le rimando il volume degli {\it Annales}; ho copiato il paragrafo in discorso per cui, anche senza tenere pi\`u a lungo il volume presso di me, potr\`o ripensare
meglio a quella dimostrazione: la quale si pu\`o forse redigere come ho tentato nel poscritto seguente.

I professori Boschi\footnote{Pietro Boschi (Rome, 7/6/1833-Bologna, 4/11/1887) taught projective and descriptive geometry in Bologna.}, Villari\footnote{Emilio
Villari (Naples, 25/9/1836-20/8/1904) was a chemist and a physicist who taught experimental physics in Bologna from 1871 to 1889, when he moved to Naples. He
discovered the variation in magnetic susceptibility of a sample subject to mechanic tension (the Villari effect).} ed Arzel\`a mi hanno incaricato di
contraccambiare i Suoi saluti; insieme a questi, accolga ch$^{\underline{\rm mo}}$ Professore, i saluti ed i sensi di affetto
\begin{flushright}
del Suo dev$^{\underline{\rm mo}}$ \\
S. Pincherle
\end{flushright}

La prego di riverire per me il prof. Beltrami

\bigskip

Tentativo di dimostrazione del teorema che ogni funzione algebrica soddisfa ad un [sic!] equaz.[ione] differenziale lineare.

Sia l'equaz[ione].
\begin{equation}\tag{1}
f(x,y)=0
\end{equation}
e sia $\xi$ una funzione razionale di $x$ ed $y$, la quale come \`e noto, si pu\`o scrivere:
\begin{equation}\tag{2}
\xi =a_{0,0}+a_{0,1}y +a_{0,2}y^2+\dots+a_{0,n-1}y^{n-1}
\end{equation}
dove le $a_{h,k}$ sono funzioni razionali di $x$. Si ha dalla (1)
\begin{equation}\tag{3}
\frac{\de f}{\de x}+\frac{\de f}{\de y}y'=0\,,
\end{equation}
onde $y'$ \`e una funzione della forma (2). Ora derivando la (2), viene ponendo per $y$ il suo valore dalla (3):
\begin{equation}\tag{4}
\left\{
\begin{array}{l}
\frac{\dd\xi}{\dd x}=a_{1,0}+a_{1,1}y+\dots a_{1,n-1}y^{n-1} \\
 \\
\frac{\dd^2\xi}{\dd x^2}=a_{2,0}+a_{2,1}y+\dots a_{2,n-1}y^{n-1} \\
 \\
\dots\dots\dots \\
 \\
\frac{\dd^n\xi}{\dd x^n}=a_{n,0}+a_{n,1}y+\dots a_{n,n-1}y^{n-1}:
\end{array}
\right.
\end{equation}
e moltiplicando le (2) e (4) per i reciproci della prima colonna del determinante
$$
\left|
\begin{array}{ccccc}
1 & a_{0,0} & a_{0,1} &\dots & a_{0,n-1} \\
1 & a_{1,0} & a_{0,1} &\dots & a_{1,n-1} \\
\dots & \dots & \dots &\dots & \dots \\
1 & a_{n,0} & a_{n,1} &\dots & a_{n,n-1}
\end{array}
\right|
$$
supposto diverso da zero, viene un'eq.[uazione] della forma
$$
A_0\xi +A_1\frac{\dd\xi}{\dd x}\dots A_n\frac{\dd^n\xi}{\dd x^n}=0.
$$
Nel caso speciale che la $\xi$ sia la stessa $y$, viene l'equaz.[ione] di Tannery.\footnote{If $A_{i,j}$ is the $n\times n$ minor associated with the $(ij)$
element of the matrix (determinante) $A$ written above, Pincherle actually considers
$$
A_{1,1}\xi +A_{2,1}\frac{\dd\xi}{\dd x}+\dots A_{n+1,1}\frac{\dd^n\xi}{\dd x^n}
$$
and, due to the identity ({\it see}, e.g. [Pincherle, 1909], p. 61)
$$
\sum_{j=1}^{n+1}a_{jr}A_{j,s}=0
$$
that holds whenever $r\neq s$ if $a_{jr}$ is an element of $A$, he concludes that
$$
A_0\xi +A_1\frac{\dd\xi}{\dd x}\dots A_n\frac{\dd^n\xi}{\dd x^n}=0:
$$
by setting $\xi=y$, Pincherle thinks that Tannery's equation \eqref{Tannery} is recovered.}

\medskip

{\it Letter 10}: Casorati to Pincherle.

\begin{flushright}
Pavia, 10 marzo 1886
\end{flushright}

Caro Pincherle

............................................

Di matematica non dovrei scrivere nulla, perch\'e non ho la testa a posto. Nondimeno osserver\`o, rispetto alla dimostrazione del Tannery che l'errore di stampa
($m-1$ invece di $m$) della linea 10 di pag. 133 (dove parmi eziandio doversi leggere $y^0$ insieme con $y^2$, $y^3$...,$y^{m-1}$) non costituisca il difetto a
cui io alludevo. Questo difetto si riferisce, se la memoria non m'inganna, alle linee 12 e 13.

Prendiamo il caso di $m=2$.
$$
f(x,y)=ay^2+by+c\,,\qquad \varphi(x)=ac-b^2=(ay^2+2by+c)\times a +(2ay+2b)\left(-\frac a2y-\frac b2\right)
$$
$$
\frac{\dd y}{\dd x}=\frac{\alpha y+\beta}{a\vap(x)}\,,\qquad\frac{\dd^2y}{\dd x^2}=\frac{\gamma y+\delta}{a^2\vap^2}
$$
dove $\alpha$, $\beta$, $\gamma$, $\delta$ sono funzioni intere di $a$, $b$, $c$ e loro derivate $1^{\rm e}$ e $2^{\rm e}$.
L'eliminazione di $y^0$ fra queste espressioni d\`a
$$
\beta\frac{\dd^2y}{\dd x^2}-\frac{\delta}{a\vap}\frac{\dd y}{\dd x}+\frac{\alpha\delta-\beta\gamma}{a^2\vap^2}y=0\,;
$$
ma non \`e cos{\`\i} {\it visibile}, come dice Tannery, che pure supponendo $a$ costante si possa ritenere costante il coeff[iciente]. di $\frac{\dd^2y}{\dd
x^2}$, che Tannery riduce ad 1\footnote{On the back of this minute, Casorati wrote: ``Noto qui, per future notizie che potr\`o dare a Pincherle, che, sulla
formaz.[ione] delle equaz.[ioni] diff.[erenziali] soddisfatte dalle radici di equaz.[ioni] algebriche, havvi da leggere in:

\noindent\underline{Besso} {\it Sull'eq.[uazione] del $5^\circ$ grado}. V.[ol] XIX dele Mem.[orie] lincee (Anno 1883-84), e i successivi lavori da lui presentati
con questo al premio minist.[eriale] scaduto il 30 aprile 1885.

\noindent\underline{Heymann} {\it Ueber die Integration der Diff[erential]gl.[eichungen]} $\frac{\dd^r y}{\dd x^r}+A_m\frac{\dd^m y}{\dd(\ell
x)^{m}}+A_{m-1}\frac{\dd^{m-1} y}{\dd(\ell x)^{m-1}}+\dots+ \frac{\dd y}{\dd(\ell x)}+A_0y=0.$ Math. Annalen XXVI B[and]. 4. H[eft]. Anno 1886.''}

\begin{flushright}
Sono il suo aff.$^{\underline{\rm mo}}$ F.C.
\end{flushright}

\medskip

{\it Letter 11}: Pincherle to Casorati.

\begin{flushright}
Bologna, 15/5/86
\end{flushright}

\begin{center}
Chiarissimo Sig. Professore
\end{center}

\bigskip

Eccomi ancora una volta a ricorrere alla di Lei gentilezza. Nel continuare le mie ricerche sulle {\it operazioni funzionali rappresentabili da integrali
definiti}, ricerche di cui la prima parte \`e in corso di stampa, mi si \`e presentata una famiglia di trasformazioni, che credo nuova e che permette di passare
da qua\-lun\-que equazione lineare differenziale a coefficienti razionali, in una equazione di forma simile alle differenze finite. Prendo la libert\`a di unirle
una breve nota\footnote{Casorati inserted in a footnote the title of the note, by writing: {\it dal titolo ``Sopra una trasformazione delle eq.[uazioni]
diff.[erenziali] lin.[eari] in equaz.[ioni] lin.[eari] alle differenze, e viceversa'' F.C.}} in proposito e Le sarei gratissimo se Ella volesse gettarvi
un'occhiata. Nel caso poi che questo risultato Le sembrasse nuovo e degno d'interesse, ed ove ci\`o non dovesse recarle disturbo, sarei ben lieto se Ella si
compiacesse di presentarla sia all'Istituto Lombardo, sia ai Lincei.

\bigskip

[...]

Nella speranza che ella vorr\`a perdonare il nuovo disturbo che le reco, e pregandola di salutare per me i Sig. professori Beltrami e Bertini, La prego,
ch$^{\rm\underline{mo}}$ sig. professore, di credere ai sensi di profondo rispetto e di vivo affetto del Suo dev$^{\rm\underline{mo}}$

\begin{flushright}
S. Pincherle
\end{flushright}

\bigskip

Possiamo sperare di vedere fra poco un ampliamento ed un seguito delle Sue ricerche sulle funzioni a pi\`u di due periodi?

\medskip

{\it Letter 12}: Pincherle to Casorati.

\begin{flushright}
Bologna, 14/6/86.
\end{flushright}

\begin{flushleft}
Chiarissimo Signor Professore.
\end{flushleft}

La ringrazio vivamente per la premura che Ella si \`e data circa alla mia nota sulle equazioni differenziali, e sono lieto che Le sia sembrata interessante. Ci\`o
che mi ha spinto a redigerla \`e stato il vedere che il Sig. Mellin, nell'ultimo fa\-sci\-co\-lo degli Acta pubblic\`o un teorema che \`e un caso molto, ma molto
particolare della trasformazione che ho indicata.

\noindent....................................

\begin{flushright}
Suo S. Pincherle
\end{flushright}

\begin{center}
{\bf Acknowledgements}
\end{center}

As already remarked in the text, most of the original documents are preserved in the Casorati {\it Nachlass} at Pavia. It is a pleasure to thank once more Prof.
Alberto Gabba for his kindness in disclosing the {\it Nachlass} to me: without him, this paper would not have been written. I would like to thank also prof.
Umberto Bottazzini, prof. Salvatore Cohen, and prof. Otto Liess for useful information on Pincherle's unpublished {\it Ricerche e Saggi} preserved in the
historical section of the library of the Department of Mathematics in Bologna; dr. Anna Biavati and dr. Claudio Cappelletti who helped me during my visit at this
library. I also acknowledge the kind assistance of dr. Alessandra Baretta and dr. Maria Piera Milani during my visits at the historical archive of the University
of Pavia.  Finally, I express my gratitude to prof. Mario Ferrari for his constant encouragement.


\begin{thebibliography}{99}

\bibitem{A37} Amaldi, U., 1937. {Della vita e delle opere di Salvatore Pincherle}. In [Pincherle, 1954], 3-16.

\bibitem{A82a} Appell, P., 1882a. Sur les fonctions uniformes d'un point analytique $(x,y)$. {\it Acta Mathematica}, {\bf 1}, 109-131.

\bibitem{A82b} Appell, P., 1882b. Sur les fonctions uniformes d'un point analytique $(x,y)$. (Seconde m\'emoire). {\it Acta Mathematica}, {\bf 1}, 132-144.

\bibitem{B59} Boole, G., 1859. {\it A treatise on differential equations}, Mc Millan, Cambridge (U.K.).

\bibitem{B80} Boole, G., 1880. {\it Calculus of Finite differences}, Fourth Edition, edited by J.F. Moulton. Chelsea Publishing Company, New York.

\bibitem{B95} Bortolotti, E., 1895. Un contributo alla teoria delle forme lineari alle differenze. {\it Annali di Matematica pura ed applicata}, {\bf 23} (S. II),
    309--344.

\bibitem{B36} Bortolotti, E., 1937. Salvatore Pincherle. {\it Annuario della Regia Universit\`a di Bologna 1936-37}, 151-156.

\bibitem{B94} Bottazzini, U., 1994. {\it Va' Pensiero. Immagini della Matematica nell'Italia dell'Ottocento}. Il Mulino, Bologna.

\bibitem{BG13} Bottazzini, U., Gray, J., 2013. {\it Hidden Harmony--Geometric Fantasies. The rise of Complex Function Theory}. Springer, New York.

\bibitem{B81} Butzer, P.L., 1981. An outline of the life and work of E. B. Christoffel (1829-1900). {\it Historia Mathematica}, {\bf 8}, 243--276.

\bibitem{C68} Casorati, F., 1868. {\it Teorica delle funzioni di variabili complesse.} Volume I. Tip. F.lli Fusi, Pavia.

%\bibitem{C75a} Casorati, F., 1875a. Alcune formole fondamentali per lo studio delle equazioni algebrico-differenziali di primo ordine e secondo grado tra due
variabili ad integrale generale algebrico. {\it Annali di Mat. Pura ed Appl.}, {\bf 7} (S. II), 197-201. In [Casorati, 1952], 3-8.

\bibitem{C75b} Casorati, F., 1875. Sulla teoria delle soluzioni singolari delle equazioni differenziali. {\it Rendiconti del Reale Istituto Lombardo di Lettere,
    Scienze ed Arti}, {\bf 8} (S. 2), 962-966. In [Casorati, 1952], 9--14.

\bibitem{C76} Casorati, F., 1876. Nuova teoria delle soluzioni singolari delle equazioni differenziali di primo ordine e secondo grado tra due variabili. {\it
    Atti della Reale Accademia dei Lincei} {\bf 3} (S. II), 160-167. In [Casorati, 1952], 15-24.

\bibitem{C78_79} Casorati, F., 1878a. Analisi superiore nello [Anno] scolastico 1878-79. Unpublished material in the Casorati {\it Nachlass}, Pavia.

\bibitem{C78} Casorati, F., 1878b. Sulla integrazione delle equazioni algebrico-differenziali di primo ordine e di primo grado per mezzo di funzioni lineari. {\it
    Rendiconti del Reale Istituto Lombardo di Lettere, Scienze ed Arti} {\bf 11} (S. II), 804-808. In [Casorati, 1952], 96-101.

\bibitem{C79_80} Casorati, F., 1879. Analisi superiore nello [Anno] scolastico 1879-80. Unpublished material in the Casorati {\it Nachlass}, Pavia.

\bibitem{C80a} Casorati, F., 1880a. Sull'equazione fondamentale nella teoria delle equazioni differenziali lineari. {\it Rendiconti del Reale Istituto Lombardo di
    Lettere, Scienze ed Arti} (S. II), {\bf 13}, 176-182. In [Casorati, 1952], 102-108.

\bibitem{C80b} Casorati, F., 1880b. Il calcolo delle differenze finite interpretato ed ac\-cre\-sciu\-to di nuovi teoremi a sussidio principalmente delle odierne
    ricerche basate sulla variabilit\`a complessa. {\it Annali di Matematica pura ed applicata}, {\bf 10} (S. II), 10-45. In [Casorati, 1951], 317-356.

\bibitem{C80c} Casorati, F., 1880c. Sopra un recentissimo scritto del sig. L. Stickelberger. {\it Annali di Matematica pura ed applicata}, {\bf 10} (S. II),
    (1880), 154--157. In [Casorati, 1951], 357--361.

\bibitem{C81a} Casorati, F., 1881a. Sur la distinction des int\'egrales des \'equations diff\'erentielles lin\'eaires en sous-groupes. {\it Comptes Rendus
    Hebdomadaires des s\'eances de l'Acad\'emie des Sciences} {\bf 92}, 175--178; 238-241. In [Casorati, 1952], 109--115.

\bibitem{C81b} Casorati, F., 1881b. Sur un \'ecrit tr\`es-r\'ecent de M. Stickelberger. Tip. F.lli Fusi, Pavia.

\bibitem{C86a} Casorati, F., 1886a. Les fonctions d'une seule variable \`a un nombre quelconque de p\'eriodes. {\it Acta Mathematica}, {\bf 8}, 345-359. In
    [Casorati, 1951], 223-238.

\bibitem{C86b} Casorati, F., 1886b. Les lieux fondamentaux des fonctions inverses des int\'egrales ab\'eliennes et en particulier des fonctions inverses des
    int\'egrales elliptiques de $2^{\text{me}}$ et $3^{\text{me}}$ esp\`ece. {\it Acta Mathematica}, {\bf 8}, 360-386. In [Casorati, 1951], 239--265.

\bibitem{C51} Casorati, F. 1951. {\it Opere} (vol. I), Cremonese, Roma, (1951).

\bibitem{C52} Casorati, F. 1952. {\it Opere} (vol. II), Cremonese, Roma, (1952).

\bibitem{C27} Cauchy, A.-L., 1827. Sur l'analogie des puissances et des diff\'erences. In {\it Exercices de Math\'ematique} II ann\'ee. In {\it {\OE}uvres
    Compl\`etes d'Augustin Cauchy}, {\bf 7} (S. II), 198-235, (1889).

\bibitem{C82} Cazzaniga, P., 1882. Il calcolo dei simboli d'operazione elementarmente esposto. {\it Giornale di Matematiche}, {\bf 20}, 48-78; 194-230.

\bibitem{C58} Christoffel, E.B., 1858. Ueber die lineare Abh\"angigkeit von Functionen einer einzigen Ver\"andlichen. {\it Journal f\"ur die reine und angewandte
    Mathematik}, {\bf 55}, 281-299.

\bibitem{C10} Christoffel, E.B., 1910. Zusatz zu der vorstehenden Abhandlung \"uber die lineare Abh\"angigkeit von Funktionen einer einzigen Ver\"anderlichen. In
    {\it Gesammelte mathematische Abhandlungen}, unter Mitwirkung von A. Krazer und G. Faber; herausgegeben von L. Maurer. I Band, Teubner, Leipzig, 106--109.

\bibitem{C62} Cremona, L., 1862. {\it Introduzione ad una teoria geometrica delle curve pia\-ne}, Tip. Gamberini e Parmeggiani, Bologna.

\bibitem{D83a} Daniels, A.L., 1883a. Note on Weierstrass' methods in the theory of elliptic functions. {\it American Journal of Mathematics}, {\bf 6}, 177-182.

\bibitem{D83b} Daniels, A.L., 1883b. Second note on Weierstrass' theory of elliptic functions. {\it American Journal of Mathematics}, {\bf 6}, 253-269.

\bibitem{D84} Daniels, A.L., 1884. Third note on Weierstrass' theory of elliptic functions. {\it American Journal of Mathematics}, {\bf 7}, 82-99.

\bibitem{DHT82} Darboux, G., Hou\"el, J. and Tannery, J., 1882. Revue des publications acad\'emiques et p\'eriodiques. {\it Bulletin des Sciences Math\'ematiques
    et Astronomiques} {\bf 6}$_2$ (S. II), 106-110.

\bibitem{DW82} Dedekind, R., Weber, H., 1882. Theorie der algebraischen Functionen einer Ver\"andlichen. {\it Journal f\"ur die reine und angewandte Mathematik},
    {\bf 92}, 181-290.

\bibitem{D89} Dugac, P., 1989. Henri Poincar\'e. La correspondance avec des math\'ematiciens de J \`a Z. {\it Cahiers du S\'eminaire d'Histoire des
    Math\'ematiques}, {\bf 10}, 83-229.

\bibitem{F83} G. Floquet. 1883. Sue les \'equations diff\'erentielles lin\'eaires \`a coefficents p\'eriodiques. {\it Annales scientifiques de l'\'Ecole Normale Superieure}, {\bf 12} (S. II), 47-88.

\bibitem{F84} G. Floquet. 1884. Sue les \'equations diff\'erentielles lin\'eaires \`a coefficents doublement p\'eriodiques. {\it Annales scientifiques de l'\'Ecole Normale Superieure}, {\bf 1} (S. III), 181-238.

\bibitem{F66} Fuchs, J.L., 1866. Zur Theorie der linearen Differentialgleichungen mit verh\"andlicher Coefficienten. {\it Journal f\"ur die reine und angewandte
    Mathematik}, {\bf 66}, 121-160.

\bibitem{G83} Goursat, E., 1883. Sur une classe des fonctions represent\'ees per des int\'egrales d\'efinies. {\it Acta Mathematica}, {\bf 2}, 1-71.

\bibitem{G86} Gray, J., 1986. Linear differential equations and group theory from Riemann to Poincar\'e. Birkh\"auser, Boston.

\bibitem{G94} Gr\'evy, A., 1894. \'Etude sur les \'equations fonctionnelles. {\it Annales scientifiques de l'\'Ecole Normale Superieure}, {\bf 11} (S. III),
    249-323.

\bibitem{GW11} Guldberg, A. and Wallenberg, G., 1911. {\it Theorie der Linearen Differenzengleichungen}. Teubner, Leipzig und Berlin.

\bibitem{H84} Halphen, G.-H., (1884). Sur la r\'eduction des \'equations diff\'erentielles lin\'eaires aux formes int\'egrables. {\it M\'emoires pr\'esent\'es par
    divers savants \`a l'Acad\'emie des Sciences}, {\bf 28} (S. II), 1884, 1-301. In {\it {\OE}uvres de G.-H. Halphen}, vol. III, Gauthier-Villars, Paris, (1921),
    1-260.

\bibitem{H73} Hamburger, H.L., 1873. Bemerkung \"uber die Form der Integrale der linearen Differentialgleichungen mit ver\"andlichen Coefficienten. {\it Journal
    f\"ur die reine und angewandte Mathematik}, {\bf 76}, 113--125.

\bibitem{H49} Hermite, C., 1849. Sur une question relative \`a la th\'eorie des nombres. {\it Journal de Math\'ematiques pures et appliqu\'ees}, {\bf 14} (S. 1),
    21-30.

\bibitem{H81} Hermite, C., 1881. Sur quelques points de la th\'eorie des fonctions. {\it Journal f\"ur die reine und angewandte Mathematik}, {\bf 91}, 54--78.

\bibitem{K82} Kronecker, L., 1882. Grundz\"uge einer arithmetische Theorie der algebraischen Gr\"ossen. {\it Journal f\"ur die reine und angewandte Mathematik},
    {\bf 92}, 1-122.

\bibitem{M86} Mellin, H., 1886. Zur Theorie der Gammafunctionen. {\it Acta Mathematica} {\bf 8}, 37-80.

\bibitem{M11} Muir, T., 1911. {\it The Theory of Determinants in the historical order of development.} Vol. II: The period 1841-1860. MacMillan, London.

\bibitem{M23} Muir, T., 1923. {\it The Theory of Determinants in the historical order of development.} Vol. IV: The period 1880-1900. MacMillan, London.

\bibitem{N78} Neuenschwander, E., 1978. Der Nachlass von Casorati (1835-1890) in Pavia. {\it Archive for the History of Exact Sciences}, {\bf 19}, 1--89.

\bibitem{N24}  N\"orlund, N.E., 1924. {\it Vorlesungen \"uber Differenzenrechnung}. Springer, Berlin.

\bibitem{P35} Pascal, B., 1835. {\it Les Pens\'ees}. {\it R\'etablies suivant le plan de l'auteur}. Lagier, Dijon.

\bibitem{P97a} Pascal, E., 1897a. {\it Calcolo delle variazioni e calcolo delle differenze finite. (III parte del Calcolo infinitesimale)}. Hoepli, Milan.

\bibitem{P97b} Pascal, E., 1897b. {\it I determinanti. Teoria ed applicazioni con tutte le pi\`u recenti ricerche}. Hoepli, Milan.

\bibitem{Pi79} Picard, \'E., 1879. Sur une g\'en\'eralisation des fonctions p\'eriodiques et sur certaines \'equations diff\'erentielles lin\'eaires. {\it Comptes
    Rendus Hebdomadaires des s\'eances de l'Acad\'emie des Sciences}, {\bf 89}, 140-144.

\bibitem{P81} Picard, \'E., 1881. Sur les \'equations diff\'erentielles lin\'eaires \`a coefficients doublement p\'eriodiques. {\it Journal f\"ur die reine und
    angewandte Mathematik}, {\bf 90}, 281-303.

\bibitem{P74} Pincherle, S., 1874. Sulle superficie di capillarit\`a. {\it Il Nuovo Cimento}, {\bf 12} (S. II), 19-64.

\bibitem{P75} Pincherle, S., 1875. Sulle costanti di capillarit\`a. {\it Il Nuovo Cimento}, {\bf 14} (S. II), 17-25.

\bibitem{P76a} Pincherle, S., 1876a. Sopra alcuni problemi relativi alle superficie d'area minima. {\it Rendiconti del Reale Istituto Lombardo di Lettere, Scienze
    ed Arti}, {\bf 9} (S. 2), 444-456.

\bibitem{P76b} Pincherle, S., 1876b. Sulle superficie d'area minima. {\it Programma del R. Istituto Foscolo di Pavia}, 23 pages, Succ. Bizzoni, Pavia.

\bibitem{P76c} Pincherle, S., 1876c. Nota sulle superficie d'area minima. {\it Giornale di Matematiche}, {\bf 14}, 75-82.

\bibitem{P77} Pincherle, S., 1877. Sulle equazioni algebrico-differenziali di primo ordine e di primo grado a primitiva generale algebrica. {\it Rendiconti del
    Reale Istituto Lombardo di Lettere, Scienze ed Arti}, {\bf 10}, (S. II), 143-151.

\bibitem{P80}  Pincherle, S., 1880. Saggio di una introduzione alla teoria delle funzioni analitiche secondo i principii del Prof. C. Weierstrass. {\it Giornale
    di Ma\-te\-ma\-ti\-che}, {\bf 18}, 178-254; 317-357.

\bibitem{P84a} Pincherle, S., 1884a. Sui sistemi di funzioni analitiche e gli sviluppi in serie formati coi medesimi. {\it Annali di Matematica pura ed
    applicata}, {\bf 12} (S. II), 11-41; 107-133.

\bibitem{P84b} Pincherle, S., 1884b. Sui gruppi lineari di funzioni di una variabile. {\it Memorie della Reale Accademia delle Scienze dell'Istituto di Bologna},
    {\bf 6} (S. IV), 101-118.

\bibitem{P85a} Pincherle, S., 1885a. Alcune osservazioni generali sui gruppi di funzioni. {\it Memorie della Reale Accademia delle Scienze dell'Istituto di
    Bologna}, {\bf 6} (S. IV), 205--214.

\bibitem{P85b} Pincherle, S., 1885b. Note sur une int\'egrale d\'efinie. {\it Acta Mathematica}, {\bf 7}, 381-386.

\bibitem{P86a} Pincherle, S., 1886a. Sopra una trasformazione delle equazioni differenziali lineari in equazioni lineari alle differenze, e viceversa. {\it
    Rendiconti del Reale Istituto Lombardo di Lettere, Scienze ed Arti}, {\bf 19} (S. II), 559-562.

\bibitem{P86b} Pincherle, S., 1886b. Studi sopra alcune operazioni funzionali. {\it Memorie della Reale Accademia delle Scienze dell'Istituto di Bologna}, {\bf
    7}, (S. IV), (1886), 391-442. In [Pincherle, 1954], 92-141.

\bibitem{P87} Pincherle, S., 1887.  Della trasformazione di Laplace e di alcune sue applicazioni. {\it Memorie della Reale Accademia delle Scienze dell'Istituto
    di Bologna}, {\bf 8} (S. IV), 125-143. In [Pincherle, 1954], 173--192.

\bibitem{P90} Pincherle, S., 1890. Saggio di una generalizzazione delle frazioni continue algebriche. {\it Memorie della Reale Accademia delle Scienze
    dell'Istituto di Bologna}, {\bf 10} (S. IV), 513-538.

\bibitem{P92} Pincherle, S., 1892. Sur la g\'en\'eration de syst\`emes r\'ecurrents au moyen d'une \'equation lin\'eaire diff\'erentielle. {\it Acta Mathematica},
    {\bf 16}, 341-363.

\bibitem{P94} Pincherle, S., 1894. Delle funzioni ipergeometriche e di varie questioni ad esse attinenti. {\it Giornale di Matematiche}, {\bf 32}, 209-291. In
    [Pincherle,1954], 273-357.

\bibitem{P95} Pincherle, S., 1895a. L'algebra delle forme lineari alle differenze. {\it Memorie della Reale Accademia delle Scienze dell'Istituto di Bologna},
    {\bf 5} (S. V), 87-126.

\bibitem{P95b} Pincherle, S., 1895b. Sulle operazioni funzionali distributive. {\it Rendiconti della Reale Accademia dei Lincei}, {\bf 4}$_1$, 142-149.

\bibitem{P95c} Pincherle, S., 1895c. Sopra alcune equazioni simboliche. {\it Memorie della Reale Accademia delle Scienze dell'Istituto di Bologna}, {\bf 5} (S.
    V), 663-675.

\bibitem{PA01} Pincherle, S. and Amaldi, U., 1901. {\it Le operazioni distributive e le loro applicazioni all'analisi}. Zanichelli, Bologna.

\bibitem{P09} Pincherle, S., 1909. {\it Lezioni di Algebra Complementare. II: Teoria delle equazioni}. Zanichelli, Bologna.

\bibitem{P26} Pincherle, S., 1926. Il calcolo delle differenze finite. {\it Bollettino dell'Unione Matematica Italiana}, {\bf 5}, 233-242.

\bibitem{P54} Pincherle, S., 1954. {\it Opere Scelte} vol. I, Cremonese, Roma.

\bibitem{Po85} Poincar\'e, H., 1885. Sur les \'equations lin\'eaires aux diff\'erentielles ordinaires et aux diff\'erences finies. {\it American Journal of
    Mathematics}, {\bf 7}, 203--258.

\bibitem{P23} Poincar\'e, H., 1923. Extrait d'une M\'emoire in\'edit de Henri Poincar\'e sur les fonctions Fuchsiennes. {\it Acta Mathematica}, {\bf 39}, 58-93.

\bibitem{P50} Puiseux, V., 1850. Recherches sur les fonctions alg\'ebriques. {\it Journal de Math\'ematiques pures et appliqu\'ees}, {\bf 15} (S. 1), 365-480.

\bibitem{R57} Riemann, B., 1857. Beitr\"age zur Theorie der Gauss'schen Reihe $F(\alpha, \beta,\gamma,x)$ darstellbaren Functionen. {\it Abhandlungen der
    K\"oniglichen Gesellschaft der Wissenschaften zu G\"ottingen}, {\bf 7}. In {\it Bernhard Riemann's Gesammelte mathematische Werke und wissenschaflicher
    Nachlass}. Teubner, Leipzig, 62-78 (1902).

\bibitem{R71} Rosanes, J., 1871. Ueber algebraische Differentialgleichungen erster Ordnung. {\it Mathematische Annalen} {\bf 3}, 535--546.

\bibitem{S81} Stickelberger, L., 1881. {\it Zur Theorie der linearen Differentialgleichungen}. Teubner, Leipzig.

\bibitem{T75} Tannery, J., 1875. Propri\'et\'es des int\'egrales des \'equations diff\'erentielles lin\'eaires \`a coefficients variables. {\it Annales
    Scientifiques de l'\'Ecole Normale Superieure}, {\bf 4} (S. II), 113-182.

\bibitem{T37} Tonelli, L., 1937. Salvatore Pincherle. {\it Annali della Scuola Normale Superiore di Pisa}, {\bf 6} (S. II), 1-10.

\bibitem{V87} Vivanti, G., 1887. Ricerche sulle funzioni uniformid'un punto analitico. {\it Giornale di Matematiche}, {\bf 25}, 54-72; 232-256.

\bibitem{V00} Volterra, V., 1902. Betti, Brioschi, Casorati, trois analystes italiens et trois mani\`eres d'envisager les questions d'analyse. In {\it Compte
    rendu du II Congr\`es international des math\'ematiciens, Paris 1900}, 43-57. In V. Volterra: {\it Opere Ma\-te\-ma\-ti\-che. Memorie e Note} vol. III
    (1900-1913), Roma, Accademia Nazionale dei Lincei, (1957), 1-11.

\bibitem{W76} Weierstrass, K. 1876. Zur Theorie der eindeutigen analytischen Functionen. {\it Abhandlungen der K\"oniglichen Akademie der Wissenschaften zu
    Berlin}, 11-66. In {\it Mathematische Werke von Karl Weierstrass}, Band II, Mayer \& M\"uller, Berlin, (1895), 77-124.

\end{thebibliography}
\end{document}